\documentclass[11pt]{amsart}
\usepackage{amsopn}
\usepackage{amssymb, amscd}
\usepackage{color}
\usepackage [all] {xy}
\usepackage[perpage,symbol*]{footmisc}
\usepackage{hyperref}
\usepackage{array}
\usepackage{multirow}
\usepackage{graphicx}

\newcommand{\nc}{\newcommand}

\nc{\fg}{\mathfrak{f} }     \nc{\vg}{\mathfrak{v} }       \nc{\wg}{\mathfrak{w} }
\nc{\zg}{\mathfrak{z} }     \nc{\ngo}{\mathfrak{n} }      \nc{\kg}{\mathfrak{k} }
\nc{\mg}{\mathfrak{m} }     \nc{\bg}{\mathfrak{b} }       \nc{\ggo}{\mathfrak{g} }
\nc{\sog}{\mathfrak{so} }
\nc{\sug}{\mathfrak{su} }   \nc{\spg}{\mathfrak{sp} }     \nc{\slg}{\mathfrak{sl} }
\nc{\glg}{\mathfrak{gl} }   \nc{\cg}{\mathfrak{c} }       \nc{\rg}{\mathfrak{r} }
\nc{\hg}{\mathfrak{h} }     \nc{\tg}{\mathfrak{t} }       \nc{\ug}{\mathfrak{u} }
\nc{\dg}{\mathfrak{d} }     \nc{\ag}{\mathfrak{a} }       \nc{\pg}{\mathfrak{p} }
\nc{\sg}{\mathfrak{s} }     \nc{\affg}{\mathfrak{aff} }
\nc{\ggob}{\overline{\mathfrak{g}} }

\nc{\pca}{\mathcal{P}}       \nc{\nca}{\mathcal{N}}       \nc{\lca}{\mathcal{L}}
\nc{\oca}{\mathcal{O}}       \nc{\mca}{\mathcal{M}}       \nc{\tca}{\mathcal{T}}
\nc{\aca}{\mathcal{A}}       \nc{\cca}{\mathcal{C}}       \nc{\gca}{\mathcal{G}}
\nc{\sca}{\mathcal{S}}       \nc{\hca}{\mathcal{H}}       \nc{\bca}{\mathcal{B}}
\nc{\dca}{\mathcal{D}}       \nc{\rca}{\mathcal{R}}

\nc{\val}{\operatorname{val}}

\nc{\vp}{\varphi}
\nc{\ddt}{\tfrac{{\rm d}}{{\rm d}t}}
\nc{\dpar}{\tfrac{\partial}{\partial t}}

\nc{\im}{\mathtt{i}}

\nc{\SO}{\mathrm{SO}}           \nc{\Spe}{\mathrm{Sp}}          \nc{\Sl}{\mathrm{SL}}
\nc{\SU}{\mathrm{SU}}           \nc{\Or}{\mathrm{O}}            \nc{\U}{\mathrm{U}}
\nc{\Se}{\mathrm{S}}            \nc{\Cl}{\mathrm{Cl}}           \nc{\Spein}{\mathrm{Spin}}
\nc{\Pin}{\mathrm{Pin}}
\nc{\Glr}{\mathrm{GL}_n(\RR)}   \nc{\Glc}{\mathrm{GL}_n(\CC)}   \nc{\Glv}{\mathrm{GL}(V)}    \nc{\Glk}{\mathrm{GL}_n(\fk)}   \nc{\Gl}{\mathrm{GL}}
\nc{\Grp}{\mathrm{G}}

\nc{\g}{\mathfrak{gl}_n(\RR)}

\nc{\RR}{{\Bbb R}} \nc{\HH}{{\Bbb H}} \nc{\CC}{{\Bbb C}} \nc{\ZZ}{{\Bbb Z}}
\nc{\FF}{{\Bbb F}} \nc{\NN}{{\Bbb N}} \nc{\QQ}{{\Bbb Q}} \nc{\PP}{{\Bbb P}}

\nc{\euler}{{\rm e}}

\nc{\vs}{\vspace{.2cm}} \nc{\vsp}{\vspace{1cm}}

\nc{\ip}{\langle\cdot,\cdot\rangle}
\nc{\ipp}{(\cdot,\cdot)}
\nc{\la}{\langle} \nc{\ra}{\rangle}

\nc{\unm}{\tfrac{1}{2}}\nc{\unc}{\tfrac{1}{4}} \nc{\und}{\tfrac{1}{16}}

\nc{\no}{\vs\noindent}

\nc{\lamn}{\Lambda^2(\RR^n)^*\otimes\RR^n} \nc{\lamp}{\Lambda^2\pg^*\otimes\pg}
\nc{\lamg}{\Lambda^2\ggo^*\otimes\ggo} \nc{\lamngo}{\Lambda^2\ngo^*\otimes\ngo}
\nc{\lamnk}{\Lambda^2(\fk^n)^*\otimes \fk^n} \nc{\lamnkt}{\Lambda^3(\fk^n)^*\otimes \fk^n}

\nc{\tangz}{{\rm T}^{\rm Zar}}
\nc{\mum}{/\!\!/} \nc{\kir}{/\!\!/\!\!/}
\nc{\lievark}{\mathfrak{L}_n(\fk)}         \nc{\lievarc}{\mathfrak{L}_n(\CC)}        \nc{\lievarr}{\mathfrak{L}_n(\RR)}
\nc{\solvvark}{\mathfrak{R}_n(\fk)}        \nc{\solvvarc}{\mathfrak{R}_n(\CC)}       \nc{\solvvarr}{\mathfrak{R}_n(\RR)}
\nc{\nilvark}{\mathfrak{N}_n(\fk)}         \nc{\nilvarc}{\mathfrak{N}_n(\CC)}        \nc{\nilvarr}{\mathfrak{N}_n(\RR)}
\nc{\cirre}{\textrm{C}}
\nc{\fk}{\mathrm{k}}

\nc{\Ri}{\tfrac{4\Ric_{\mu}}{||\mu||^2}}

\nc{\ds}{\displaystyle}

\nc{\ben}{\begin{enumerate}} \nc{\een}{\end{enumerate}}

\nc{\f}{\frac}

\nc{\lb}{[\cdot,\cdot]}

\nc{\isn}{\tfrac{1}{||v||^2}}

\nc{\gkp}{(\ggo=\kg\oplus\pg,\ip)} \nc{\ukh}{(\ug=\kg\oplus\hg,\ip)}

\nc{\Hess}{\operatorname{Hess}}
\nc{\diag}{\operatorname{Diag}}
\nc{\ad}{\operatorname{ad}}       \nc{\Ad}{\operatorname{Ad}}        
\nc{\rank}{\operatorname{rank}}   \nc{\codim}{\operatorname{codim}}  
\nc{\Irr}{\operatorname{Irr}}     \nc{\End}{\operatorname{End}}
\nc{\Aut}{\operatorname{Aut}}     \nc{\Inn}{\operatorname{Inn}}
\nc{\lRad}{\operatorname{Rad}}
\nc{\Der}{\operatorname{Der}}     \nc{\Ker}{\operatorname{Ker}}
\nc{\Iso}{\operatorname{I}}       \nc{\Diff}{\operatorname{Diff}}
\nc{\Lie}{\operatorname{Lie}}     \nc{\tr}{\operatorname{tr}}
\nc{\dif}{\operatorname{d}}
\nc{\sen}{\operatorname{sen}}     \nc{\tang}{\operatorname{T}}
\nc{\modu}{\operatorname{mod}}
\nc{\Riem}{\operatorname{Rm}}     \nc{\Ric}{\operatorname{Ric}}
\nc{\sym}{\operatorname{sym}}     \nc{\symac}{\operatorname{sym^{ac}}}   \nc{\symc}{\operatorname{sym^{c}}}
\nc{\scalar}{\operatorname{sc}}
\nc{\grad}{\operatorname{grad}}
\nc{\ricci}{\operatorname{ric}}   \nc{\nr}{\operatorname{nr}}            \nc{\riccic}{\operatorname{ric^{c}}}
\nc{\riccig}{\operatorname{ric^{\gamma}}}
\nc{\Rin}{\operatorname{M}}
\nc{\Le}{\operatorname{L}}
\nc{\level}{\operatorname{level}} \nc{\rad}{\operatorname{r}}
\nc{\abel}{\operatorname{ab}} \nc{\CH}{\operatorname{CH}}
\nc{\mcc}{\operatorname{mcc}} \nc{\Adj}{\operatorname{Adj}}
\nc{\Order}{\operatorname{O}} \nc{\Ricg}{\operatorname{Ric^{\gamma}}}

\nc{\rhov}{\operatorname{\rho_{v}}}
\nc{\mm}{\operatorname{m}}
\nc{\mmt}{\widetilde{\operatorname{m}}}
\nc{\F}{\operatorname{F}}

\theoremstyle{plain}

\theoremstyle{definition}

\theoremstyle{remark}

\begin{document}
\title{Classification of 7-dimensional Einstein nilradicals \textrm{II}}

\author{EDISON ALBERTO FERN\'ANDEZ CULMA}

\address{Current affiliation: CIEM, FaMAF, Universidad Nacional de C\'ordoba, \newline \indent Ciudad Universitaria, \newline \indent (5000) C\'ordoba, \newline \indent Argentina}

\email{efernandez@famaf.unc.edu.ar}

\thanks{Fully supported by a CONICET fellowship (Argentina)}



\begin{abstract}

This paper contains all computations supporting the classification of $7$-dimensional Einstein nilradicals given in \cite{FERNANDEZ-CULMA1}. Each algebra is analyzed in detail here.

\end{abstract}

\maketitle

\section{Classification of $7$-dimensional Einstein nilradicals}

To get the classification of $7$-dimensional Einstein nilradicals, we follow Carles' classification of low-dimensional
complex nilpotent Lie algebra (with corrections by Magnin in \cite{MAGNIN1}). In \cite{CARLES1}, Carles compares his classification with previous classifications given by Romdhani, Safiullina and Seeley, we recall this work here. The notation that we use is:

\

\begin{center}
\begin{tabular}{cccccccccccc}
\hline
\multicolumn{3}{||m{1.5cm}| }  {\centering CARLES } &
\multicolumn{3}{  m{1.85cm}| }  {\centering \tiny{ROMDHANI} } &
\multicolumn{3}{  m{1.9cm}| }  {\centering {\tiny{SAFIULLINA}} } &
\multicolumn{3}{  m{4.25cm}||} {\centering SELLEY} \tabularnewline
\hline
\hline
\multicolumn{3}{||m{1.5cm}|  }  {\centering $\dim$ $\Der$ } &
\multicolumn{6}{  m{3.75cm}  |  }  {\centering $\dim$ {Derived Series} } &
\multicolumn{3}{  m{4.25cm}|| }  {\centering $\dim$ {Desc. C. Series} } \tabularnewline
\hline
\end{tabular}
\end{center}

\

The differences between the lists in \cite{CARLES1} and \cite{MAGNIN1} are in algebras $1.01(i)$, $1.01(ii)$, $1.8$, $1.9$, $1.12$, $2.2$, $2.10$, $2.45$, $2.46$, $3.1(ii)$ in Carles' list. For instance, in \cite{CARLES1}, $2.45$ is isomorphic to $2.38$ and $2.46$ is isomorphic to $2.36$, so Magnin omits $2.46$ and corrects $2.45$.

In \cite{MAGNIN1}, the only $7$-dimensional algebras of $\rank \geq 1$ that are not written in a nice basis are $1.2(ii)$, $1.2(iv)$, $1.3(i_\lambda)$, $1.3(ii)$, $1.3(v)$, $1.11$, $1.12$, $1.13$, $1.14$, $1.15$, $1.16$, $1.17$, $1.18$, $1.21$, $2.2$, $2.11$, $2.24$, $2.25$, $2.26$, $2.27$, $2.37$. For the remaining algebras, we apply Nikolayevsky's nice basis criterium \cite[Theorem 3.]{NIKOLAYEVSKY1}.

In cases where we give an explicit nilsoliton to prove that certain algebra is an Einstein nilradical, we have also provided a corresponding isomorphism.

Since any nilpotent Lie algebra of dimension $\leq 6$ is an Einstein nilradical and a direct sum of Einstein nilradicals is again an Einstein nilradical, we focus on studying indecomposable algebras.

Let $\mm$ be the unnormalized moment map for the action of $\Glr$ on $\lamn$ by \textquotedblleft{change of basis}\textquotedblright (cf. \cite{LAURET1}). Let $\mu \in V$ be a nilpotent Lie algebra law and we denote by $\Ric_{\mu}$ the Ricci operator of the nilmanifold $(N_{\mu},\ip)$, where $N_{\mu}$ is the simply connected Lie group with $\Lie(N_{\mu})=(\RR^n,\mu)$ and $\ip$ is the canonical inner product of $\RR^n$, then
\begin{equation}\label{momentmapricci}
\mm(\mu)=4\Ric_{\mu}.
\end{equation}
By closeness with GIT, we work with the moment map $\mm$ instead of $\Ric$ so calculations are made using $\mm$.

\section{Rank zero}

Recall that a $\rank$-zero nilpotent Lie algebra (also called characteristically nilpotent Lie algebra) can not be Einstein nilradical since these do not admit an $\NN$-gradation. The following are not Einstein nilradical

\ \newline

\begin{center}
\begin{tabular}{cccccccccccc}
\hline
\multicolumn{3}{||m{1.5cm}| } {\centering $0.1$  } &
\multicolumn{3}{  m{1.5cm}| } {\centering $\mathfrak{n}_{7,2}^{ }$ } &
\multicolumn{3}{  m{1.5cm}| } {\centering  } &
\multicolumn{3}{  m{4.25cm}||} {\centering $(1, 2, 3, 4, 5, 7)F$} \tabularnewline
\hline
\hline
\multicolumn{3}{||m{1.5cm}|  }  {\centering $\dim$ $\Der$ \\      $10$            } &
\multicolumn{6}{  m{3.75cm}  |  }  {\centering $\dim$ {Derived Series}  \\   $(7, 5, 1, 0)$     } &
\multicolumn{3}{  m{4.25cm}|| }  {\centering $\dim$ {Desc. C. Series} \\   $(7, 5, 4, 3, 2, 1, 0)$    } \tabularnewline
\hline
\end{tabular}
\end{center}

\

\noindent $[e_1,e_2]=e_3, [e_1,e_3]=e_4, [e_1,e_4]=e_5, [e_1,e_5]=e_6, [e_1,e_6]=e_7, [e_2, e_3]=e_6, [e_2, e_4]=e_7, [e_2, e_5]=e_7, [e_3, e_4]=-e_7$

\

\begin{center}
\begin{tabular}{cccccccccccc}
\hline
\multicolumn{3}{||m{1.5cm}| } {\centering $0.2$  } &
\multicolumn{3}{  m{1.5cm}| } {\centering $\mathfrak{n}_{7,4}^{ }$ } &
\multicolumn{3}{  m{1.5cm}| } {\centering $L_{75}^{33}$ } &
\multicolumn{3}{  m{4.25cm}||} {\centering $(1, 2, 3, 4, 5, 7)H$ } \tabularnewline
\hline
\hline
\multicolumn{3}{||m{1.5cm}|  }  {\centering $\dim$ $\Der$ \\  10              } &
\multicolumn{6}{  m{3.75cm}  |  }  {\centering $\dim$ {Derived Series}  \\  $(7, 5, 0)$      } &
\multicolumn{3}{  m{4.25cm}|| }  {\centering $\dim$ {Desc. C. Series} \\  $(7, 5, 4, 3, 2, 1, 0)$     } \tabularnewline
\hline
\end{tabular}
\end{center}

\

\noindent $[e_1,e_2]=e_3, [e_1,e_3]=e_4, [e_1,e_4]=e_5, [e_1,e_5]=e_6, [e_1,e_6]=e_7,[e_2, e_3]=e_5 + e_7, [e_2, e_4]=e_6, [e_2, e_5]=e_7$

\

\begin{center}
\begin{tabular}{cccccccccccc}
\hline
\multicolumn{3}{||m{1.5cm}| } {\centering $0.3 $  } &
\multicolumn{3}{  m{1.5cm}| } {\centering $\mathfrak{n}_{7,7 }^{ }$ } &
\multicolumn{3}{  m{1.5cm}| } {\centering $L_{75 }^{34 }$ } &
\multicolumn{3}{  m{4.25cm}||} {\centering $(1, 2, 3, 4, 5, 7)E$ } \tabularnewline
\hline
\hline
\multicolumn{3}{||m{1.5cm}|  }  {\centering $\dim$ $\Der$ \\              $11$    } &
\multicolumn{6}{  m{3.75cm}  |  }  {\centering $\dim$ {Derived Series}  \\   $(7, 5, 0)$     } &
\multicolumn{3}{  m{4.25cm}|| }  {\centering $\dim$ {Desc. C. Series} \\   $(7, 5, 4, 3, 2, 1, 0)$     } \tabularnewline
\hline
\end{tabular}
\end{center}

\

\noindent $[e_1,e_2]=e_3, [e_1,e_3]=e_4, [e_1,e_4]=e_5, [e_1, e_5]=e_6, [e_1,e_6]=e_7, [e_2, e_3]=e_6 + e_7, [e_2, e_4]=e_7$

\

\begin{center}
\begin{tabular}{cccccccccccc}
\hline
\multicolumn{3}{||m{1.5cm}| } {\centering $0.4_{\lambda}$  } &
\multicolumn{3}{  m{1.5cm}| } {\centering $\mathfrak{n}_{7,14}^{\lambda }$ } &
\multicolumn{3}{  m{1.5cm}| } {\centering $L_{ 74}^{25 }( \mathfrak{n}_3 )$ \\ $(\alpha=\lambda)$ } &
\multicolumn{3}{  m{4.25cm}||} {\centering $(1, 2, 4, 5, 7)N$ \\ $(\xi=\frac{\lambda+2\im}{\lambda-2\im} \neq 1$, $\lambda \neq \pm 2\im$) \\
                                          $(1, 2, 4, 5, 7)M$ \\ $(\xi=0$, $\lambda=-2\im$)} \tabularnewline
\hline
\hline
\multicolumn{3}{||m{1.5cm}|  }  {\centering $\dim$ $\Der$ \\        $10$          } &
\multicolumn{6}{  m{3.75cm}  |  }  {\centering $\dim$ {Derived Series}  \\   $(7, 5, 1, 0)$     } &
\multicolumn{3}{  m{4.25cm}|| }  {\centering $\dim$ {Desc. C. Series} \\     $(7, 5, 4, 2, 1, 0)$   } \tabularnewline
\hline
\end{tabular}
\end{center}

\

\noindent $[e_1, e_2]=e_3, [e_1, e_3]=e_4, [e_1, e_4]=\lambda e_7 + e_6, [e_1, e_5]=e_7, [e_1, e_6]=e_7, [e_2, e_3]=e_5, [e_2, e_4]=e_7, [e_2, e_5]=e_6, [e_3, e_5]=e_7$

\

\begin{center}
\begin{tabular}{cccccccccccc}
\hline
\multicolumn{3}{||m{1.5cm}| } {\centering $0.5$  } &
\multicolumn{3}{  m{1.5cm}| } {\centering $\mathfrak{n}_{7,16}^{ }$ } &
\multicolumn{3}{  m{1.5cm}| } {\centering $L_{74}^{17}(\mathfrak{n}_3)$ } &
\multicolumn{3}{  m{4.25cm}||} {\centering $(1, 2, 4, 5, 7)N$ \\ $\xi=1$} \tabularnewline
\hline
\hline
\multicolumn{3}{||m{1.5cm}|  }  {\centering $\dim$ $\Der$ \\          $10$         } &
\multicolumn{6}{  m{3.75cm}  |  }  {\centering $\dim$ {Derived Series}  \\   $(7, 5, 1, 0)$      } &
\multicolumn{3}{  m{4.25cm}|| }  {\centering $\dim$ {Desc. C. Series} \\  $(7, 5, 4, 2, 1, 0)$      } \tabularnewline
\hline
\end{tabular}
\end{center}

\

\noindent $[e_1,e_2]=e_3, [e_1,e_3]=e_4, [e_1,e_4]=e_6 + e_7, [e_1,e_6]=e_7, [e_2,e_3]=e_5, [e_2, e_5]=e_6, [e_3, e_5]=e_7$

\

\begin{center}
\begin{tabular}{cccccccccccc}
\hline
\multicolumn{3}{||m{1.5cm}|  }  {\centering  $0.6$  } &
\multicolumn{3}{  m{1.5cm}|  }  {\centering  $\mathfrak{n}_{7,20}^{}$ } &
\multicolumn{3}{  m{1.5cm}|  }  {\centering  $L_{74}^{26}(\mathfrak{n}_3)$ \\ $(\alpha\neq0)$} &
\multicolumn{3}{  m{4.25cm}||}  {\centering  $(1, 2, 4, 5, 7)J$  } \tabularnewline
\hline
\hline
\multicolumn{3}{||m{1.5cm}|   }  {\centering $\dim$ $\Der$ \\               $10$    } &
\multicolumn{6}{  m{3.75cm}  |   }  {\centering $\dim$ {Derived Series}  \\     $(7, 5, 1, 0)$    } &
\multicolumn{3}{  m{4.25cm}|| }  {\centering $\dim$ {Desc. C. Series} \\     $(7, 5, 4, 2, 1, 0)$    } \tabularnewline
\hline
\end{tabular}
\end{center}

\

\noindent $[e_1,e_2]=e_3, [e_1,e_3]=e_4, [e_1,e_4]=e_7, [e_1,e_5]=e_6,[e_1,e_6]=e_7,[e_2, e_3]=e_5, [e_2, e_4]=e_6, [e_2, e_5]=e_7, [e_3, e_4]=e_7$

\

\begin{center}
\begin{tabular}{cccccccccccc}
\hline
\multicolumn{3}{||m{1.5cm}|  }  {\centering  $0.7$  } &
\multicolumn{3}{  m{1.5cm}|  }  {\centering  $\mathfrak{n}_{7,25}^{}$ } &
\multicolumn{3}{  m{1.5cm}|  }  {\centering   } &
\multicolumn{3}{  m{4.25cm}||}  {\centering  $(1, 3, 4, 5, 7)I$  } \tabularnewline
\hline
\hline
\multicolumn{3}{||m{1.5cm}|   }  {\centering $\dim$ $\Der$ \\               $10$    } &
\multicolumn{6}{  m{3.75cm}  |   }  {\centering $\dim$ {Derived Series}  \\     $(7, 5, 1, 0)$    } &
\multicolumn{3}{  m{4.25cm}|| }  {\centering $\dim$ {Desc. C. Series} \\     $(7, 5, 4, 2, 1, 0)$    } \tabularnewline
\hline
\end{tabular}
\end{center}

\

\noindent  $[e_1,e_2]=e_3, [e_1,e_3]=e_4, [e_1,e_4]=e_7, [e_1,e_5]=e_7,[e_1,e_6]=e_7,[e_2, e_3]=e_5, [e_2, e_4]=e_7, [e_2, e_5]=e_6, [e_3, e_5]=e_7$

\

\begin{center}
\begin{tabular}{cccccccccccc}
\hline
\multicolumn{3}{||m{1.5cm}|  }  {\centering  $0.8$  } &
\multicolumn{3}{  m{1.5cm}|  }  {\centering  $\mathfrak{n}_{7,34}^{}$ } &
\multicolumn{3}{  m{1.5cm}|  }  {\centering  $L_{74}^{22}(\mathfrak{n}_3)$ } &
\multicolumn{3}{  m{4.25cm}||}  {\centering  $(1, 2, 4, 5, 7)G$  } \tabularnewline
\hline
\hline
\multicolumn{3}{||m{1.5cm}|   }  {\centering $\dim$ $\Der$ \\               $10$    } &
\multicolumn{6}{  m{3.75cm}  |   }  {\centering $\dim$ {Derived Series}  \\     $(7, 4, 1, 0)$    } &
\multicolumn{3}{  m{4.25cm}|| }  {\centering $\dim$ {Desc. C. Series} \\     $(7, 4, 3, 2, 1, 0)$    } \tabularnewline
\hline
\end{tabular}
\end{center}

\

\noindent  $[e_1,e_2]=e_4, [e_1,e_3]=e_7, [e_1,e_4]=e_5, [e_1,e_5]=e_6, [e_2, e_3]=e_6, [e_2, e_4]=e_6, [e_2, e_6]=e_7, [e_4, e_5]=-e_7$

\
\newpage

\section{Rank one}

\

\begin{center}
\begin{tabular}{cccccccccccc}
\hline
\multicolumn{3}{||m{1.5cm}|  }  {\centering  $1.01(i)$  } &
\multicolumn{3}{  m{1.5cm}|  }  {\centering  $\mathfrak{n}_{7,43}^{}$ } &
\multicolumn{3}{  m{1.5cm}|  }  {\centering  $L_{75}^{12}$ } &
\multicolumn{3}{  m{4.25cm}||}  {\centering  $(1, 2, 3, 5, 7)B$  } \tabularnewline
\hline
\hline
\multicolumn{3}{||m{1.5cm}|   }  {\centering $\dim$ $\Der$ \\               $11$    } &
\multicolumn{6}{  m{3.75cm}  |   }  {\centering $\dim$ {Derived Series}  \\     $(7, 4, 0)$    } &
\multicolumn{3}{  m{4.25cm}|| }  {\centering $\dim$ {Desc. C. Series} \\     $(7, 4, 3, 2, 1, 0)$    } \tabularnewline
\hline
\end{tabular}
\end{center}

\

\noindent  $[e_1,e_2]=e_4,[e_1,e_4]=e_5,[e_1,e_5]=e_6,[e_1,e_6]=e_7,[e_2,e_3]=e_5+e_7,[e_3,e_4]=-e_6,[e_3,e_5]=-e_7$

\

\noindent Pre-Einstein Derivation: $\phi = \diag(0, 1, 0, 1, 1, 1, 1)$. It is not an Einstein nilradical by \cite[Theorem 2.5]{FERNANDEZ-CULMA1}; $\phi \ngtr 0$.

\

\begin{center}
\begin{tabular}{cccccccccccc}
\hline
\multicolumn{3}{||m{1.5cm}|  }  {\centering  $1.01(ii)$  } &
\multicolumn{3}{  m{1.5cm}|  }  {\centering  $\mathfrak{n}_{7,47}^{}$ } &
\multicolumn{3}{  m{1.5cm}|  }  {\centering  $L_{75}^{13}$ } &
\multicolumn{3}{  m{4.25cm}||}  {\centering  $(1, 2, 4, 5, 7)B$  } \tabularnewline
\hline
\hline
\multicolumn{3}{||m{1.5cm}|   }  {\centering $\dim$ $\Der$ \\               $12$    } &
\multicolumn{6}{  m{3.75cm}  |   }  {\centering $\dim$ {Derived Series}  \\     $(7, 4, 0)$    } &
\multicolumn{3}{  m{4.25cm}|| }  {\centering $\dim$ {Desc. C. Series} \\     $(7, 4, 3, 2, 1, 0)$    } \tabularnewline
\hline
\end{tabular}
\end{center}

\

\noindent  $[e_1,e_2]=e_4,[e_1,e_4]=e_5,[e_1,e_5]=e_6,[e_1,e_6]=e_7,[e_2,e_3]=e_6+e_7, [e_3,e_4]=-e_7$

\

\noindent Pre-Einstein Derivation: $\phi = \diag(0, 1, 0, 1, 1, 1, 1)$. It is not an Einstein nilradical by \cite[Theorem 2.5]{FERNANDEZ-CULMA1}; $\phi \ngtr 0$.

\

\begin{center}
\begin{tabular}{cccccccccccc}
\hline
\multicolumn{3}{||m{1.5cm}|  }  {\centering  $1.02$  } &
\multicolumn{3}{  m{1.5cm}|  }  {\centering  $\mathfrak{n}_{7,23}^{}$ } &
\multicolumn{3}{  m{1.5cm}|  }  {\centering   } &
\multicolumn{3}{  m{4.25cm}||}  {\centering  $(1, 2, 4, 5, 7)K$  } \tabularnewline
\hline
\hline
\multicolumn{3}{||m{1.5cm}|   }  {\centering $\dim$ $\Der$ \\               $11$    } &
\multicolumn{6}{  m{3.75cm}  |   }  {\centering $\dim$ {Derived Series}  \\     $(7, 5, 1, 0)$    } &
\multicolumn{3}{  m{4.25cm}|| }  {\centering $\dim$ {Desc. C. Series} \\     $(7, 5, 4, 2, 1, 0)$    } \tabularnewline
\hline
\end{tabular}
\end{center}

\

\noindent  $[e_1,e_2]=e_3,[e_1,e_3]=e_4,[e_1,e_5]=e_6,[e_2,e_3]=e_5,[e_2,e_4]=e_6,[e_2,e_5]=e_7,[e_2,e_6]=e_7,[e_3,e_5]=-e_7$

\

\noindent Pre-Einstein Derivation: $\phi =\frac{1}{2}\diag(0, 1, 1, 1, 2, 2, 3)$. It is not an Einstein nilradical by \cite[Theorem 2.5]{FERNANDEZ-CULMA1}; $\phi \ngtr 0$.

\

\begin{center}
\begin{tabular}{cccccccccccc}
\hline
\multicolumn{3}{||m{1.5cm}|  }  {\centering  $1.03$  } &
\multicolumn{3}{  m{1.5cm}|  }  {\centering  $\mathfrak{n}_{7,26}^{}$ } &
\multicolumn{3}{  m{1.5cm}|  }  {\centering  $L_{74}^{13}(\mathfrak{n}_3)$ } &
\multicolumn{3}{  m{4.25cm}||}  {\centering  $(1, 3, 4, 5, 7)G$  } \tabularnewline
\hline
\hline
\multicolumn{3}{||m{1.5cm}|   }  {\centering $\dim$ $\Der$ \\               $11$    } &
\multicolumn{6}{  m{3.75cm}  |   }  {\centering $\dim$ {Derived Series}  \\     $(7, 5, 1, 0)$    } &
\multicolumn{3}{  m{4.25cm}|| }  {\centering $\dim$ {Desc. C. Series} \\     $(7, 5, 4, 2, 1, 0)$    } \tabularnewline
\hline
\end{tabular}
\end{center}

\

\noindent  $[e_1,e_2]=e_3,[e_1,e_3]=e_4,[e_1,e_4]=e_5,[e_1,e_6]=e_7,[e_2,e_3]=e_6,[e_2,e_4]=e_7,[e_2,e_5]=e_7,[e_3,e_4]=-e_7$

\

\noindent Pre-Einstein Derivation: $\phi =\frac{2}{3}\diag(0, 1, 1, 1, 1, 2, 2)$. It is not an Einstein nilradical by \cite[Theorem 2.5]{FERNANDEZ-CULMA1}; $\phi \ngtr 0$.

\

\begin{center}
\begin{tabular}{cccccccccccc}
\hline
\multicolumn{3}{||m{1.5cm}|  }  {\centering  $1.1(i_\lambda)$ \\ $\lambda\neq 0,1$} &
\multicolumn{3}{  m{1.5cm}|  }  {\centering  $\mathfrak{n}_{7,1}^{\lambda}$ \\ $\lambda\neq 0,1$ } &
\multicolumn{3}{  m{1.5cm}|  }  {\centering  } &
\multicolumn{3}{  m{4.25cm}||}  {\centering  $(1,2,3,4,5,7)I$ \\ $(\xi=\lambda,$ $\xi\neq0)$ } \tabularnewline
\hline
\hline
\multicolumn{3}{||m{1.5cm}|   }  {\centering $\dim$ $\Der$ \\               $10$    } &
\multicolumn{6}{  m{3.75cm}  |   }  {\centering $\dim$ {Derived Series}  \\     $(7, 5, 1, 0)$    } &
\multicolumn{3}{  m{4.25cm}|| }  {\centering $\dim$ {Desc. C. Series} \\     $(7, 5, 4, 3, 2, 1, 0)$    } \tabularnewline
\hline
\end{tabular}
\end{center}

\

\noindent  $[e_1,e_2]=e_3,[e_1,e_3]=e_4,[e_1,e_4]=e_5,[e_1,e_5]=e_6,[e_1,e_6]=e_7,[e_2,e_3]=e_5,[e_2,e_4]=e_6,[e_2,e_5]=\lambda e_7,[e_3,e_4]=(1-\lambda)e_7$

\

\noindent  $\{e_1...e_7\}$ is a nice basis.

$$U=\tiny{\left[ \begin {array}{ccccccccc} 3&0&1&1&1&0&1&1&-1 \\ \noalign{\medskip}0&3&0&1&1&1&-1&0&0\\ \noalign{\medskip}1&0&3&0&1&1&1&-1&1\\ \noalign{\medskip}1&1&0&3&0&-1&1&1&0\\ \noalign{\medskip}1&1&1&0&3&0&-1&1&1\\ \noalign{\medskip}0&1&1&-1&0&3&1&0&1\\ \noalign{\medskip}1&-1&1&1&-1&1&3&1&1\\ \noalign{\medskip}1&0&-1&1&1&0&1&3&1\\ \noalign{\medskip}-1&0&1&0&1&1&1&1&3\end {array} \right]}$$

\noindent A solution to $Ux=[1]$: $x=\frac{1}{10}(1,2,2,2,1,2,1,2,1)^T$. It is an Einstein nilradical.

\

\noindent Pre-Einstein Derivation: $\frac{1}{5}\diag(1, 2, 3, 4, 5, 6, 7)$ \newline
\noindent $||\mathcal{S}_{\beta}||^2=\frac{5}{7}\thickapprox 0.7142857143$

\

\begin{center}
\begin{tabular}{cccccccccccc}
\hline
\multicolumn{3}{||m{1.5cm}|  }  {\centering  $1.1(i_\lambda)$ \\ $\lambda = 0$} &
\multicolumn{3}{  m{1.5cm}|  }  {\centering  $\mathfrak{n}_{7,1}^{\lambda}$ \\ $\lambda = 0$ } &
\multicolumn{3}{  m{1.5cm}|  }  {\centering  } &
\multicolumn{3}{  m{4.25cm}||}  {\centering  $(1,2,3,4,5,7)I$ \\ $(\xi = 0)$ } \tabularnewline
\hline
\hline
\multicolumn{3}{||m{1.5cm}|   }  {\centering $\dim$ $\Der$ \\               $10$    } &
\multicolumn{6}{  m{3.75cm}  |   }  {\centering $\dim$ {Derived Series}  \\     $(7, 5, 1, 0)$    } &
\multicolumn{3}{  m{4.25cm}|| }  {\centering $\dim$ {Desc. C. Series} \\     $(7, 5, 4, 3, 2, 1, 0)$    } \tabularnewline
\hline
\end{tabular}
\end{center}

\

\noindent $[e_1,e_2]=e_3,[e_1,e_3]=e_4,[e_1,e_4]=e_5,[e_1,e_5]=e_6,[e_1,e_6]=e_7,[e_2,e_3]=e_5,[e_2,e_4]=e_6,[e_3,e_4]=e_7$

\

\noindent $\{e_1...e_7\}$ is a nice basis.

$$U= \tiny{\left[ \begin {array}{cccccccc} 3&0&1&1&1&0&1&-1\\ \noalign{\medskip}0
&3&0&1&1&1&-1&0\\ \noalign{\medskip}1&0&3&0&1&1&1&1
\\ \noalign{\medskip}1&1&0&3&0&-1&1&0\\ \noalign{\medskip}1&1&1&0&3&0&
-1&1\\ \noalign{\medskip}0&1&1&-1&0&3&1&1\\ \noalign{\medskip}1&-1&1&1
&-1&1&3&1\\ \noalign{\medskip}-1&0&1&0&1&1&1&3\end {array} \right]}$$

\noindent General solution to $Ux=[1]$: $x=(t_2,-\frac{1}{5}+t_1+t_2,0,\frac{3}{5}-t_1-t_2,\frac{2}{5}-t_2,\frac{3}{5}-t_1-t_2,t_1,t_2)^T$. It is not an Einstein nilradical.

\

\noindent Pre-Einstein Derivation: $\frac{1}{5}\diag(1, 2, 3, 4, 5, 6, 7)$

\

\begin{center}
\begin{tabular}{cccccccccccc}
\hline
\multicolumn{3}{||m{1.5cm}|  }  {\centering  $1.1(i_\lambda)$ \\ $\lambda = 1$} &
\multicolumn{3}{  m{1.5cm}|  }  {\centering  $\mathfrak{n}_{7,1}^{\lambda}$ \\ $\lambda = 1$ } &
\multicolumn{3}{  m{1.5cm}|  }  {\centering  $L_{75}^{30}$} &
\multicolumn{3}{  m{4.25cm}||}  {\centering  $(1,2,3,4,5,7)G$ } \tabularnewline
\hline
\hline
\multicolumn{3}{||m{1.5cm}|   }  {\centering $\dim$ $\Der$ \\               $11$    } &
\multicolumn{6}{  m{3.75cm}  |   }  {\centering $\dim$ {Derived Series}  \\     $(7, 5, 0)$    } &
\multicolumn{3}{  m{4.25cm}|| }  {\centering $\dim$ {Desc. C. Series} \\     $(7, 5, 4, 3, 2, 1, 0)$    } \tabularnewline
\hline
\end{tabular}
\end{center}

\

\noindent $[e_1,e_2]=e_3,[e_1,e_3]=e_4,[e_1,e_4]=e_5,[e_1,e_5]=e_6,[e_1,e_6]=e_7,[e_2,e_3]=e_5,[e_2,e_4]=e_6,[e_2,e_5]=e_7$

\

\noindent $\{e_1...e_7\}$ is a nice basis.

$$U=\tiny{\left[ \begin {array}{cccccccc} 3&0&1&1&1&0&1&1\\ \noalign{\medskip}0
&3&0&1&1&1&-1&0\\ \noalign{\medskip}1&0&3&0&1&1&1&-1
\\ \noalign{\medskip}1&1&0&3&0&-1&1&1\\ \noalign{\medskip}1&1&1&0&3&0&
-1&1\\ \noalign{\medskip}0&1&1&-1&0&3&1&0\\ \noalign{\medskip}1&-1&1&1
&-1&1&3&1\\ \noalign{\medskip}1&0&-1&1&1&0&1&3\end {array} \right]}$$

\noindent General solution to $Ux=[1]$: $x=(0,\frac{2}{5}-t_1,\frac{3}{5}-t_1-t_2,t_1,-\frac{1}{5}+t_1+t_2,t_1,t_2,\frac{3}{5}-t_1-t_2)^T$. It is not an Einstein nilradical.

\

\noindent Pre-Einstein Derivation: $\frac{1}{5}\diag(1, 2, 3, 4, 5, 6, 7)$

\
\begin{center}
\begin{tabular}{cccccccccccc}
\hline
\multicolumn{3}{||m{1.5cm}|  }  {\centering  $1.1(ii)$  } &
\multicolumn{3}{  m{1.5cm}|  }  {\centering  $\mathfrak{n}_{7,3}^{}$ } &
\multicolumn{3}{  m{1.5cm}|  }  {\centering   } &
\multicolumn{3}{  m{4.25cm}||}  {\centering  $(1,2,3,4,5,7)C$  } \tabularnewline
\hline
\hline
\multicolumn{3}{||m{1.5cm}|   }  {\centering $\dim$ $\Der$ \\               $11$    } &
\multicolumn{6}{  m{3.75cm}  |   }  {\centering $\dim$ {Derived Series}  \\     $(7, 5, 1, 0)$    } &
\multicolumn{3}{  m{4.25cm}|| }  {\centering $\dim$ {Desc. C. Series} \\     $(7, 5, 4, 3, 2, 1, 0)$    } \tabularnewline
\hline
\end{tabular}
\end{center}

\

\noindent  $[e_1,e_2]=e_3,[e_1,e_3]=e_4,[e_1,e_4]=e_5,[e_1,e_5]=e_6,[e_1,e_6]=e_7,[e_2,e_5]=e_7,[e_3,e_4]=-e_7$

\

\noindent  $\{e_1...e_7\}$ is a nice basis.

$$U=\tiny{\left[ \begin {array}{ccccccc} 3&0&1&1&1&1&-1\\ \noalign{\medskip}0&3
&0&1&1&0&0\\ \noalign{\medskip}1&0&3&0&1&-1&1\\ \noalign{\medskip}1&1&0
&3&0&1&0\\ \noalign{\medskip}1&1&1&0&3&1&1\\ \noalign{\medskip}1&0&-1&
1&1&3&1\\ \noalign{\medskip}-1&0&1&0&1&1&3\end {array} \right]}$$

\noindent  General solution to $Ux=[1]$: $x=(t_1,\frac{2}{5},\frac{3}{5}-t_1,0,-\frac{1}{5},\frac{3}{5}-t_1,t_1)^T$. It is not an Einstein nilradical.

\

\noindent  Pre-Einstein Derivation: $\frac{1}{5}\diag(1, 2, 3, 4, 5, 6, 7)$

\

\begin{center}
\begin{tabular}{cccccccccccc}
\hline
\multicolumn{3}{||m{1.5cm}|  }  {\centering  $1.1(iii)$  } &
\multicolumn{3}{  m{1.5cm}|  }  {\centering  $\mathfrak{n}_{7,11}^{}$ } &
\multicolumn{3}{  m{1.5cm}|  }  {\centering  $L_{74}^{41}(\mathfrak{n}_3)$ } &
\multicolumn{3}{  m{4.25cm}||}  {\centering  $(2,3,4,5,7)G$  } \tabularnewline
\hline
\hline
\multicolumn{3}{||m{1.5cm}|   }  {\centering $\dim$ $\Der$ \\               $10$    } &
\multicolumn{6}{  m{3.75cm}  |   }  {\centering $\dim$ {Derived Series}  \\     $(7, 5, 1, 0)$    } &
\multicolumn{3}{  m{4.25cm}|| }  {\centering $\dim$ {Desc. C. Series} \\     $(7, 5, 4, 3, 2, 0)$    } \tabularnewline
\hline
\end{tabular}
\end{center}

\

\noindent  $[e_1,e_2]=e_3,[e_1,e_3]=e_4,[e_1,e_4]=e_5,[e_1,e_5]=e_6,[e_2,e_3]=e_5,[e_2,e_4]=e_6,[e_2,e_5]=-e_7,[e_3,e_4]=e_7$

\

\noindent  $\{e_1...e_7\}$ is a nice basis.

$$U=\tiny{ \left[ \begin {array}{cccccccc} 3&0&1&1&0&1&1&-1\\ \noalign{\medskip}0
&3&0&1&1&-1&0&0\\ \noalign{\medskip}1&0&3&0&1&1&-1&1
\\ \noalign{\medskip}1&1&0&3&-1&1&1&0\\ \noalign{\medskip}0&1&1&-1&3&1
&0&1\\ \noalign{\medskip}1&-1&1&1&1&3&1&1\\ \noalign{\medskip}1&0&-1&1
&0&1&3&1\\ \noalign{\medskip}-1&0&1&0&1&1&1&3\end {array} \right]}$$

\noindent  A solution to $Ux=[1]$: $x=\frac{1}{10}(2,3,2,1,1,1,2,2)^T$. It is an Einstein nilradical.

\

\noindent  Pre-Einstein Derivation: $\frac{1}{5}\diag(1, 2, 3, 4, 5, 6, 7)$ \newline
\noindent $||\mathcal{S}_{\beta}||^2=\frac{5}{7}\thickapprox 0.714$

\

\begin{center}
\begin{tabular}{cccccccccccc}
\hline
\multicolumn{3}{||m{1.5cm}|  }  {\centering  $1.1(iv)$  } &
\multicolumn{3}{  m{1.5cm}|  }  {\centering  $\mathfrak{n}_{7,22}^{}$ } &
\multicolumn{3}{  m{1.5cm}|  }  {\centering  $L_{74}^{27}(\mathfrak{n}_3)$ } &
\multicolumn{3}{  m{4.25cm}||}  {\centering  $(1,2,4,5,7)I$  } \tabularnewline
\hline
\hline
\multicolumn{3}{||m{1.5cm}|   }  {\centering $\dim$ $\Der$ \\               $11$    } &
\multicolumn{6}{  m{3.75cm}  |   }  {\centering $\dim$ {Derived Series}  \\     $(7, 5, 1, 0)$    } &
\multicolumn{3}{  m{4.25cm}|| }  {\centering $\dim$ {Desc. C. Series} \\     $(7, 5, 4, 2, 1, 0)$    } \tabularnewline
\hline
\end{tabular}
\end{center}

\

\noindent  $[e_1,e_2]=e_3,[e_1,e_3]=e_4,[e_1,e_5]=e_6,[e_1,e_6]=e_7,[e_2,e_3]=e_5,[e_2,e_4]=e_6,[e_2,e_5]=e_7,[e_3,e_4]=e_7$

\

\noindent  $\{e_1...e_7\}$ is a nice basis.

$$U=\tiny{\left[ \begin {array}{cccccccc} 3&0&1&1&0&1&1&-1\\ \noalign{\medskip}0
&3&1&1&1&-1&0&0\\ \noalign{\medskip}1&1&3&0&-1&1&1&0
\\ \noalign{\medskip}1&1&0&3&0&-1&1&1\\ \noalign{\medskip}0&1&-1&0&3&1
&0&1\\ \noalign{\medskip}1&-1&1&-1&1&3&1&1\\ \noalign{\medskip}1&0&1&1
&0&1&3&1\\ \noalign{\medskip}-1&0&0&1&1&1&1&3\end {array} \right]}$$

\noindent  General solution to $Ux=[1]$: $x=(\frac{2}{5}-t_1,\frac{1}{5}-t_1+t_2,\frac{1}{5}-t_2+t_1,t_1,\frac{1}{5}-t_2+t_1,t_2,0,\frac{2}{5}-t_1 )^T$. It is not an Einstein nilradical.

\

\noindent  Pre-Einstein Derivation: $\frac{1}{5}\diag(1, 2, 3, 4, 5, 6, 7)$

\

\begin{center}
\begin{tabular}{cccccccccccc}
\hline
\multicolumn{3}{||m{1.5cm}|  }  {\centering  $1.1(v)$  } &
\multicolumn{3}{  m{1.5cm}|  }  {\centering  $\mathfrak{n}_{7,42}^{}$ } &
\multicolumn{3}{  m{1.5cm}|  }  {\centering  $L_{74}^{2}(3\CC)$ } &
\multicolumn{3}{  m{4.25cm}||}  {\centering  $(1,2,3,5,7)C$  } \tabularnewline
\hline
\hline
\multicolumn{3}{||m{1.5cm}|   }  {\centering $\dim$ $\Der$ \\               $10$    } &
\multicolumn{6}{  m{3.75cm}  |   }  {\centering $\dim$ {Derived Series}  \\     $(7, 4, 0)$    } &
\multicolumn{3}{  m{4.25cm}|| }  {\centering $\dim$ {Desc. C. Series} \\     $(7, 4, 3, 2, 1, 0)$    } \tabularnewline
\hline
\end{tabular}
\end{center}

\

\noindent  $[e_1,e_3]=e_4,[e_1,e_4]=e_5,[e_1,e_5]=e_6,[e_1,e_6]=e_7,[e_2,e_3]=e_5,[e_2,e_4]=e_6,[e_2,e_5]=e_7,[e_3,e_4]=-e_7$

\

\noindent  $\{e_1...e_7\}$ is a nice basis.

$$U=\tiny{ \left[ \begin {array}{cccccccc} 3&0&1&1&1&-1&0&0\\ \noalign{\medskip}0
&3&0&1&1&1&-1&1\\ \noalign{\medskip}1&0&3&0&-1&1&1&0
\\ \noalign{\medskip}1&1&0&3&0&-1&1&1\\ \noalign{\medskip}1&1&-1&0&3&1
&0&1\\ \noalign{\medskip}-1&1&1&-1&1&3&1&1\\ \noalign{\medskip}0&-1&1&
1&0&1&3&1\\ \noalign{\medskip}0&1&0&1&1&1&1&3\end {array} \right]}$$

\noindent General solution to $Ux=[1]$: $x=(\frac{2}{5}-t_1,\frac{3}{5}-t_1-t_2,t_1,-\frac{1}{5}+t_1+t_2,t_1,t_2,\frac{3}{5}-t_1-t_2,0)^T$. It is not an Einstein nilradical.

\

\noindent Pre-Einstein Derivation: $\frac{1}{5}\diag(1, 2, 3, 4, 5, 6, 7)$

\

\begin{center}
\begin{tabular}{cccccccccccc}
\hline
\multicolumn{3}{||m{1.5cm}|  }  {\centering  $1.1(vi)$  } &
\multicolumn{3}{  m{1.5cm}|  }  {\centering  $\mathfrak{n}_{7,38}^{}$ } &
\multicolumn{3}{  m{1.5cm}|  }  {\centering  $L_{74}^{7}(\mathfrak{n}_3)$ } &
\multicolumn{3}{  m{4.25cm}||}  {\centering  $(1,3,4,5,7)E$  } \tabularnewline
\hline
\hline
\multicolumn{3}{||m{1.5cm}|   }  {\centering $\dim$ $\Der$ \\               $11$    } &
\multicolumn{6}{  m{3.75cm}  |   }  {\centering $\dim$ {Derived Series}  \\     $(7, 4, 1, 0)$    } &
\multicolumn{3}{  m{4.25cm}|| }  {\centering $\dim$ {Desc. C. Series} \\     $(7, 4, 3, 2, 1, 0)$    } \tabularnewline
\hline
\end{tabular}
\end{center}

\

\noindent  $[e_1,e_2]=e_3,[e_1,e_3]=e_4,[e_1,e_4]=e_5,[e_1,e_6]=e_7,[e_2,e_3]=e_5,[e_2,e_5]=e_7,[e_3,e_4]=-e_7$

\

\noindent  $\{e_1...e_7\}$ is a nice basis.

$$U=\tiny{ \left[ \begin {array}{ccccccc} 3&0&1&1&0&1&-1\\ \noalign{\medskip}0&3
&0&1&1&0&0\\ \noalign{\medskip}1&0&3&1&1&-1&1\\ \noalign{\medskip}1&1&
1&3&0&1&1\\ \noalign{\medskip}0&1&1&0&3&0&1\\ \noalign{\medskip}1&0&-1
&1&0&3&1\\ \noalign{\medskip}-1&0&1&1&1&1&3\end {array} \right]}$$

\noindent  General solution to $Ux=[1]$: $x=(t_1,\frac{2}{5},\frac{3}{5}-t_1,-\frac{1}{5},0,\frac{3}{5}-t_1,t_1 )^T$. It is not an Einstein nilradical.

\

\noindent  Pre-Einstein Derivation: $\frac{1}{5}\diag(1, 2, 3, 4, 5, 6, 7)$

\

\begin{center}
\begin{tabular}{cccccccccccc}
\hline
\multicolumn{3}{||m{1.5cm}|  }  {\centering  $1.2(i_\lambda)$ \\ $\lambda \neq 0 , 1$ } &
\multicolumn{3}{  m{1.5cm}|  }  {\centering  $\mathfrak{n}_{7,52}^{\lambda}$ \\ \cite{CARLES1}} &
\multicolumn{3}{  m{1.5cm}|  }  {\centering  $L_{74}^{8}(3\CC)$ \\ \cite{CARLES1}} &
\multicolumn{3}{  m{4.25cm}||}  {\centering  $(1,3,5,7)S$ \\ \cite{CARLES1}  } \tabularnewline
\hline
\hline
\multicolumn{3}{||m{1.5cm}|   }  {\centering $\dim$ $\Der$ \\               $12$    } &
\multicolumn{6}{  m{3.75cm}  |   }  {\centering $\dim$ {Derived Series}  \\     $(7, 4, 0)$    } &
\multicolumn{3}{  m{4.25cm}|| }  {\centering $\dim$ {Desc. C. Series} \\     $(7, 4, 3, 1, 0)$    } \tabularnewline
\hline
\end{tabular}
\end{center}

\

\noindent  $[e_1,e_2]=e_4,[e_1,e_3]=e_6,[e_1,e_4]=e_5,[e_1,e_5]=e_7,[e_2,e_3]=\lambda e_5,[e_2,e_4]=e_6,[e_2,e_6]=e_7,[e_3,e_4]=(1-\lambda)e_7$

\

\noindent  $\{e_1...e_7\}$ is a nice basis.

$$U=\tiny{\left[ \begin {array}{cccccccc} 3&1&0&1&1&0&1&-1\\ \noalign{\medskip}
1&3&1&1&1&1&-1&1\\ \noalign{\medskip}0&1&3&0&1&1&0&1
\\ \noalign{\medskip}1&1&0&3&-1&0&1&1\\ \noalign{\medskip}1&1&1&-1&3&1
&1&1\\ \noalign{\medskip}0&1&1&0&1&3&0&1\\ \noalign{\medskip}1&-1&0&1&
1&0&3&1\\ \noalign{\medskip}-1&1&1&1&1&1&1&3\end {array} \right]}$$

\noindent A solution to $Ux=[1]$: $x=\frac{1}{44}(7,6,8,7,3,8,10,3)^T$. It is an Einstein nilradical.

\

\noindent Pre-Einstein Derivation: $\frac{4}{11}\diag(1, 1, 2, 2, 3, 3, 4)$ \newline
\noindent $||\mathcal{S}_{\beta}||^2=\frac{11}{13}\thickapprox 0.846$
\

\begin{center}
\begin{tabular}{cccccccccccc}
\hline
\multicolumn{3}{||m{1.5cm}|  }  {\centering  $1.2(i_\lambda)$ \\ $\lambda = 0$ } &
\multicolumn{3}{  m{1.5cm}|  }  {\centering  $\mathfrak{n}_{7,52}^{\lambda}$ \\ \cite{CARLES1}} &
\multicolumn{3}{  m{1.5cm}|  }  {\centering  $L_{74}^{8}(3\CC)$ \\ \cite{CARLES1}} &
\multicolumn{3}{  m{4.25cm}||}  {\centering  $(1,3,5,7)S$ \\ \cite{CARLES1}  } \tabularnewline
\hline
\hline
\multicolumn{3}{||m{1.5cm}|   }  {\centering $\dim$ $\Der$ \\               $12$    } &
\multicolumn{6}{  m{3.75cm}  |   }  {\centering $\dim$ {Derived Series}  \\     $(7, 4, 0)$    } &
\multicolumn{3}{  m{4.25cm}|| }  {\centering $\dim$ {Desc. C. Series} \\     $(7, 4, 3, 1, 0)$    } \tabularnewline
\hline
\end{tabular}
\end{center}

\

\noindent  $[e_1,e_2]=e_4,[e_1,e_3]=e_6,[e_1,e_4]=e_5,[e_1,e_5]=e_7,[e_2,e_4]=e_6,[e_2,e_6]=e_7,[e_3,e_4]=e_7$

\

\noindent  $\{e_1...e_7\}$ is a nice basis.

$$U=\tiny{ \left[ \begin {array}{ccccccc} 3&1&0&1&0&1&-1\\ \noalign{\medskip}1&3
&1&1&1&-1&1\\ \noalign{\medskip}0&1&3&0&1&0&1\\ \noalign{\medskip}1&1&0
&3&0&1&1\\ \noalign{\medskip}0&1&1&0&3&0&1\\ \noalign{\medskip}1&-1&0&
1&0&3&1\\ \noalign{\medskip}-1&1&1&1&1&1&3\end {array} \right]}$$

\noindent A solution to $Ux=[1]$: $x=\frac{1}{22}(5,3,4,2,4,5,3 )^T$. It is an Einstein nilradical.

\

\noindent Pre-Einstein Derivation: $\frac{4}{11}\diag(1, 1, 2, 2, 3, 3, 4)$ \newline
\noindent $||\mathcal{S}_{\beta}||^2=\frac{11}{13}\thickapprox 0.846$

\

\begin{center}
\begin{tabular}{cccccccccccc}
\hline
\multicolumn{3}{||m{1.5cm}|  }  {\centering  $1.2(i_\lambda)$ \\ $\lambda = 1$ } &
\multicolumn{3}{  m{1.5cm}|  }  {\centering  $\mathfrak{n}_{7,52}^{\lambda}$ \\ \cite{CARLES1}} &
\multicolumn{3}{  m{1.5cm}|  }  {\centering  $L_{74}^{8}(3\CC)$ \\ \cite{CARLES1}} &
\multicolumn{3}{  m{4.25cm}||}  {\centering  $(1,3,5,7)S$ \\ \cite{CARLES1}  } \tabularnewline
\hline
\hline
\multicolumn{3}{||m{1.5cm}|   }  {\centering $\dim$ $\Der$ \\               $12$    } &
\multicolumn{6}{  m{3.75cm}  |   }  {\centering $\dim$ {Derived Series}  \\     $(7, 4, 0)$    } &
\multicolumn{3}{  m{4.25cm}|| }  {\centering $\dim$ {Desc. C. Series} \\     $(7, 4, 3, 1, 0)$    } \tabularnewline
\hline
\end{tabular}
\end{center}

\

\noindent  $[e_1,e_2]=e_4,[e_1,e_3]=e_6,[e_1,e_4]=e_5,[e_1,e_5]=e_7,[e_2,e_3]=e_5,[e_2,e_4]=e_6,[e_2,e_6]=e_7$

\

\noindent  $\{e_1...e_7\}$ is a nice basis.

$$U=\tiny{ \left[ \begin {array}{ccccccc} 3&1&0&1&1&0&1\\ \noalign{\medskip}1&3&
1&1&1&1&-1\\ \noalign{\medskip}0&1&3&0&1&1&0\\ \noalign{\medskip}1&1&0
&3&-1&0&1\\ \noalign{\medskip}1&1&1&-1&3&1&1\\ \noalign{\medskip}0&1&1
&0&1&3&0\\ \noalign{\medskip}1&-1&0&1&1&0&3\end {array} \right]}$$

\noindent  A solution to $Ux=[1]$: $x=\frac{1}{22}(2,3,4,5,3,4,5)^T$. It is an Einstein nilradical.

\

\noindent  Pre-Einstein Derivation: $\frac{4}{11}\diag(1, 1, 2, 2, 3, 3, 4)$. \newline
\noindent $||\mathcal{S}_{\beta}||^2=\frac{11}{13}\thickapprox 0.846$

\

\begin{center}
\begin{tabular}{cccccccccccc}
\hline
\multicolumn{3}{||m{1.5cm}|  }  {\centering  $1.2(ii)$  } &
\multicolumn{3}{  m{1.5cm}|  }  {\centering  } &
\multicolumn{3}{  m{1.5cm}|  }  {\centering  $L_{74}^{11}(\mathfrak{n}_3)$ } &
\multicolumn{3}{  m{4.25cm}||}  {\centering  $(1,3,5,7)S$  \\ $\xi=0$} \tabularnewline
\hline
\hline
\multicolumn{3}{||m{1.5cm}|   }  {\centering $\dim$ $\Der$ \\               $12$    } &
\multicolumn{6}{  m{3.75cm}  |   }  {\centering $\dim$ {Derived Series}  \\     $(7, 4, 0)$    } &
\multicolumn{3}{  m{4.25cm}|| }  {\centering $\dim$ {Desc. C. Series} \\     $(7, 4, 3, 1, 0)$    } \tabularnewline
\hline
\end{tabular}
\end{center}

\

\noindent  $\mu:=\left\{ \begin{array}{l}
[e_1,e_2]=e_4, [e_1,e_4]=e_5, [e_1,e_5]=e_7, [e_1,e_6]=e_7, [e_2,e_3]=e_6, \\
{[e_2,e_4]=e_6}, [e_2,e_5]=e_7, [e_3,e_4]=-e_7\end{array}\right.$

\

\noindent It is not an Einstein nilradical by \cite[Theorem 2.6]{FERNANDEZ-CULMA1}. In fact,

\

\noindent  Einstein Pre-Derivation $\phi=\frac{4}{11}\diag(1, 1, 2, 2, 3, 3, 4)$. \newline
Let $X=\diag(1, -1, 0, 0, 1, -1, 0) $. $X\in \ggo_{\phi}$ therefore $g_{t}=\exp(tX) \in \Grp_{\phi}$.

\

\noindent  $g_t\cdot\mu=\left\{\begin{array}{l}
[e_1,e_2]=e_4, [e_1,e_4]=e_5, [e_1,e_5]=\euler^{-2t}e_7, [e_1,e_6]=e_7, [e_2,e_3]=e_6, \\
{[e_2,e_4]=e_6}, [e_2,e_5]=e_7, [e_3,e_4]=-e_7 \end{array}
\right.$

\

\noindent The $\Grp_{\phi}-$orbit of $\mu$ is not closed because  $g_t\cdot\mu \rightarrow \widetilde{\mu}$ as $t \rightarrow \infty$ and $(\RR^7,\widetilde{\mu})$ is non-isomorphic to $(\RR^7,\mu)$; $(\RR^7,\mu)$ has $\rank=1$ and $(\RR^7,\widetilde{\mu})$ has $\rank=2$ with maximal torus generated by $D1=\diag(1,0,1,1,2,1,2)$ and $D2=\diag(0,1,1,1,1,2,2)$

\

\begin{center}
\begin{tabular}{cccccccccccc}
\hline
\multicolumn{3}{||m{1.5cm}|  }  {\centering  $1.2(iii)$  } &
\multicolumn{3}{  m{1.5cm}|  }  {\centering  $\mathfrak{n}_{7,54}^{}$ } &
\multicolumn{3}{  m{1.5cm}|  }  {\centering  $L_{74}^{39}(3\CC)$ } &
\multicolumn{3}{  m{4.25cm}||}  {\centering  $(2,3,5,7)D$  } \tabularnewline
\hline
\hline
\multicolumn{3}{||m{1.5cm}|   }  {\centering $\dim$ $\Der$ \\               $12$    } &
\multicolumn{6}{  m{3.75cm}  |   }  {\centering $\dim$ {Derived Series}  \\     $(7, 4, 0)$    } &
\multicolumn{3}{  m{4.25cm}|| }  {\centering $\dim$ {Desc. C. Series} \\     $(7, 4, 3, 1, 0)$    } \tabularnewline
\hline
\end{tabular}
\end{center}

\

\noindent  $[e_1,e_2]=e_4,[e_1,e_3]=e_6,[e_1,e_4]=e_5,[e_1,e_5]=e_7,[e_2,e_3]=e_5,[e_2,e_4]=e_6,[e_3,e_4]=-e_7$

\

\noindent  $\{e_1...e_7\}$ is a nice basis.

$$U=\tiny{\left[ \begin {array}{ccccccc} 3&1&0&1&1&0&-1\\ \noalign{\medskip}1&3
&1&1&1&1&1\\ \noalign{\medskip}0&1&3&0&1&1&1\\ \noalign{\medskip}1&1&0
&3&-1&0&1\\ \noalign{\medskip}1&1&1&-1&3&1&1\\ \noalign{\medskip}0&1&1
&0&1&3&1\\ \noalign{\medskip}-1&1&1&1&1&1&3\end {array} \right]}$$

\noindent  General solution to $Ux=[1]$: $x=(\frac{1}{11}+t_1,-\frac{1}{11},\frac{2}{11},\frac{5}{11}-t_1,\frac{4}{11}-t_1,\frac{2}{11},t_1)^T$. It is not an Einstein nilradical.

\

\noindent  Pre-Einstein Derivation: $\frac{4}{11}\diag(1, 1, 2, 2, 3, 3, 4)$

\

\begin{center}
\begin{tabular}{cccccccccccc}
\hline
\multicolumn{3}{||m{1.5cm}|  }  {\centering  $1.2(iv)$  } &
\multicolumn{3}{  m{1.5cm}|  }  {\centering  $\mathfrak{n}_{7,72}^{}$ } &
\multicolumn{3}{  m{1.5cm}|  }  {\centering  $L_{74}^{16}(3\CC)$ } &
\multicolumn{3}{  m{4.25cm}||}  {\centering  $(1,3,5,7)H$  } \tabularnewline
\hline
\hline
\multicolumn{3}{||m{1.5cm}|   }  {\centering $\dim$ $\Der$ \\               $12 $    } &
\multicolumn{6}{  m{3.75cm}  |   }  {\centering $\dim$ {Derived Series}  \\     $(7, 4, 0)$    } &
\multicolumn{3}{  m{4.25cm}|| }  {\centering $\dim$ {Desc. C. Series} \\     $(7, 4, 2, 1, 0)$    } \tabularnewline
\hline
\end{tabular}
\end{center}

\

\noindent $\mu:=\left\{\begin{array}{l}
[e_1,e_2]=e_4, [e_1,e_4]=e_6,  [e_1,e_5]=-e_7, [e_1,e_6]=e_7, [e_2,e_3]=e_5,\\
{[e_2,e_5]=e_7},  [e_3,e_4]=e_7
\end{array}
\right.$

\

\noindent It is not an Einstein nilradical by \cite[Theorem 2.6]{FERNANDEZ-CULMA1}. In fact,

\

\noindent  Einstein Pre-Derivation $\phi=\frac{4}{11}\diag (1, 1, 2, 2, 3, 3, 4)$. \newline
Let $X=\diag (3, -10, 18, -8, 7, -6,-4 ) $. $X\in \ggo_{\phi}$ therefore $g_{t}=\exp(tX) \in \Grp_{\phi}$.

\noindent $g_t\cdot\mu=\left\{\begin{array}{l}
[e_1,e_2]=\euler^{-t}e_4, [e_1,e_4]=\euler^{-t}e_6, [e_1,e_5]=-\euler^{-14t}e_7, [e_1,e_6]=\euler^{-t}e_7,\\
{[e_2,e_3]}=\euler^{-t}e_5, [e_2,e_5]=\euler^{-t}e_7, [e_3,e_4]=\euler^{-14t}e_7
\end{array} \right.$

\

\noindent The $\Grp_{\phi}-$orbit of $\mu$ is not closed because  $g_t\cdot\mu \rightarrow 0$ as $t \rightarrow \infty$

\

\begin{center}
\begin{tabular}{cccccccccccc}
\hline
\multicolumn{3}{||m{1.5cm}|  }  {\centering  $1.3(i_\lambda)$ \\ $\lambda \neq 0$ } &
\multicolumn{3}{  m{1.5cm}|  }  {\centering  $\mathfrak{n}_{7,62}^{\lambda}$ } &
\multicolumn{3}{  m{1.5cm}|  }  {\centering  $L_{74}^{15}(3\CC)$ \\ $(\alpha=\lambda)$ } &
\multicolumn{3}{  m{4.25cm}||}  {\centering  $(1,3,5,7)N$ \\ $\xi = \lambda$} \tabularnewline
\hline
\hline
\multicolumn{3}{||m{1.5cm}|   }  {\centering $\dim$ $\Der$ \\               $13$    } &
\multicolumn{6}{  m{3.75cm}  |   }  {\centering $\dim$ {Derived Series}  \\     $(7, 4, 0)$    } &
\multicolumn{3}{  m{4.25cm}|| }  {\centering $\dim$ {Desc. C. Series} \\     $(7, 4, 2, 1, 0)$    } \tabularnewline
\hline
\end{tabular}
\end{center}

\

\noindent $\mu:= \left\{ \begin{array}{l}[e_1,e_2]=e_4,[e_1,e_3]=e_5,[e_1,e_4]=e_6,[e_1,e_6]=e_7,[e_2,e_3]=e_6,\\
{[e_2,e_4]}=\lambda e_7,[e_2,e_5]=e_7,[e_3,e_5]=e_7 \end{array} \right.$

\

For any $\lambda \neq 0$, $\ggo_{1.3(i_{\lambda})}$ is an Einstein nilradical. We prove this by contradiction; assume that $\ggo_{1.3(i_{\lambda})}$ is not an Einstein nilradical. The derivation $\phi$ given by the diagonal matrix $\frac{5 }{17 }\diag( 1, 2, 2, 3, 3, 4, 5 )$ with respect to the basis $\{e_i\}$ is pre-Einstein. It follows from \cite[Theorem 2.6]{FERNANDEZ-CULMA1} that the orbit $\mathrm{G}_{\phi}\cdot\mu$ is not closed and so there exist a symmetric matrix $A\in \ggo_{\phi}$ such that $\mu$ degenerates by the action of the one-parameter subgroup $\exp(tA)$ as $t\rightarrow \infty$. As $A\in \ggo_{\phi}$ and $A$ is a symmetric matrix, then there exist $X=\diag(a_1,...,a_7) \in \ggo_{\phi}$ and $A(\alpha)$, $B(\beta)$ in $\SO_2(\RR)$ such that
$$
A=\diag(1,A(-\alpha),B(-\beta),1,1) X \diag(1,A(\alpha),B(\beta),1,1).
$$
As the action is continuous, it follows that $\mu$ also degenerates by the action of
$$
g_t:=\exp(tX) \diag(1,A(\alpha),B(\beta),1,1)
$$
as $t\rightarrow \infty$. The contradiction will be found in this last fact. The action of $g_t$ in $\mu$ is
$$
g_t\cdot\mu =\left\{\begin{array}{l}
{[e_1,e_2]}={\euler^{-t \left( a_1+a_2-a_4 \right) }}\cos \left( \alpha-\beta \right) e_4-{\euler^{-t \left( a_1+a_2-a_5 \right) }}\sin \left( \alpha-\beta \right) e_5,\\
{[e_1,e_3]}={\euler^{-t \left( a_1+a_3-a_4 \right) }}\sin \left( \alpha-\beta \right) e_4+{\euler^{-t \left( a_1+a_3-a_5 \right) }}\cos \left( \alpha-\beta \right) e_5,\\
{[e_1,e_4]}={\euler^{-t \left( a_1+a_4-a_6 \right) }}\cos \left( \beta \right) e_6,[e_1,e_5]={\euler^{-t\left( a_1+a_5-a_6 \right) }}\sin \left( \beta \right) e_6,\\
{[e_1,e_6]}={\euler^{-t \left( a_1+a_6-a_7 \right) }}e_7,[e_2,e_3]={\euler^{-t \left( a_2+a_3-a_6 \right) }}e_6,\\
{[e_2,e_4]}={\euler^{-t \left( a_2+a_4-a_7 \right) }} f_{2,4}(\alpha,\beta)e_7,\\
{[e_2,e_5]}={\euler^{-t \left( a_2+a_5-a_7 \right) }} f_{2,5}(\alpha,\beta)e_7,\\
{[e_3,e_4]}={\euler^{-t \left( a_3+a_4-a_7 \right) }} f_{3,4}(\alpha,\beta)e_7,\\
{[e_3,e_5]}={\euler^{-t \left( a_3+a_5-a_7 \right) }} f_{3,5}(\alpha,\beta)e_7,
\end{array}
\right.
$$
where
$$
\begin{array}{l}
f_{2,4}(\alpha, \beta)=(\cos \left( \alpha \right) \cos \left( \beta \right) \lambda + \sin \left( \alpha \right) \sin \left( \beta \right) -\cos \left( \alpha \right) \sin \left( \beta \right) )\\
f_{2,5}(\alpha, \beta)=(\cos \left( \alpha \right) \sin \left( \beta \right) \lambda - \sin \left( \alpha \right) \cos \left( \beta \right) +\cos \left( \alpha \right) \cos \left( \beta \right) )\\
f_{3,4}(\alpha, \beta)=(\sin \left( \alpha \right) \cos \left( \beta \right) \lambda - \cos \left( \alpha \right) \sin \left( \beta \right) -\sin \left( \alpha \right) \sin \left( \beta \right) )\\
f_{3,5}(\alpha, \beta)=(\sin \left( \alpha \right) \sin \left( \beta \right) \lambda + \cos \left( \alpha \right) \cos \left( \beta \right) +\sin \left( \alpha \right) \cos \left( \beta \right) ).
\end{array}
$$
Depending on the values of $\alpha$ and $\beta$, some terms are zero in the Lie algebra law $g_t\cdot \mu$ and since the degeneration is determinated for non-zero terms, our attention should be at the exponents of the exponential factor of such terms; when $t>0$, such exponent must be non negative. We consider
$$
\begin{array}{ll}
\multicolumn{2}{l}{p_1 := a_1+a_2+a_3+a_4+a_5+a_6+a_7,} \tabularnewline
\multicolumn{2}{l}{p_2 := a_1+2a_2+2a_3+3a_4+3a_5+4a_6+5a_7,} \tabularnewline
q_1 := a_1+a_2-a_4-b_1^2,       & q_2 := a_1+a_2-a_5-b_2^2,\\
q_3 := a_1+a_3-a_4-b_3^2,       & q_4 := a_1+a_3-a_5-b_4^2,\\
q_5 := a_1+a_4-a_6-b_5^2,       & q_6 := a_1+a_5-a_6-b_6^2,\\
q_7 := a_1+a_6-a_7-b_7^2,       & q_8 := a_2+a_3-a_6-b_8^2,\\
q_9 := a_2+a_4-a_7-b_9^2,       & q_{10}:= a_2+a_5-a_7-b_{10}^2,\\
q_{11} := a_3+a_4-a_7-b_{11}^2, & q_{12}:= a_3+a_5-a_7-b_{12}^2,
\end{array}
$$
where $p_1$ and $p_2$ correspond to that $X\in \ggo_{\phi}$ and $q_1,...,q_{12}$ correspond to the observation above.

It is easy to see that pairs of functions $\{f_{2,4},f_{2,5}\}$, $\{f_{2,4},f_{3,4}\}$, $\{f_{2,5},f_{3,5}\}$ and $\{f_{3,4},f_{3,5}\}$ do not vanish simultaneously (as $\sin$ and $\cos$). We have the following cases depending on which terms are non zero.\newline

\noindent I) $\cos(\beta)$ and $\sin(\beta)$ are non zero

\

\noindent 1. $\cos(\alpha-\beta),\,f_{2,4},\,f_{3,4}$ are non zero \newline

In this case, the degeneration gives a non-trivial solution to the system of polynomial equations $\{p_1=0, p_2=0, q_7=0, q_8=0, q_5=0, q_6=0\} \cup \{q_1=0,q_4=0\} \cup \{q_9=0\} \cup \{q_{11}=0\}$.

\

By solving this system, we have $b_1 = \sqrt{b_8^2-b_{11}^2+b_7^2}$,  \newline $b_4 = \frac{1}{2}\sqrt{-12b_7^2-10b_8^2-4b_9^2-8b_6^2}$, $b_5 = \frac{1}{2}\sqrt{-16b_7^2-14b_8^2-4b_6^2}$ and  $b_6, b_7, b_8,$ $b_9, b_{11} \in \RR$ .

Therefore $b_6 = 0$, $b_7 = 0$, $b_8 = 0$, $b_9 = 0$ and it follows that $b_4 = 0$, $b_5 = 0$ and $b_1=\pm \im b_{11}$. So it must be that $b_{11}=0$ and then $b_1=0$. There is no non-trivial degeneration.

\

\noindent 2. $\cos(\alpha-\beta),\,f_{2,4},\,f_{3,5}$ are non zero \newline

In this case, the degeneration gives a non-trivial solution to the system of polynomial equations $\{p_1=0, p_2=0, q_7=0, q_8=0, q_5=0, q_6=0\} \cup \{q_1=0,q_4=0\} \cup \{q_9=0\} \cup \{q_{12=0}\}$.

\

By solving this system, we have $b_1 = \frac{1}{2}\sqrt{20b_7^2+18b_8^2+8b_6^2-4b_{12}^2}$, $b_4 = \frac{1}{2}\sqrt{-8b_6^2-12b_7^2-10b_8^2-4b_9^2}$, $b_5 = \frac{1}{2}\sqrt{-4b_6^2-16b_7^2-14b_8^2}$ and $b_6, b_7, b_8, b_9,$ $b_{12} \in \RR$.

Therefore $b_6 = 0$, $b_7 = 0$, $b_8 = 0$, $b_9 = 0$ and it follows that $b_4 = 0$, $b_5 = 0$ and $b_1=\pm \im b_{12}$. So it must be that $b_{12}=0$ and then $b_1=0$. There is no non-trivial degeneration.

\

\noindent 3. $\cos(\alpha-\beta),\,f_{2,5},\,f_{3,4}$ are non zero \newline

In this case, the degeneration gives a non-trivial solution to the system of polynomial equations $\{p_1=0, p_2=0, q_7=0, q_8=0, q_5=0, q_6=0\} \cup \{q_1=0,q_4=0\} \cup \{q_{10}=0\} \cup \{q_{11}=0\}$.

\

By solving this system, we have $b_1 = \sqrt{-b_{11}^2+b_{7}^2+b_8^2}$, $b_{10} = \sqrt{-b_4^2+b_7^2+b_8^2}$,  $b_5 = \frac{1}{2}\sqrt{-16b_7^2-14b_8^2-4b_6^2}$ and
$b_6, b_7, b_8, b_{11}, b_4 \in \RR$.

Therefore $b_6 = 0$, $b_7 = 0$, $b_8 = 0$ and it follows that $b_5 = 0$, $b_1=\pm \im b_{11}$, $b_{10}=\pm \im b_{4}$. So it must be that $b_{11}=0$, $b_{4}=0$ and then $b_1=0$ and $b_{10}=0$. There is no non-trivial degeneration.

\

\noindent 4. $\cos(\alpha-\beta),\,f_{2,5},\,f_{3,5}$ are non zero \newline

Since $f_{2,4}$ and $f_{3,4}$ do not vanish simultaneously, this case is included in the case I) 2. and I) 3. There is no non-trivial degeneration.

\

\noindent 5. $\sin(\alpha-\beta),\,f_{2,4},\,f_{3,4}$ are non zero \newline

In this case, the degeneration gives a non-trivial solution to the system of polynomial equations $\{p_1=0, p_2=0, q_7=0, q_8=0, q_5=0, q_6=0\} \cup \{q_2=0,q_3=0\} \cup \{q_9=0\} \cup \{q_{11}=0\}$.

\

By solving this system, we have $b_{11} = \frac{1}{2}\sqrt{-12b_7^2-10b_8^2-4b_2^2-8b_6^2}$,  $b_3 = \sqrt{b_8^2-b_9^2+b_7^2}$, $b_5 = \frac{1}{2}\sqrt{-16b_7^2-14b_8^2-4b_6^2}$ and $b_6, b_7, b_8, b_9, b_2 \in \RR$.

Therefore $b_7 = 0$, $b_8 = 0$, $b_2 = 0$ $b_6 = 0$ and it follows that $b_{11}=0$, $b_3=\pm \im b_9$, $b_5=0$. So it must be that $b_9=0$, and then $b_3=0$. There is no non-trivial degeneration.

\

\noindent 6. $\sin(\alpha-\beta),\,f_{2,4},\,f_{3,5}$ are non zero \newline

In this case, the degeneration gives a non-trivial solution to the system of polynomial equations $\{p_1=0, p_2=0, q_7=0, q_8=0, q_5=0, q_6=0\} \cup \{q_2=0,q_3=0\} \cup \{q_9=0\} \cup \{q_{12}=0\}$.

\

By solving this system, we have $b_{12} =\sqrt{b_7^2+b_8^2-b_2^2}$, $b_3 = \sqrt{b_7^2+b_8^2-b_9^2}$, $b_5 = \frac{1}{2}\sqrt{-16b_7^2-14b_8^2-4b_6^2}$ and $b_2, b_6, b_7, b_8, b_9 \in \RR$.

Therefore $b_6 = 0$, $b_7 = 0$, $b_8 = 0$ and it follows that $b_{12}=\pm \im b_2$, $b_3=\pm \im b_9$, $b_5=0$. So it must be that $b_9=0$, $b_2=0$, and then $b_{12}=0$ and $b_{b_3}$. There is no non-trivial degeneration.

\

\noindent 7. $\sin(\alpha-\beta),\,f_{2,5},\,f_{3,4}$ are non zero \newline

In this case, the degeneration gives a non-trivial solution to the system of polynomial equations $\{p_1=0, p_2=0, q_7=0, q_8=0, q_5=0, q_6=0\} \cup \{q_2=0,q_3=0\} \cup \{q_{10}=0\} \cup \{q_{11}=0\}$.

\

By solving this system, we have $b_{10} = \frac{1}{2}\sqrt{20b_7^2+18b_8^2-4b_3^2+8b_6^2}$, $b_{11} = \frac{1}{2}\sqrt{-12b_7^2-10b_8^2-4b_2^2-8b_6^2}$, $b_5 = \frac{1}{2}\sqrt{-16b_7^2-14b_8^2-4b_6^2}$ and $b_6, b_7, b_8, b_2,$ $b_3 \in \RR$.

Therefore $b_6 = 0$, $b_7 = 0$, $b_8 = 0$, $b_2 = 0$ and it follows that $b_{11} = 0$, $b_5 = 0$, $b_{10} = \pm \im b_3$. So it must be that $b_3=0$, and then $b_{10}=0$. There is no non-trivial degeneration.

\

\noindent 8. $\sin(\alpha-\beta),\,f_{2,5},\,f_{3,5}$ are non zero \newline

Since $f_{2,4}$ and $f_{3,4}$ do not vanish simultaneously, this case is included in the case I) 6. and I) 7. There is no non-trivial degeneration.

\

\noindent II) $\cos(\beta)=0$ ($\sin(\beta)=\pm 1$)
$$
g_t\cdot\mu =\left\{\begin{array}{l}
{[e_1,e_2]}={\euler^{-t \left( a_1+a_2-a_4 \right) }}\pm \sin(\alpha) e_4-{\euler^{-t \left( a_1+a_2-a_5 \right) }}\mp \cos(\alpha) e_5,\\
{[e_1,e_3]}={\euler^{-t \left( a_1+a_3-a_4 \right) }}\mp \cos(\alpha) e_4+{\euler^{-t \left( a_1+a_3-a_5 \right) }}\pm \sin(\alpha) e_5,\\
{[e_1,e_5]}={\euler^{-t \left( a_1+a_5-a_6 \right) }}\pm e_6,\\
{[e_1,e_6]}={\euler^{-t \left( a_1+a_6-a_7 \right) }}e_7,[e_2,e_3]={\euler^{-t \left( a_2+a_3-a_6 \right) }}e_6,\\
{[e_2,e_4]}={\euler^{-t \left( a_2+a_4-a_7 \right) }} (\pm\sin(\alpha)\mp\cos(\alpha))e_7,\\
{[e_2,e_5]}={\euler^{-t \left( a_2+a_5-a_7 \right) }} \pm \lambda \cos(\alpha) e_7,\\
{[e_3,e_4]}={\euler^{-t \left( a_3+a_4-a_7 \right) }} (\mp \cos(\alpha) \mp \sin(\alpha)) e_7,\\
{[e_3,e_5]}={\euler^{-t \left( a_3+a_5-a_7 \right) }} \mp \lambda \sin(\alpha)e_7
\end{array}
\right.
$$

\noindent 1. $\sin(\alpha)$ and $\cos(\alpha)$ are non zero \newline

We recall that $f_{2,4}$ and $f_{3,4}$ do not vanish simultaneously. \newline

\noindent a. $f_{2,4}$ is no null \newline

In this case, the degeneration gives a non-trivial solution to the system of polynomial equations $\{p_1=0, p_2=0, q_7=0, q_8=0, q_6=0\} \cup \{q_1=0,q_2=0,q_3=0,q_4=0\} \cup \{q_{10}=0,q_{12}=0\} \cup \{q_{9}=0\}$.

\

By solving this system, we have
$$
\begin{array}{rl}
i)   & 8b_7^2 +7b_8^2 +4b_6^2-2b_1^2+2b_2^2=0 \\
ii)  & 10b_7^2+9b_8^2 +4b_6^2-2b_{10}^2-2b_3^2=0 \\
iii) & 2b_7^2 +2b_8^2 -2b_{10}^2-2b_4^2=0\\
iv)  & 10b_7^2+9b_8^2 +4b_6^2-2b_1^2-2b_{12}^2=0\\
v)   & 8b_7^2 +7b_8^2 +4b_6^2-2b_{10}^2+2b_9^2=0
\end{array}
$$
we solve the equation $iii)$ for $2b_{10}^2$  and substitute this expression in $v)$
$$
8b_7^2 +7b_8^2 +4b_6^2+(2b_4^2-2b_7^2-2b_8^2)+2b_9^2=0
$$
Therefore $b_7=0$, $b_8=0$, $b_6=0$, $b_4=0$ and $b_9=0$ and it follows that $b_{10}=0$, $b_{10}^2+b_3^2=0$, $b_1^2+b_{12}^2=0$, $b_1^2-b^2=0$. So it must be that $b_{10}=0$, $b_3=0$, $b_1=0$, $b_{12}=0$ and $b_2=0$. There is no non-trivial degeneration.

\

\noindent b. $f_{3,4}$ is no null \newline

In this case, the degeneration gives a non-trivial solution to the system of polynomial equations $\{p_1=0, p_2=0, q_7=0, q_8=0, q_6=0\} \cup \{q_1=0,q_2=0,q_3=0,q_4=0\} \cup \{q_{10}=0,q_{12}=0\} \cup \{q_{11}=0\}$.

\

By solving this system, we have
$$
\begin{array}{rl}
i)   & 8b_7^2 +7b_8^2 +4b_6^2-2b_1^2+2b_2^2=0 \\
ii)  & 10b_7^2+9b_8^2 +4b_6^2-2b_{10}^2-2b_3^2=0 \\
iii) & 2b_7^2 +2b_8^2 -2b_{10}^2-2b_4^2=0\\
iv)  & 10b_7^2+9b_8^2 +4b_6^2-2b_1^2-2b_{12}^2=0\\
v)   &  2b_7^2+2b_8^2 -2b_1^2-2b_{11}^2=0
\end{array}
$$
we solve the equation $v)$ for $2b_{1}^2$  and substitute this expression in $i)$
$$
8b_7^2 +7b_8^2 +4b_6^2+(2b_{11}^2-2b_7^2-2b_8^2)+2b_2^2=0
$$
Therefore $b_7=0$, $b_8=0$, $b_6=0$, $b_{11}=0$, $b_2=0$ and it follows that $b_1^2+b_{11}^2=0$, $b_1^2+b_{12}^2=0$, $b_{10}^2+b_{4}^2=0$, $b_{10}^2+b_3^2=0$ and $b_1^2-b_2^2=0$. So it must be that $b_1=0$, $b_{12}=0$, $b_{10}^2=0$, $b_{4}=0$, $b_3^2=0$ and $b_2=0$. There is no non-trivial degeneration.

\

\noindent 2. $\sin(\alpha)=0$ ($\cos(\alpha)=\pm 1$) \newline

In this case, the degeneration gives a non-trivial solution to the system of polynomial equations $\{p_1=0, p_2=0, q_7=0, q_8=0, q_6=0\} \cup \{ q_2=0,q_3=0\} \cup \{q_9=0,q_{10}=0,q_{11}=0\}$.

\

By solving this system, we have $b_{10} = \frac{1}{2}\sqrt{16b_7^2+14b_8^2+4b_9^2+8b_6^2}$, $b_2 = \frac{1}{2}\sqrt{-4b_{11}^2-8b_6^2-10b_8^2-12b_7^2}$, $b_3 = \sqrt{-b_9^2+b_8^2+b_7^2}$ and $b_6, b_7, b_8, b_9, b_{11} \in \RR$. Therefore $b_6 = 0$, $b_7 = 0$, $b_8 = 0$, $b_9 = 0$ and it follows that $b_{10}=0$, $b_3 = 0$ and $b_2 = \pm \im b_{11}$. So it must be that $b_{11}=0$. There is no non-trivial degeneration.

\

\noindent 3. $\cos(\alpha)=0$ ($\sin(\alpha)=\pm 1$) \newline

In this case, the degeneration gives a non-trivial solution to the system of polynomial equations $\{p_1=0, p_2=0, q_7=0, q_8=0, q_6=0\} \cup \{ q_1=0,q_4=0\} \cup \{q_9=0,q_{11}=0,q_{12}=0\}$.

\

By solving this system, we have $b_1 = \sqrt{-b_{11}^2+b_8^2+b_7^2}$, \newline $b_{12} = \frac{1}{2}\sqrt{16b_7^2+14b_8^2+4b_{11}^2+8b_6^2}$,  $b_9 = \frac{1}{2}\sqrt{-4b_4^2-8b_6^2-10b_8^2-12b_7^2}$ and $b_{11}, b_4, b_6, b_7, b_8 \in \RR$. Therefore $b_4 = 0$ $b_6 = 0$, $b_7 = 0,$ $b_8 = 0$ and it follows that $b_9 = 0$, $b_{12} = b_{11}$, and $b_1 = \pm \im b_{11}$. So it must be that $b_{11}=0$ and then $b_{12}=0$, $b_1=0$. There is no non-trivial degeneration.

\

\noindent III) $\sen(\beta)=0$ ($\cos(\beta)=\pm 1$)

$$
g_t\cdot\mu =\left\{\begin{array}{l}
{[e_1,e_2]}={\euler^{-t \left( a_1+a_2-a_4 \right) }}\pm\cos \left( \alpha\right) e_4-{\euler^{-t \left( a_1+a_2-a_5 \right) }}\pm\sin\left( \alpha\right) e_5,\\
{[e_1,e_3]}={\euler^{-t \left( a_1+a_3-a_4 \right) }}\pm\sin \left( \alpha\right) e_4+{\euler^{-t \left( a_1+a_3-a_5 \right) }}\pm\cos\left( \alpha\right) e_5,\\
{[e_1,e_4]}={\euler^{-t \left( a_1+a_4-a_6 \right) }}\pm e_6,\\
{[e_1,e_6]}={\euler^{-t \left( a_1+a_6-a_7 \right) }}e_7,[e_2,e_3]={\euler^{-t \left( a_2+a_3-a_6 \right) }}e_6,\\
{[e_2,e_4]}={\euler^{-t \left( a_2+a_4-a_7 \right) }} \pm \cos(\alpha)\lambda e_7,\\
{[e_2,e_5]}={\euler^{-t \left( a_2+a_5-a_7 \right) }} (\mp \sin(\alpha) \pm \cos(\alpha)) e_7,\\
{[e_3,e_4]}={\euler^{-t \left( a_3+a_4-a_7 \right) }} \pm \sin(\alpha)\lambda e_7,\\
{[e_3,e_5]}={\euler^{-t \left( a_3+a_5-a_7 \right) }} (\pm \cos(\alpha) \pm \sin(\alpha)) e_7
\end{array}
\right.
$$

1. $\sin(\alpha)$ and $\cos(\alpha)$ are non zero \newline

We recall that $f_{2,5}$ and $f_{3,5}$ do not vanish simultaneously.

\

\noindent a. $f_{2,5}$ is no null \newline

In this case, the degeneration gives a non-trivial solution to the system of polynomial equations $\{p_1=0, p_2=0, q_7=0, q_8=0, q_5=0\} \cup \{q_1=0,q_2=0,q_3=0,q_4=0\} \cup \{q_{9}=0,q_{11}=0\} \cup \{q_{10}=0\}$.
$$
\begin{array}{rl}
  i) & 8b_7^2+7b_8^2+4b_5^2-2b_2^2+2b_1^2 =0\\
 ii) & 2b_7^2+2b_8^2-2b_9^2-2b_3^2=0\\
iii) & 10b_7^2+9b_8^2+4b_5^2-2b_9^2-2b_4^2=0\\
 iv) & b_7^2+b_8^2-b_1^2-b_{11}^2=0\\
  v) & 8b_7^2+7b_8^2-2b_9^2+4b_5^2+2b_{10}^2=0
\end{array}
$$
We solve the equation $ii)$ for $2b_{9}^2$  and substitute this expression in $v)$

$$8b_7^2+7b_8^2+(2b_3^2-2b_7^2-2b_8^2)+4b_5^2+2b_{10}^2=0$$

Therefore $b_7=0$, $b_8=0$, $b_3=0$, $b_5=0$, $b_{10}=0$  and it follows that $b_1^2+b_{11}^2=0$, $b_9^2+b_4^2$, $b_9^2+b_3^2=0$ and $b_2^2+b_1^2=0$.
So it must be that $b_1^2=0$, $b_{11}^2=0$, $b_9^2=0$, $b_4^2=0$, $b_2^2=0$. There is no non-trivial degeneration.

\

\noindent b. $f_{3,5}$ is no null \newline

In this case, the degeneration gives a non-trivial solution to the system of polynomial equations $\{p_1=0, p_2=0, q_7=0, q_8=0, q_5=0\} \cup \{q_1=0,q_2=0,q_3=0,q_4=0\} \cup \{q_{9}=0,q_{11}=0\} \cup \{q_{10}=0\}$.

\

By solving this system, we have
$$
\begin{array}{rl}
  i) & 8b_7^2+7b_8^2+4b_5^2-2b_2^2+2b_1^2 =0\\
 ii) & 2b_7^2+2b_8^2-2b_9^2-2b_3^2=0\\
iii) & 10b_7^2+9b_8^2+4b_5^2-2b_9^2-2b_4^2=0\\
 iv) & b_7^2+b_8^2-b_1^2-b_{11}^2=0\\
  v) & 6b_7^2+5b_8^2+2b_1^2+4b_5^2+2b_{12}^2=0
\end{array}
$$

By $v)$ it is follows that  $b_7=0$, $b_8=0$, $b_1=0$, $b_5=0$, $b_{12}=0$ and so $b_2^2+b_1^2=0$, $b_9^2+b_3^2=0$, $b_9^2+b_4^2=0$, $b_{11}^2=0$. So it must be that $b_1=0$, $b_2=0$, $b_3=0$, $b_4=0$, $b_9=0$ and $b_{11}=0$. There is no non-trivial degeneration.

\

\noindent 2. $\sin(\alpha)=0$ ($\cos(\alpha)=\pm 1$) \newline

In this case, the degeneration gives a non-trivial solution to the system of polynomial equations $\{p_1=0, p_2=0, q_7=0, q_8=0, q_5=0\} \cup \{ q_1=0,q_4=0\} \cup \{q_9=0,q_{10}=0,q_{12}=0\}$.

\

By solving this system, we have  $b_{12} = \frac{1}{2}\sqrt{-4b_1^2-8b_5^2-10b_8^2-12b_7^2}$, $b_4 = \sqrt{-b_{10}^2+b_8^2+b_7^2}$,  $b_9 = \frac{1}{2}\sqrt{16b_7^2+14b_8^2+4b_{10}^2+8b_5^2}$, $b_1, b_{10}, b_5, b_7, b_8 \in \RR$. Therefore $b_5 = 0$, $b_7 = 0$, $b_8 = 0$, $b_1 = 0$ and it follows that $b_{12} = 0$, $b_4 = \pm \im b_{10}$, $b_9 = b_{10}$. So it must be that $b_{10}=0$ and then $b_4=0$, $b_9=0$. There is no non-trivial degeneration.

\

\noindent 3. $\cos(\alpha)=0$ ($\sin(\alpha)=\pm 1$) \newline

In this case, the degeneration gives a non-trivial solution to the system of polynomial equations $\{p_1=0, p_2=0, q_7=0, q_8=0, q_5=0\} \cup \{ q_2=0,q_3=0\} \cup \{q_{10}=0,q_{12}=0,q_{11}=0\}$.

By solving this system, we have $b_{11} = \frac{1}{2}\sqrt{16b_7^2+14b_8^2+4b_{12}^2+8b_5^2}$, $b_2 = \sqrt{-b_{12}^2+b_8^2+b_7^2}$, $b_3 = \frac{1}{2}\sqrt{-4b_{10}^2-8b_5^2-10b_8^2-12b_7^2}$, $b_5, b_7, b_8, b_{12}, b_{10}$. Therefore $b_5 = 0$, $b_7 = 0$ $b_8 = 0$, $b_{10} = 0$ and it follows that $b_{12} = b_{11}$, $b_2 = \pm \im b_{11}$, $b_3 = 0$. So it must be that $b_{11}=0$ and then $b_{12}=0$, $b_2=0$. There is no non-trivial degeneration.

\

In any case, There is no non-trivial degeneration. So, $1.3(i_\lambda)$ must be Einstein nilradical for any $\lambda \neq 0$

\

\noindent Pre-Einstein Derivation: $\frac{5 }{17 }\diag( 1, 2, 2, 3, 3, 4, 5 )$ \newline
\noindent $||\mathcal{S}_{\beta}||^2=\frac{17}{19}\thickapprox 0.895$

\

\begin{center}
\begin{tabular}{cccccccccccc}
\hline
\multicolumn{3}{||m{1.5cm}|  }  {\centering  $1.3(i_\lambda)$ \\ $\lambda = 0$ } &
\multicolumn{3}{  m{1.5cm}|  }  {\centering  $\mathfrak{n}_{7,62}^{\lambda}$ \\ $\lambda = 0$ } &
\multicolumn{3}{  m{1.5cm}|  }  {\centering  $L_{74}^{15}(3\CC)$ \\ \cite{CARLES1}} &
\multicolumn{3}{  m{4.25cm}||}  {\centering  $(1,3,5,7)N$ \\ $(\xi=0)$} \tabularnewline
\hline
\hline
\multicolumn{3}{||m{1.5cm}|   }  {\centering $\dim$ $\Der$ \\               $13$    } &
\multicolumn{6}{  m{3.75cm}  |   }  {\centering $\dim$ {Derived Series}  \\     $(7, 4, 0)$    } &
\multicolumn{3}{  m{4.25cm}|| }  {\centering $\dim$ {Desc. C. Series} \\     $(7, 4, 2, 1, 0)$    } \tabularnewline
\hline
\end{tabular}
\end{center}

\

\noindent  $[e_1,e_2]=e_4,[e_1,e_3]=e_5,[e_1,e_4]=e_6,[e_1,e_6]=e_7,[e_2,e_3]=e_6,[e_2,e_5]=e_7,[e_3,e_5]=e_7$

\

\noindent It is not an Einstein nilradical by \cite[Theorem 2.6]{FERNANDEZ-CULMA1}. In fact,

\

\noindent  Einstein Pre-Derivation $\frac{5 }{17 }\diag ( 1, 2, 2, 3, 3, 4, 5 )$. \newline
Let $X=\diag \left(1,\left(\begin{array}{cc} 0 & 2\\ 2 & 0 \end{array}\right) , -1, -1, 0, 1 \right) $. $X\in \ggo_{\phi}$ therefore $g_{t}=\exp(tX) \in \Grp_{\phi}$.

\

\noindent  $g_t\cdot\mu=\left\{\begin{array}{l}
{[e_1,e_2]}=\frac{1}{2} ( 1 + {{\euler}^{-4t}} ) e_4 +\frac{1}{2} ( {{\euler}^{-4t}} - 1 ) e_5,\\
{[e_1,e_3]}=\frac{1}{2} ( {{\euler}^{-4t}} - 1 ) e_4 +\frac{1}{2} ( 1 + {{\euler}^{-4t}} ) e_5,\\
{[e_1,e_4]}=e_6, [e_1,e_6]=e_7,[e_2,e_3]=e_6,[e_2,e_5]=e_7,[e_3,e_5]=e_7
\end{array}
\right.$

\

\noindent The $\Grp_{\phi}-$orbit of $\mu$ is not closed because  $g_t\cdot\mu \rightarrow \widetilde{\mu}$ as $t \rightarrow \infty$ and $(\RR^7,\widetilde{\mu})$ is non-isomorphic to $(\RR^7,\mu)$; $\dim \Der (\RR^7,\mu) = 13$ and $\dim \Der (\RR^7,\widetilde{\mu}) = 14$

\

\begin{center}
\begin{tabular}{cccccccccccc}
\hline
\multicolumn{3}{||m{1.5cm}|  }  {\centering  $1.3(ii)$  } &
\multicolumn{3}{  m{1.5cm}|  }  {\centering  $\mathfrak{n}_{7,63}^{ }$ } &
\multicolumn{3}{  m{1.5cm}|  }  {\centering   } &
\multicolumn{3}{  m{4.25cm}||}  {\centering  $(1,3,5,7)L$  } \tabularnewline
\hline
\hline
\multicolumn{3}{||m{1.5cm}|   }  {\centering $\dim$ $\Der$ \\               $14$    } &
\multicolumn{6}{  m{3.75cm}  |   }  {\centering $\dim$ {Derived Series}  \\     $(7, 4, 0)$    } &
\multicolumn{3}{  m{4.25cm}|| }  {\centering $\dim$ {Desc. C. Series} \\     $(7, 4, 2, 1, 0)$    } \tabularnewline
\hline
\end{tabular}
\end{center}

\

\noindent $\mu:=\left\{\begin{array}{l}
[e_1,e_2]=e_4, [e_1,e_3]=e_5,  [e_1,e_4]=e_6,  [e_1,e_6]=e_7, {[e_2,e_3]}=e_6, \\
{[e_2,e_4]}=e_7,  [e_2,e_5]=\frac{1}{2}e_7,  [e_3,e_4]=-\frac{1}{2}e_7
\end{array}\right.$

\

\noindent It is not an Einstein nilradical by \cite[Theorem 2.6]{FERNANDEZ-CULMA1}. In fact,

\

\noindent  Einstein Pre-Derivation $\phi=\frac{5}{17}\diag(1,2,2,3,3,4,5) $. \newline
Let $X=\diag (-145, 730, -294, 293, -731, 148, -1)$. $X\in \ggo_{\phi}$ therefore $g_{t}=\exp(tX) \in \Grp_{\phi}$.

\

\noindent  $g_t\cdot\mu=\left\{\begin{array}{l}
[e_1,e_2]={\euler^{-292\,t}} e_4, [e_1,e_3]={\euler^{-292\,t}} e_5, [e_1,e_4]= e_6, \\
{[e_1,e_6]}={\euler^{-4\,t}} e_7, [e_2,e_3]={\euler^{-288\,t}} e_6, [e_2,e_4]={\euler^{-1024\,t}} e_7,\\
{[e_2,e_5]}=\frac{1}{2} e_7, [e_3,e_4]=-\frac{1}{2} e_7
\end{array}\right.$

\

\noindent The $\Grp_{\phi}-$orbit of $\mu$ is not closed because  $g_t\cdot\mu \rightarrow \widetilde{\mu}$ as $t \rightarrow \infty$ and $(\RR^7,\widetilde{\mu})$ is non-isomorphic to $(\RR^7,\mu)$; $\dim \Der (\RR^7,\mu) = 14$ and $\dim \Der (\RR^7,\widetilde{\mu}) = 21$

\

\begin{center}
\begin{tabular}{cccccccccccc}
\hline
\multicolumn{3}{||m{1.5cm}|  }  {\centering  $1.3(iii)$  } &
\multicolumn{3}{  m{1.5cm}|  }  {\centering  $\mathfrak{n}_{7,65}^{}$ } &
\multicolumn{3}{  m{1.5cm}|  }  {\centering  $L_{74}^{18}(3\CC)$ } &
\multicolumn{3}{  m{4.25cm}||}  {\centering  $(1,3,5,7)F$  } \tabularnewline
\hline
\hline
\multicolumn{3}{||m{1.5cm}|   }  {\centering $\dim$ $\Der$ \\               $13$    } &
\multicolumn{6}{  m{3.75cm}  |   }  {\centering $\dim$ {Derived Series}  \\     $(7, 4, 0)$    } &
\multicolumn{3}{  m{4.25cm}|| }  {\centering $\dim$ {Desc. C. Series} \\     $(7, 4, 2, 1, 0)$    } \tabularnewline
\hline
\end{tabular}
\end{center}

\

\noindent  $[e_1,e_2]=e_4,[e_1,e_3]=e_5,[e_1,e_4]=e_6,[e_1,e_6]=e_7,[e_2,e_4]=e_7,[e_3,e_5]=e_7$

\

\noindent  $\{e_1...e_7\}$ is a nice basis.

$$U=\tiny{ \left[ \begin {array}{cccccc} 3&1&0&1&0&0\\ \noalign{\medskip}1&3&1&1
&0&0\\ \noalign{\medskip}0&1&3&0&1&0\\ \noalign{\medskip}1&1&0&3&1&1
\\ \noalign{\medskip}0&0&1&1&3&1\\ \noalign{\medskip}0&0&0&1&1&3
\end {array} \right]}$$

\noindent  General solution to $Ux=[1]$: $x=\frac{1}{34}(9,5,8,2,5,9)^T$. It is an Einstein nilradical.

\

\noindent  Pre-Einstein Derivation: $\frac{5}{17} \diag(1, 2, 2, 3, 3, 4, 5)$. \newline
\noindent $||\mathcal{S}_{\beta}||^2=\frac{17}{19}\thickapprox 0.895$

\

\begin{center}
\begin{tabular}{cccccccccccc}
\hline
\multicolumn{3}{||m{1.5cm}|  }  {\centering  $1.3(iv)$  } &
\multicolumn{3}{  m{1.5cm}|  }  {\centering  $\mathfrak{n}_{7,87}^{}$ } &
\multicolumn{3}{  m{1.5cm}|  }  {\centering  $L_{74}^{35}(3\CC)$ } &
\multicolumn{3}{  m{4.25cm}||}  {\centering  $(2,4,7)R$  } \tabularnewline
\hline
\hline
\multicolumn{3}{||m{1.5cm}|   }  {\centering $\dim$ $\Der$ \\               $13$    } &
\multicolumn{6}{  m{3.75cm}  |   }  {\centering $\dim$ {Derived Series}  \\     $(7, 4, 0)$    } &
\multicolumn{3}{  m{4.25cm}|| }  {\centering $\dim$ {Desc. C. Series} \\     $(7, 4, 2, 0)$    } \tabularnewline
\hline
\end{tabular}
\end{center}

\

\noindent  $[e_1,e_2]=e_4,[e_1,e_3]=e_5,[e_1,e_4]=e_6,[e_2,e_3]=e_6,[e_2,e_4]=e_7,[e_3,e_5]=e_7$

\

\noindent  $\{e_1...e_7\}$ is a nice basis.

$$U=\tiny{ \left[ \begin {array}{cccccc} 3&1&0&1&0&0\\ \noalign{\medskip}1&3&1&1
&0&0\\ \noalign{\medskip}0&1&3&1&1&0\\ \noalign{\medskip}1&1&1&3&1&1
\\ \noalign{\medskip}0&0&1&1&3&1\\ \noalign{\medskip}0&0&0&1&1&3
\end {array} \right]}$$

\noindent  General solution to $Ux=[1]$: $x=\frac{1}{17}(5,3,4,-1,3,5)^T$. It is not an Einstein nilradical.

\noindent  Pre-Einstein Derivation: $\frac{5}{17} (1, 2, 2, 3, 3, 4, 5) $

\

\begin{center}
\begin{tabular}{cccccccccccc}
\hline
\multicolumn{3}{||m{1.5cm}|  }  {\centering  $1.3(v)$  } &
\multicolumn{3}{  m{1.5cm}|  }  {\centering  $\mathfrak{n}_{7,99}^{}$ } &
\multicolumn{3}{  m{1.5cm}|  }  {\centering  $L_{74}^{20}(\mathfrak{n}_3)$ } &
\multicolumn{3}{  m{4.25cm}||}  {\centering  $(1,3,5,7)C$  } \tabularnewline
\hline
\hline
\multicolumn{3}{||m{1.5cm}|   }  {\centering $\dim$ $\Der$ \\               $13 $    } &
\multicolumn{6}{  m{3.75cm}  |   }  {\centering $\dim$ {Derived Series}  \\     $(7, 3, 0)$    } &
\multicolumn{3}{  m{4.25cm}|| }  {\centering $\dim$ {Desc. C. Series} \\     $(7, 3, 2, 1, 0)$    } \tabularnewline
\hline
\end{tabular}
\end{center}

\

\noindent  $\mu:=\left\{
\begin{array}{l}
[e_1,e_2]=e_4, [e_1,e_4]=e_6, [e_1,e_6]=e_7, [e_2,e_3]=e_6, [e_2,e_4]=e_7, \\
{[e_3,e_4]}=-e_7,  [e_3,e_5]=-e_7
\end{array}\right.$

\

\noindent It is not an Einstein nilradical by \cite[Theorem 2.6]{FERNANDEZ-CULMA1}. In fact,

\

\noindent  Einstein Pre-Derivation $\phi=\frac{5}{17} \diag ( 1, 2, 2, 3, 3, 4, 5 ) $. \newline
Let $X=\diag (-1,2,-2,1,1,0,-1) $. $X\in \ggo_{\phi}$ therefore $g_{t}=\exp(tX) \in \Grp_{\phi}$.

\

\noindent  $g_t\cdot\mu=\left\{\begin{array}{l}
[e_1,e_2]=e_4,  [e_1,e_4]=e_6,  [e_1,e_6]=e_7, [e_2,e_3]=e_6, \\
{[e_2,e_4]}=\euler^{-4t}e_7,  [e_3,e_4]=-e_7,  [e_3,e_5]=-e_7
\end{array}\right.$

\

\noindent The $\Grp_{\phi}-$orbit of $\mu$ is not closed because  $g_t\cdot\mu \rightarrow \widetilde{\mu}$ as $t \rightarrow \infty$ and $(\RR^7,\widetilde{\mu})$ is non-isomorphic to $(\RR^7,\mu)$; $(\RR^7,\mu)$ has $\rank=1$ and $(\RR^7,\widetilde{\mu})$ has $\rank= 2$ with maximal torus generated by $D1=\diag(1,0,2,1,1,2,3) $ and $D2=\diag(0,1,0,1,1,1,1)$

\

\begin{center}
\begin{tabular}{cccccccccccc}
\hline
\multicolumn{3}{||m{1.5cm}|  }  {\centering  $1.4$  } &
\multicolumn{3}{  m{1.5cm}|  }  {\centering  $\mathfrak{n}_{7,8}^{}$ } &
\multicolumn{3}{  m{1.5cm}|  }  {\centering  $L_{75}^{31}$ } &
\multicolumn{3}{  m{4.25cm}||}  {\centering  $(1, 2, 3, 4, 5, 7)$  } \tabularnewline
\hline
\hline
\multicolumn{3}{||m{1.5cm}|   }  {\centering $\dim$ $\Der$ \\               $12$    } &
\multicolumn{6}{  m{3.75cm}  |   }  {\centering $\dim$ {Derived Series}  \\     $(7, 5, 0)$    } &
\multicolumn{3}{  m{4.25cm}|| }  {\centering $\dim$ {Desc. C. Series} \\     $(7, 5, 4, 3, 2, 1, 0)$    } \tabularnewline
\hline
\end{tabular}
\end{center}

\

\noindent  $[e_1,e_2]=e_3,[e_1,e_3]=e_4,[e_1,e_4]=e_5,[e_1,e_5]=e_6,[e_1,e_6]=e_7,[e_2,e_3]=e_6,[e_2,e_4]=e_7$

\

\noindent  $\{e_1...e_7\}$ is a nice basis.

$$U=\tiny{\left[ \begin {array}{ccccccc} 3&0&1&1&1&0&1\\ \noalign{\medskip}0&3&0
&1&1&1&-1\\ \noalign{\medskip}1&0&3&0&1&0&1\\ \noalign{\medskip}1&1&0&
3&0&1&0\\ \noalign{\medskip}1&1&1&0&3&-1&1\\ \noalign{\medskip}0&1&0&1
&-1&3&1\\ \noalign{\medskip}1&-1&1&0&1&1&3\end {array} \right]}$$

\noindent A solution to $Ux=[1]$: $x=\frac{1}{100}(10,21,18,16,18,21,18)^T$. It is an Einstein nilradical.

\

\noindent Pre-Einstein Derivation: $\frac{17}{100}\diag(1, 3, 4, 5, 6, 7, 8)$. \newline
$||\mathcal{S}_{\beta}||^2 = \frac{50}{61} \thickapprox 0.8196721312$

\

\begin{center}
\begin{tabular}{cccccccccccc}
\hline
\multicolumn{3}{||m{1.5cm}|  }  {\centering  $1.5$  } &
\multicolumn{3}{  m{1.5cm}|  }  {\centering  $\mathfrak{n}_{7,12}^{}$ } &
\multicolumn{3}{  m{1.5cm}|  }  {\centering  $L_{74}^{36}(\mathfrak{n}_3)$ } &
\multicolumn{3}{  m{4.25cm}||}  {\centering  $(2,3,4,5,7)D$  } \tabularnewline
\hline
\hline
\multicolumn{3}{||m{1.5cm}|   }  {\centering $\dim$ $\Der$ \\               $ 11$    } &
\multicolumn{6}{  m{3.75cm}  |   }  {\centering $\dim$ {Derived Series}  \\     $(7, 5, 1, 0)$    } &
\multicolumn{3}{  m{4.25cm}|| }  {\centering $\dim$ {Desc. C. Series} \\     $(7, 5, 4, 3, 2, 0)$    } \tabularnewline
\hline
\end{tabular}
\end{center}

\

\noindent  $[e_1,e_2]=e_3,[e_1,e_3]=e_4,[e_1,e_4]=e_5,[e_1,e_5]=e_6,[e_2,e_3]=e_6,[e_2,e_5]=-e_7,[e_3,e_4]=e_7$

\

\noindent  $\{e_1...e_7\}$ is a nice basis.

$$U=\tiny{\left[ \begin {array}{ccccccc} 3&0&1&1&0&1&-1\\ \noalign{\medskip}0&3
&0&1&1&0&0\\ \noalign{\medskip}1&0&3&0&0&-1&1\\ \noalign{\medskip}1&1&0
&3&1&1&0\\ \noalign{\medskip}0&1&0&1&3&1&1\\ \noalign{\medskip}1&0&-1&
1&1&3&1\\ \noalign{\medskip}-1&0&1&0&1&1&3\end {array} \right]}$$

\noindent  A solution to $Ux=[1]$: $x=\frac{1}{31}(7,9,8,2,2,7,7 )^T$. It is an Einstein nilradical.

\noindent  Pre-Einstein Derivation: $\frac{5}{31}\diag(1, 3, 4, 5, 6, 7, 9) $. \newline
$||\mathcal{S}_{\beta}||^2 = \frac{31}{42} \thickapprox 0.7380952381$

\

\begin{center}
\begin{tabular}{cccccccccccc}
\hline
\multicolumn{3}{||m{1.5cm}|  }  {\centering  $1.6$  } &
\multicolumn{3}{  m{1.5cm}|  }  {\centering  $\mathfrak{n}_{7,9}^{}$ } &
\multicolumn{3}{  m{1.5cm}|  }  {\centering  $L_{75}^{32}$ } &
\multicolumn{3}{  m{4.25cm}||}  {\centering  $(1,2,3,4,5,7)B$  } \tabularnewline
\hline
\hline
\multicolumn{3}{||m{1.5cm}|   }  {\centering $\dim$ $\Der$ \\               $12$    } &
\multicolumn{6}{  m{3.75cm}  |   }  {\centering $\dim$ {Derived Series}  \\     $(7, 5, 0)$    } &
\multicolumn{3}{  m{4.25cm}|| }  {\centering $\dim$ {Desc. C. Series} \\     $(7, 5, 4, 3, 2, 1, 0)$    } \tabularnewline
\hline
\end{tabular}
\end{center}

\

\noindent  $[e_1,e_2]=e_3,[e_1,e_3]=e_4,[e_1,e_4]=e_5,[e_1,e_5]=e_6,[e_1,e_6]=e_7,[e_2,e_3]=e_7$

\

\noindent  $\{e_1...e_7\}$ is a nice basis.

$$U=\tiny{\left[ \begin {array}{cccccc} 3&0&1&1&1&0\\ \noalign{\medskip}0&3&0&1
&1&1\\ \noalign{\medskip}1&0&3&0&1&0\\ \noalign{\medskip}1&1&0&3&0&0
\\ \noalign{\medskip}1&1&1&0&3&1\\ \noalign{\medskip}0&1&0&0&1&3
\end {array} \right]}$$

\noindent General solution to $Ux=[1]$: $x=\frac{1}{34}(5,5,9,8,2,9)^T$. It is an Einstein nilradical.

\

\noindent Pre-Einstein Derivation: $\frac{5}{34} \diag(1, 4, 5, 6, 7, 8, 9)$. \newline
\noindent $||\mathcal{S}_{\beta}||^2 = \frac{17}{19} \thickapprox 0.8947368421$

\

\begin{center}
\begin{tabular}{cccccccccccc}
\hline
\multicolumn{3}{||m{1.5cm}|  }  {\centering  $1.7$  } &
\multicolumn{3}{  m{1.5cm}|  }  {\centering   } &
\multicolumn{3}{  m{1.5cm}|  }  {\centering  $L_{75}^{56}$ } &
\multicolumn{3}{  m{4.25cm}||}  {\centering  $(2, 4, 7)O$  } \tabularnewline
\hline
\hline
\multicolumn{3}{||m{1.5cm}|   }  {\centering $\dim$ $\Der$ \\               $15$    } &
\multicolumn{6}{  m{3.75cm}  |   }  {\centering $\dim$ {Derived Series}  \\     $(7, 4, 0)$    } &
\multicolumn{3}{  m{4.25cm}|| }  {\centering $\dim$ {Desc. C. Series} \\     $(7, 4, 2, 0)$    } \tabularnewline
\hline
\end{tabular}
\end{center}

\

\noindent  $[e_1,e_2]=e_4,[e_1,e_3]=e_5,[e_1,e_4]=e_6,[e_1,e_5]=e_7,[e_2,e_3]=e_6,[e_2,e_4]=e_7$

\

\noindent  $\{e_1...e_7\}$ is a nice basis.

$$U=\tiny{\left[ \begin {array}{cccccc} 3&1&0&1&1&0\\ \noalign{\medskip}1&3&1&0
&1&0\\ \noalign{\medskip}0&1&3&1&1&1\\ \noalign{\medskip}1&0&1&3&0&1
\\ \noalign{\medskip}1&1&1&0&3&1\\ \noalign{\medskip}0&0&1&1&1&3
\end {array} \right]}$$

\noindent  General solution to $Ux=[1]$: $x=\frac{1}{29}(5,6,3,5,3,6)^T$. It is an Einstein nilradical.

\noindent Pre-Einstein Derivation: $\frac{5}{29} \diag(2, 3, 4, 5, 6, 7, 8)$. \newline
\noindent \noindent $||\mathcal{S}_{\beta}||^2=\frac{29}{28}\thickapprox 1.04$

\

\begin{center}
\begin{tabular}{cccccccccccc}
\hline
\multicolumn{3}{||m{1.5cm}|  }  {\centering  $1.8$  } &
\multicolumn{3}{  m{1.5cm}|  }  {\centering  $\mathfrak{n}_{7,70}^{}$ } &
\multicolumn{3}{  m{1.5cm}|  }  {\centering  $L_{74}^{14}(3\CC)$ } &
\multicolumn{3}{  m{4.25cm}||}  {\centering  $(1,3,5,7)J$  } \tabularnewline
\hline
\hline
\multicolumn{3}{||m{1.5cm}|   }  {\centering $\dim$ $\Der$ \\               $11$    } &
\multicolumn{6}{  m{3.75cm}  |   }  {\centering $\dim$ {Derived Series}  \\     $(7, 4, 0)$    } &
\multicolumn{3}{  m{4.25cm}|| }  {\centering $\dim$ {Desc. C. Series} \\     $(7, 4, 2, 1, 0)$    } \tabularnewline
\hline
\end{tabular}
\end{center}

\

\noindent  $[e_1,e_2]=e_4,[e_1,e_4]=e_6,[e_1,e_6]=e_7,[e_2,e_3]=e_5,[e_2,e_4]=e_7,[e_3,e_5]=e_7$

\

\noindent  $\{e_1...e_7\}$ is a nice basis.

$$U=\tiny{\left[ \begin {array}{cccccc} 3&0&1&1&0&0\\ \noalign{\medskip}0&3&0&0
&1&0\\ \noalign{\medskip}1&0&3&0&1&1\\ \noalign{\medskip}1&0&0&3&1&0
\\ \noalign{\medskip}0&1&1&1&3&1\\ \noalign{\medskip}0&0&1&0&1&3
\end {array} \right]}$$

\noindent  General solution to $Ux=[1]$: $x=\frac{1}{139}(24,48,27,40,-5,39)^T$. It is not an Einstein nilradical.

\

\noindent Pre-Einstein Derivation: $ \frac{20}{139}\diag(2, 4, 3, 6, 7, 8, 10)$

\

\begin{center}
\begin{tabular}{cccccccccccc}
\hline
\multicolumn{3}{||m{1.5cm}|  }  {\centering  $1.9$  } &
\multicolumn{3}{  m{1.5cm}|  }  {\centering  $\mathfrak{n}_{7,59}^{}$ } &
\multicolumn{3}{  m{1.5cm}|  }  {\centering  $L_{75}^{55}$ } &
\multicolumn{3}{  m{4.25cm}||}  {\centering  $(2,4,5,7)K$  } \tabularnewline
\hline
\hline
\multicolumn{3}{||m{1.5cm}|   }  {\centering $\dim$ $\Der$ \\               $14$    } &
\multicolumn{6}{  m{3.75cm}  |   }  {\centering $\dim$ {Derived Series}  \\     $(7, 4, 0)$    } &
\multicolumn{3}{  m{4.25cm}|| }  {\centering $\dim$ {Desc. C. Series} \\     $(7, 4, 3, 1, 0)$    } \tabularnewline
\hline
\end{tabular}
\end{center}

\

\noindent  $[e_1,e_2]=e_4,[e_1,e_3]=e_6,[e_1,e_4]=e_5,[e_1,e_5]=e_7,[e_2,e_3]=e_7,[e_2,e_4]=e_6$

\

\noindent  $\{e_1...e_7\}$ is a nice basis.

$$U=\tiny{ \left[ \begin {array}{cccccc} 3&1&0&1&1&0\\ \noalign{\medskip}1&3&1&1
&1&1\\ \noalign{\medskip}0&1&3&0&0&1\\ \noalign{\medskip}1&1&0&3&1&0
\\ \noalign{\medskip}1&1&0&1&3&1\\ \noalign{\medskip}0&1&1&0&1&3
\end {array} \right]}$$

\noindent General solution to $Ux=[1]$: $x=\frac{1}{67}(15,-1,18,15,8,14 )^T$. It is not an Einstein nilradical.

\

\noindent Pre-Einstein Derivation: $ \frac{10}{67}\diag(2, 3, 6, 5, 7, 8, 9) $

\

\begin{center}
\begin{tabular}{cccccccccccc}
\hline
\multicolumn{3}{||m{1.5cm}|  }  {\centering  $1.10$  } &
\multicolumn{3}{  m{1.5cm}|  }  {\centering  $\mathfrak{n}_{7,29}^{}$ } &
\multicolumn{3}{  m{1.5cm}|  }  {\centering  $L_{75}^{36}$ } &
\multicolumn{3}{  m{4.25cm}||}  {\centering  $(1, 3, 4, 5, 7)F$  } \tabularnewline
\hline
\hline
\multicolumn{3}{||m{1.5cm}|   }  {\centering $\dim$ $\Der$ \\               $11$    } &
\multicolumn{6}{  m{3.75cm}  |   }  {\centering $\dim$ {Derived Series}  \\     $(7, 5, 0)$    } &
\multicolumn{3}{  m{4.25cm}|| }  {\centering $\dim$ {Desc. C. Series} \\     $(7, 5, 4, 2, 1, 0)$    } \tabularnewline
\hline
\end{tabular}
\end{center}

\

\noindent  $[e_1,e_2]=e_3,[e_1,e_3]=e_4,[e_1,e_4]=e_6,[e_1,e_6]=e_7,[e_2,e_3]=e_5,[e_2,e_5]=e_7$

\

\noindent  $\{e_1...e_7\}$ is a nice basis.

$$U=\tiny{\left[ \begin {array}{cccccc} 3&0&1&1&0&1\\ \noalign{\medskip}0&3&0&1
&1&0\\ \noalign{\medskip}1&0&3&0&0&0\\ \noalign{\medskip}1&1&0&3&0&1
\\ \noalign{\medskip}0&1&0&0&3&0\\ \noalign{\medskip}1&0&0&1&0&3
\end {array} \right]}$$

\noindent  General solution to $Ux=[1]$: $x=\frac{1}{353}(35,68,106,54,95,88)^T$. It is  Einstein nilradical.

\

\noindent  Pre-Einstein Derivation: $ \frac{45}{353}\diag(2, 3, 5, 7, 8, 9, 11)$. \newline
\noindent $||\mathcal{S}_{\beta}||^2=\frac{353}{446}\thickapprox 0.792$

\

\begin{center}
\begin{tabular}{cccccccccccc}
\hline
\multicolumn{3}{||m{1.5cm}|  }  {\centering  $1.11$  } &
\multicolumn{3}{  m{1.5cm}|  }  {\centering  $\mathfrak{n}_{7,46}^{}$ } &
\multicolumn{3}{  m{1.5cm}|  }  {\centering  $L_{74}^{3}(3\CC)$ } &
\multicolumn{3}{  m{4.25cm}||}  {\centering  $(1, 2, 4, 5, 7)E$  } \tabularnewline
\hline
\hline
\multicolumn{3}{||m{1.5cm}|   }  {\centering $\dim$ $\Der$ \\               $11 $    } &
\multicolumn{6}{  m{3.75cm}  |   }  {\centering $\dim$ {Derived Series}  \\     $(7, 4, 0)$    } &
\multicolumn{3}{  m{4.25cm}|| }  {\centering $\dim$ {Desc. C. Series} \\     $(7, 4, 3, 2, 1, 0)$    } \tabularnewline
\hline
\end{tabular}
\end{center}

\

\noindent $\mu:=\left\{\begin{array}{l}
[e_1,e_2]=e_4  ,   [e_1,e_4]=e_5,   [e_1,e_5]=e_6,   [e_1,e_6]=e_7, {[e_2,e_3]}=e_6,  \\
{[e_2,e_4]}=e_6,   [e_2,e_5]=e_7,   [e_3,e_4]=-e_7
\end{array}\right.$

\

\noindent Let $g \in \mathrm{GL}_7(\RR)$ and $g^{-1}$ its inverse

\noindent \begin{tabular}{cc}
  $g=$ & $g^{-1}=$ \\
 \scalebox{0.59}{$\left( \begin {array}{ccccccc}
 1&0&0&0&0&0&0\\
 0&{\frac {\sqrt{2170}}{155}}&0&0&0&0&0\\
0&0&-{\frac {\sqrt{3990}}{1767}} & {\frac {7\sqrt{3990}}{8835}}&0&0&0\\
0&0&{\frac {\sqrt{42}}{93}}&{\frac {\sqrt{42}}{93}}&0&0&0\\
0&0&0&0&{\frac{28\sqrt{5890}}{91295}}&0&0\\
0&0&0&0&0&{\frac {56\sqrt {95}}{91295}}&0\\
0&0&0&0&0&0&{\frac {28\sqrt{23870}}{2830145}}
\end {array} \right)$} &
\scalebox{0.54}{$\left( \begin {array}{ccccccc}
1&0&0&0&0&0&0\\
0&\frac{\sqrt{2170}}{14}&0&0&0&0&0\\
0&0&-{\frac {31\sqrt{3990}}{168}}&{\frac {31\sqrt {42}}{24}}&0&0&0\\
0&0& {\frac {31\sqrt{3990}}{168}}&{\frac {155\sqrt {42}}{168}}&0&0&0\\
0&0&0&0&{\frac {31\sqrt{5890}}{56}}&0&0\\
0&0&0&0&0&{\frac {961\sqrt {95}}{56}}&0\\
0&0&0&0&0&0&{\frac {18259\sqrt{23870}}{4312}}
\end {array} \right)$}\end{tabular}

\

\noindent The action of $g$ over $\mu$ gives a isomorphic Lie algebra law, $g\cdot \mu=\widetilde{\mu}$

\

\noindent $\widetilde{\mu}:=\left\{
\begin{array}{l}
[e_1, e_2] = \frac{7\sqrt{1767}}{1767} e_3  + \frac{\sqrt{465}}{93} e_4, [e_1, e_3] = \frac{\sqrt{651}}{93}e_5,  [e_1, e_4] = \frac{\sqrt{61845}}{1767} e_5,\\
{[e_1, e_5]} = \frac{\sqrt{62} }{31}e_6, [e_1, e_6] = \frac{\sqrt{90706}}{1178} e_7, [e_2, e_4] = \frac{4\sqrt{1767}}{589} e_6,\\
{[e_2, e_5]} = \frac{\sqrt{64790}}{1178} e_7, {[e_3, e_4]} = \frac{\sqrt{90706}}{1178} e_7 \end{array}\right.$

\

\noindent By straight calculation, it is easy to see that the moment map of $\widetilde{\mu}$ with respect to the ordered basis $e_1,...,e_7$ is:

\noindent \begin{eqnarray*}
  \mm(\widetilde{\mu}) &=& \diag(-\frac{19}{31}, -\frac{13}{31}, -\frac{7}{31}, -\frac{7}{31}, -\frac{1}{31}, \frac{5}{31}, \frac{11}{31})\\
                       &=& -\frac{25}{31} \mathrm{Id} + \frac{6}{31} \diag( 1,2,3,3,4,5,6)
\end{eqnarray*}

\

\noindent Since $\diag(1,2,3,3,4,5,6) $ is a derivation of the Lie algebra $(\RR^7, \widetilde{\mu})$ then by \cite[Theorem 2.2]{FERNANDEZ-CULMA1}, it is an Einstein Nilradical. \newline
\noindent $||\mathcal{S}_{\beta}||^2=\frac{25}{31}\thickapprox 0.806$

\

\begin{center}
\begin{tabular}{cccccccccccc}
\hline
\multicolumn{3}{||m{1.5cm}|  }  {\centering  $1.12$  } &
\multicolumn{3}{  m{1.5cm}|  }  {\centering  $\mathfrak{n}_{7,49}^{}$ } &
\multicolumn{3}{  m{1.5cm}|  }  {\centering  $L_{75}^{44}$ } &
\multicolumn{3}{  m{4.25cm}||}  {\centering  $(1,3,4,5,7)D$  } \tabularnewline
\hline
\hline
\multicolumn{3}{||m{1.5cm}|   }  {\centering $\dim$ $\Der$ \\               $12$    } &
\multicolumn{6}{  m{3.75cm}  |   }  {\centering $\dim$ {Derived Series}  \\     $(7, 4, 0)$    } &
\multicolumn{3}{  m{4.25cm}|| }  {\centering $\dim$ {Desc. C. Series} \\     $(7, 4, 3, 2, 1, 0)$    } \tabularnewline
\hline
\end{tabular}
\end{center}

\

\noindent  $\mu:=\left\{\begin{array}{l}
[e_1,e_2]=e_4,  [e_1,e_4]=e_5,  [e_1,e_5]=e_6, [e_1,e_6]=e_7, {[e_2,e_3]}=e_7,  \\
{[e_2,e_4]}=e_6,  [e_2,e_5]=e_7
\end{array}
\right.$

\

\noindent Let $g \in \mathrm{GL}_7(\RR)$ and $g^{-1}$ its inverse

\noindent	\begin{tabular}{cc}
  $g=$ & $g^{-1}=$ \\
\scalebox{0.56}{$\left( \begin {array}{ccccccc} 1&0&0&0&0&0&0\\
0&{\frac {\sqrt {14322}}{434}}&0&0&0&0&0\\
0&0&{\frac {33\sqrt {1085}}{53816}}&0&{\frac {33\sqrt {1085}}{53816}}&0&0\\
0&0&0&{\frac {3\sqrt {385}}{868}}&0&0&0\\
0&0&-{\frac {3\sqrt {23870}}{26908}}&0&{\frac {3\sqrt {23870}}{53816}}&0&0\\
0&0&0&0&0&{\frac {99\sqrt {35}}{107632}}&0\\
0&0&0&0&0&0&{\frac {99\sqrt {11935}}{6673184}}
\end {array} \right)$}
 &
\scalebox{0.5}{$\left( \begin {array}{ccccccc}
1&0&0&0&0&0&0\\
0&\frac{\sqrt {14322}}{33}&0&0&0&0&0\\
0&0&{\frac {248\sqrt {1085}}{495}}&0&-{\frac {124\sqrt {23870}}{495}}&0&0\\
0&0&0&{\frac {124\sqrt {385}}{165}}&0&0&0\\
0&0&{\frac {496\sqrt{1085}}{495}}&0&{\frac {124\sqrt {23870}}{495}}&0&0\\
0&0&0&0&0&{\frac {15376\sqrt {35}}{495}}&0\\
0&0&0&0&0&0&{\frac {30752\sqrt {11935}}{5445}}
\end {array} \right)$}\end{tabular}

\

\noindent The action of $g$ over $\mu$ gives a isomorphic Lie algebra law, $g\cdot \mu=\widetilde{\mu}$

\

\noindent $\widetilde{\mu}:=\left\{\begin{array}{l}
[e_1, e_2] = \frac{\sqrt{930}}{124}e_4,    [e_1, e_3] = \frac{\sqrt{31}}{31}e_6,  [e_1, e_4] = \frac{\sqrt{341}}{62}e_3 + \frac{\sqrt{62}}{62}e_5,\\
{[e_1, e_5]} = \frac{\sqrt{682}}{124}e_6,  [e_1, e_6] = \frac{\sqrt{341}}{62}e_7,  [e_2, e_3] = \frac{\sqrt{1302}}{124}e_7,
{[e_2, e_4]} = \frac{\sqrt{1302}}{124}e_6
\end{array}\right.$

\

\noindent By straight calculation, it is easy to see that the moment map of $\widetilde{\mu}$ with respect to the ordered basis $e_1,...,e_7$ is:

\noindent \begin{eqnarray*}
  \mm(\widetilde{\mu}) &=& \diag (-{\frac {41}{62}},-{\frac {57}{124}},-{\frac {7}{124}},-{\frac {8}{31}},-{\frac {7}{124}},{\frac {9}{62}},{\frac {43}{124}})\\
                       &=& -\frac{107}{124}\mathrm{Id} + \frac{25}{124}\diag( 1,2,4,3,4,5,6)
\end{eqnarray*}

\

\noindent Since $( 1,2,4,3,4,5,6) $ is a derivation of the Lie algebra $(\RR^7, \widetilde{\mu})$ then by \cite[Theorem 2.2]{FERNANDEZ-CULMA1}, it is an Einstein Nilradical.\newline
\noindent $||\mathcal{S}_{\beta}||^2=\frac{107}{124}\thickapprox 0.863$

\

\begin{center}
\begin{tabular}{cccccccccccc}
\hline
\multicolumn{3}{||m{1.5cm}|  }  {\centering  $1.13$  } &
\multicolumn{3}{  m{1.5cm}|  }  {\centering  $\mathfrak{n}_{7,30}^{}$ } &
\multicolumn{3}{  m{1.5cm}|  }  {\centering  $L_{75}^{37}$ } &
\multicolumn{3}{  m{4.25cm}||}  {\centering  $(2,3,4,5,7)E$  } \tabularnewline
\hline
\hline
\multicolumn{3}{||m{1.5cm}|   }  {\centering $\dim$ $\Der$ \\               $12$    } &
\multicolumn{6}{  m{3.75cm}  |   }  {\centering $\dim$ {Derived Series}  \\     $(7, 5, 0)$    } &
\multicolumn{3}{  m{4.25cm}|| }  {\centering $\dim$ {Desc. C. Series} \\     $(7, 5, 4, 2, 1, 0)$    } \tabularnewline
\hline
\end{tabular}
\end{center}

\

\noindent $\mu:=\left\{\begin{array}{l}
[e_1,e_2]=e_3, [e_1,e_3]=e_4, [e_1,e_4]=e_6,  [e_1,e_5]=e_7, {[e_1,e_6]}=e_7,  \\
{[e_2,e_3]}=e_5,  [e_2,e_4]=e_7
\end{array}\right.$

\

\noindent Let $g \in \mathrm{GL}_7(\RR)$ and $g^{-1}$ its inverse

\noindent	\begin{tabular}{cc}
  $g=$ & $g^{-1}=$ \\
\scalebox{0.68}{$\left( \begin {array}{ccccccc}
1&0&0&0&0&0&0\\
0&\frac{\sqrt {17}}{17}&0&0&0&0&0\\
0&0&{\frac {3\sqrt {2}}{68}}&0&0&0&0\\
0&0&0&{\frac {3\sqrt {238}}{2312}}&0&0&0\\
0&0&0&0&{\frac {9\sqrt {154}}{25432}} & {\frac {9\sqrt {154}}{25432}}&0\\
0&0&0&0&{\frac {3\sqrt {77}}{12716}}  & -{\frac {21\sqrt {77}}{50864}}&0\\
0&0&0&0&0&0&{\frac {9\sqrt {119}}{78608}}
\end {array} \right)$}
 &
\scalebox{0.59}{$\left( \begin {array}{ccccccc}
1&0&0&0&0&0&0\\
0&\sqrt {17}&0&0&0&0&0\\
0&0&{\frac {34\sqrt {2}}{3}}&0&0&0&0\\
0&0&0&{\frac {68\sqrt {238}}{21}}&0&0&0\\
0&0&0&0&{\frac {1156\sqrt {154}}{99}}  &  {\frac {4624\sqrt {77}}{231}}&0\\
0&0&0&0&{\frac {4624\sqrt {154}}{693}} & -{\frac {4624\sqrt {77}}{231}}&0\\
0&0&0&0&0&0&{\frac {4624\sqrt {119}}{63}}
\end {array} \right)$}
		\end{tabular}

\

\noindent The action of $g$ over $\mu$ gives a isomorphic Lie algebra law, $g\cdot \mu=\widetilde{\mu}$

\

\noindent $\widetilde{\mu}:=\left\{
\begin{array}{l}
[e_1,e_2]={\frac {3\sqrt {34}}{68}}e_3, [e_1,e_3]=\frac{\sqrt {119}}{34} e_4, [e_1,e_4]={\frac {3\sqrt {187}}{187}} e_5 - {\frac {7\sqrt {374}}{748}}e_6,\\
{[e_1,e_5]}={\frac {\sqrt {374}}{68}}e_7, [e_2,e_3]={\frac {3\sqrt {1309}}{374}}e_5 + {\frac {\sqrt {2618}}{374}}e_6,  [e_2,e_4]={\frac {3\sqrt {34}}{68}}e_7
\end{array}
\right.$

\

\noindent By straight calculation, it is easy to see that the moment map of $\widetilde{\mu}$ with respect to the ordered basis $e_1,...,e_7$ is:

\noindent \begin{eqnarray*}
  \mm(\widetilde{\mu}) &=& \diag(-{\frac {45}{68}},-{\frac {8}{17}},-{\frac {19}{68}},-{\frac {3}{34}}
,{\frac {7}{68}},{\frac {7}{68}},{\frac {5}{17}})\\
                       &=& -\frac{29}{34}\mathrm{Id} + \frac{13}{68} \diag(1,2,3,4,5,5,6)
\end{eqnarray*}

\

\noindent Since $\diag(1,2,3,4,5,5,6)$ is a derivation of the Lie algebra $(\RR^7, \widetilde{\mu})$ then by \cite[Theorem 2.2]{FERNANDEZ-CULMA1}, it is an Einstein Nilradical.\newline
\noindent $||\mathcal{S}_{\beta}||^2=\frac{29}{34}\thickapprox 0.853$

\

\begin{center}
\begin{tabular}{cccccccccccc}
\hline
\multicolumn{3}{||m{1.5cm}|  }  {\centering  $1.14$  } &
\multicolumn{3}{  m{1.5cm}|  }  {\centering  $\mathfrak{n}_{7,27}^{}$ } &
\multicolumn{3}{  m{1.5cm}|  }  {\centering  $L_{74}^{32}(\mathfrak{n}_3)$ } &
\multicolumn{3}{  m{4.25cm}||}  {\centering  $(2,3,4,5,7)F$  } \tabularnewline
\hline
\hline
\multicolumn{3}{||m{1.5cm}|   }  {\centering $\dim$ $\Der$ \\               $11 $    } &
\multicolumn{6}{  m{3.75cm}  |   }  {\centering $\dim$ {Derived Series}  \\     $(7, 5, 1, 0)$    } &
\multicolumn{3}{  m{4.25cm}|| }  {\centering $\dim$ {Desc. C. Series} \\     $(7, 5, 4, 2, 1, 0)$    } \tabularnewline
\hline
\end{tabular}
\end{center}

\

\noindent  $\mu:=\left\{
\begin{array}{l}
[e_1,e_2]=e_3, [e_1,e_3]=e_4, [e_1,e_4]=e_6, [e_2,e_3]=e_5, {[e_2,e_5]}=-e_7, \\
{[e_2,e_6]}=-e_7,  [e_3,e_4]=e_7 \\
\end{array}\right.$

\

\noindent Let $g \in \mathrm{GL}_7(\RR)$ and $g^{-1}$ its inverse

\noindent	\begin{tabular}{cc}
  $g=$ & $g^{-1}=$ \\
\scalebox{0.41}{$\left( \begin {array}{ccccccc}
1&0&0&0&0&0&0\\
0&{\frac {3\sqrt {293683}}{2378}}&0&0&0&0&0\\
0&0&-{\frac {57\sqrt {299}}{4756}}&0&0&0&0\\
0&0&0&-{\frac {741\sqrt {667}}{275848}}&0&0&0\\
0&0&0&0&-{\frac {51129\sqrt {908314}}{2397670816}}&-{\frac {51129\sqrt {908314}}{2397670816}}&0\\
0&0&0&0&{\frac {126711\sqrt {99958}}{2397670816}}&-{\frac {741\sqrt {99958}}{58479776}}&0\\
0&0&0&0&0&0&{\frac {2914353\sqrt {5863}}{53789256608}}
\end {array} \right)
$}
 &
\scalebox{0.38}{$\left( \begin {array}{ccccccc}
1&0&0&0&0&0&0\\
0&{\frac {2\sqrt {293683}}{741}}&0&0&0&0&0\\
0&0&-{\frac {4756\sqrt {299}}{17043}}&0&0&0&0\\
0&0&0&-{\frac {9512\sqrt {667}}{17043}}&0&0&0\\
0&0&0&0&-{\frac {5654884\sqrt {908314}}{566355933}}& {\frac {137924\sqrt {99958}}{903279}}&0\\
0&0&0&0&-{\frac {137924\sqrt {908314}}{3312023}}&-{\frac {137924\sqrt {99958}}{903279}}&0\\
0&0&0&0&0&0&{\frac {1311933088\sqrt {5863}}{416752479}}
\end {array} \right)$}
		\end{tabular}

\

\noindent The action of $g$ over $\mu$ gives a isomorphic Lie algebra law, $g\cdot \mu=\widetilde{\mu}$

\

\noindent $\widetilde{\mu}:=\left\{\begin{array}{l}
[e_1,e_2]=-\frac{\sqrt {519593}}{2378}e_3,  [e_1,e_3]=\frac {\sqrt {377}}{58}e_4,\\
{[e_1,e_4]}=\frac {3\sqrt {605845438}}{252068}e_5 + \frac {\sqrt {126034}}{6148}e_6,
{[e_2,e_3]}=\frac {\sqrt {777722}}{6148}e_5 - \frac {3 \sqrt {58406}}{6148}e_6,\\
{[e_2,e_5]}=\frac {\sqrt {126034}}{1189}e_7, [e_3,e_4]=\frac {3\sqrt {13079}}{1189}e_7
\end{array}\right.$

\

\noindent By straight calculation, it is easy to see that the moment map of $\widetilde{\mu}$ with respect to the ordered basis $e_1,...,e_7$ is:

\noindent \begin{eqnarray*}
  \mm(\widetilde{\mu}) &=& \diag(-{\frac {17}{29}},-{\frac {25}{58}},-{\frac {8}{29}},-{\frac {7}{58}}
,\frac{1}{29},\frac{1}{29},{\frac {10}{29}})\\
                       &=& -\frac{43}{58}\mathrm{Id} + \frac{9}{58} \diag (1,2,3,4,5,5,7)
\end{eqnarray*}

\

\noindent Since $\diag (1,2,3,4,5,5,7) $ is a derivation of the Lie algebra $(\RR^7, \widetilde{\mu})$ then by \cite[Theorem 2.2]{FERNANDEZ-CULMA1}, it is an Einstein Nilradical.\newline
\noindent $||\mathcal{S}_{\beta}||^2=\frac{43}{58}\thickapprox 0.741$

\

\begin{center}
\begin{tabular}{cccccccccccc}
\hline
\multicolumn{3}{||m{1.5cm}|  }  {\centering  $1.15$  } &
\multicolumn{3}{  m{1.5cm}|  }  {\centering  $\mathfrak{n}_{7,50}^{}$ } &
\multicolumn{3}{  m{1.5cm}|  }  {\centering  $L_{75}^{41}$ } &
\multicolumn{3}{  m{4.25cm}||}  {\centering  $(1,3,4,5,7)B$  } \tabularnewline
\hline
\hline
\multicolumn{3}{||m{1.5cm}|   }  {\centering $\dim$ $\Der$ \\               $13 $    } &
\multicolumn{6}{  m{3.75cm}  |   }  {\centering $\dim$ {Derived Series}  \\     $(7, 4, 0)$    } &
\multicolumn{3}{  m{4.25cm}|| }  {\centering $\dim$ {Desc. C. Series} \\     $(7, 4, 3, 2, 1, 0)$    } \tabularnewline
\hline
\end{tabular}
\end{center}

\

\noindent  $\mu:=\left\{\begin{array}{l}
[e_1,e_2]=e_4,   [e_1,e_4]=e_5,  [e_1,e_5]=e_6,  [e_1,e_6]=e_7, {[e_2,e_3]}=e_7, \\
{[e_2,e_4]}=e_7
\end{array}\right.$

\

\noindent Let $g \in \mathrm{GL}_7(\RR)$ and $g^{-1}$ its inverse

\noindent	\begin{tabular}{cc}
  $g=$ & $g^{-1}=$ \\
\scalebox{0.59}{$\left( \begin {array}{ccccccc}
1&0&0&0&0&0&0\\0&\frac{\sqrt {11}}{41}&0&0&0&0&0\\0&0&-{\frac {5\sqrt {47355}}{282408}}  & {\frac {2\sqrt {47355}}{35301}}&0&0&0\\
0&0&{\frac{\sqrt {1231230}}{70602}}& {\frac {\sqrt {1231230}}{70602}}&0&0&0\\0&0&0&0&{\frac {\sqrt{1430}}{6724}}&0&0\\
0&0&0&0&0&{\frac {11\sqrt {2665}}{275684}}&0\\0&0&0&0&0&0&{\frac {11\sqrt {130}}{275684}}\end {array} \right)$}
 &
\scalebox{0.47}{$\left( \begin {array}{ccccccc}
1&0&0&0&0&0&0\\0&{\frac {41\sqrt {11}}{11}}&0&0&0&0&0\\0&0&-{\frac{328\sqrt{47355}}{1155}}& {\frac {656\sqrt{1231230}}{15015}} &0&0&0\\
0&0&{\frac {328\sqrt{47355}}{1155}}& {\frac {41 \sqrt{1231230}}{3003}}&0&0&0\\0&0&0&0&{\frac {3362\sqrt{1430}}{715}} &0&0\\
0&0&0&0&0&{\frac {6724\sqrt{2665}}{715}} &0\\0&0&0&0&0&0&{\frac {137842\sqrt{130}}{715}}\end {array} \right)$}
		\end{tabular}

\

\noindent The action of $g$ over $\mu$ gives a isomorphic Lie algebra law, $g\cdot \mu=\widetilde{\mu}$

\

\noindent $\widetilde{\mu}:=\left\{
\begin{array}{l}
[e_1,e_2]=\frac {2\sqrt {4305}}{861} e_3 + \frac {\sqrt {111930}}{1722}e_4,  [e_1,e_3]=\frac {2\sqrt {22386}}{861}e_5,  [e_1,e_4]=\frac {5\sqrt {861}}{1722}e_5,\\
{[e_1,e_5]}=\frac {\sqrt {902}}{82}e_6, [e_1,e_6]=\frac{\sqrt {82}}{41}e_7, [e_2,e_4]=\frac {\sqrt {861}}{82}e_7
\end{array}\right.$

\

\noindent By straight calculation, it is easy to see that the moment map of $\widetilde{\mu}$ with respect to the ordered basis $e_1,...,e_7$ is:

\noindent \begin{eqnarray*}
  \mm(\widetilde{\mu}) &=&  \diag (-{\frac {61}{82}},-{\frac {31}{82}},-{\frac {8}{41}},-{\frac {8}{41}},-{\frac {1}{82}},{\frac {7}{41}},{\frac {29}{82}})\\
                       &=& -\frac{38}{41}\mathrm{Id} + \frac{15}{82}\diag (1,3,4,4,5,6,7)
\end{eqnarray*}

\

\noindent Since $ \diag(1,3,4,4,5,6,7) $ is a derivation of the Lie algebra $(\RR^7, \widetilde{\mu})$ then by \cite[Theorem 2.2]{FERNANDEZ-CULMA1}, it is an Einstein Nilradical.\newline
\noindent $||\mathcal{S}_{\beta}||^2=\frac{38}{41}\thickapprox 0.927$

\

\begin{center}
\begin{tabular}{cccccccccccc}
\hline
\multicolumn{3}{||m{1.5cm}|  }  {\centering  $1.16$  } &
\multicolumn{3}{  m{1.5cm}|  }  {\centering          } &
\multicolumn{3}{  m{1.5cm}|  }  {\centering  $L_{75}^{50}$ } &
\multicolumn{3}{  m{4.25cm}||}  {\centering  $(2,4,5,7)D$  } \tabularnewline
\hline
\hline
\multicolumn{3}{||m{1.5cm}|   }  {\centering $\dim$ $\Der$ \\               $ 15$    } &
\multicolumn{6}{  m{3.75cm}  |   }  {\centering $\dim$ {Derived Series}  \\     $(7, 4, 0)$    } &
\multicolumn{3}{  m{4.25cm}|| }  {\centering $\dim$ {Desc. C. Series} \\     $(7, 4, 2, 1, 0)$    } \tabularnewline
\hline
\end{tabular}
\end{center}

\

\noindent  $\mu:=\left\{
\begin{array}{l}
[e_1, e_2] = e_3,    [e_1, e_3] = e_5, [e_1, e_4] = e_6, [e_1, e_5] = e_7, {[e_1, e_6]} = e_7,  \\
{[e_2, e_3]} = e_7
\end{array}
\right.$

\

\noindent Let $g \in \mathrm{GL}_7(\RR)$ and $g^{-1}$ its inverse

\noindent	\begin{tabular}{cc}
  $g=$ & $g^{-1}=$ \\
\scalebox{0.72}{$\left( \begin {array}{ccccccc}
1&0&0&0&0&0&0\\
0&\frac{\sqrt {95}}{38}&0&0&0&0&0\\
0&0&{\frac {\sqrt{30}}{76}} & {\frac {\sqrt{30}}{76}}&0&0&0\\
0&0&{\frac {\sqrt{15}}{76}}  & -{\frac {\sqrt{15}}{76}}&0&0&0\\
0&0&0&0&{\frac {\sqrt{285}}{1444}}  & -{\frac {\sqrt{285}}{1444}}&0\\
0&0&0&0&{\frac {5\sqrt {57}}{1444}} &  {\frac {5\sqrt {57}}{1444}}&0\\
0&0&0&0&0&0&{\frac {15\sqrt {2}}{2888}}
\end {array} \right)$
}
 &
\scalebox{0.6}{$\left( \begin {array}{ccccccc}
1&0&0&0&0&0&0\\
0&\frac{2\sqrt {95}}{5}&0&0&0&0&0\\
0&0&{\frac {19\sqrt{30}}{15}} & {\frac {38\sqrt {15}}{15}}&0&0&0\\
0&0&{\frac {19\sqrt{30}}{15}} & -{\frac {38\sqrt {15}}{15}}&0&0&0\\
0&0&0&0&{\frac {38\sqrt{285}}{15}}  & {\frac {38\sqrt{57}}{15}}&0\\
0&0&0&0&-{\frac {38\sqrt{285}}{15}} & {\frac {38\sqrt{57}}{15}}&0\\
0&0&0&0&0&0&{\frac {1444\sqrt {2}}{15}}
\end {array} \right)$}
		\end{tabular}

\

\noindent The action of $g$ over $\mu$ gives a isomorphic Lie algebra law, $g\cdot \mu=\widetilde{\mu}$

\

\noindent $\widetilde{\mu}:=\left\{
\begin{array}{l}
[e_1,e_2]=\frac{\sqrt{114}}{38} e_3 + \frac{\sqrt {57}}{38}e_4, [e_1,e_3]=\frac{\sqrt{190}}{38}e_6, [e_1,e_4]=\frac{\sqrt {19}}{19}e_5,\\
{[e_1,e_6]}=\frac{\sqrt{114}}{38}e_7, [e_2,e_3]=\frac{\sqrt{57}}{38}e_7, [e_2,e_4]=\frac{\sqrt{114}}{38}e_7
\end{array}\right.$

\

\noindent By straight calculation, it is easy to see that the moment map of $\widetilde{\mu}$ with respect to the ordered basis $e_1,...,e_7$ is:

\noindent \begin{eqnarray*}
  \mm(\widetilde{\mu}) &=& \diag(-{\frac {29}{38}},-{\frac {9}{19}},-{\frac {7}{38}},-{\frac {7}{38}},
\frac{2}{19},\frac{2}{19},{\frac {15}{38}})\\
                       &=&-\frac{20}{19}\mathrm{Id} + \frac{11}{38} \diag(1,2,3,3,4,4,5)
\end{eqnarray*}

\

\noindent Since $\diag(1,2,3,3,4,4,5)$ is a derivation of the Lie algebra $(\RR^7, \widetilde{\mu})$ then by \cite[Theorem 2.2]{FERNANDEZ-CULMA1}, it is an Einstein Nilradical.\newline
\noindent $||\mathcal{S}_{\beta}||^2=\frac{20}{19}\thickapprox 1.05$

\

\begin{center}
\begin{tabular}{cccccccccccc}
\hline
\multicolumn{3}{||m{1.5cm}|  }  {\centering  $1.17 $  } &
\multicolumn{3}{  m{1.5cm}|  }  {\centering  $\mathfrak{n}_{7,18}^{}$ } &
\multicolumn{3}{  m{1.5cm}|  }  {\centering  $L_{74}^{18}(\mathfrak{n}_3)$ } &
\multicolumn{3}{  m{4.25cm}||}  {\centering  $(1,2,4,5,7)L$  } \tabularnewline
\hline
\hline
\multicolumn{3}{||m{1.5cm}|   }  {\centering $\dim$ $\Der$ \\               $11 $    } &
\multicolumn{6}{  m{3.75cm}  |   }  {\centering $\dim$ {Derived Series}  \\     $(7, 5, 1, 0)$    } &
\multicolumn{3}{  m{4.25cm}|| }  {\centering $\dim$ {Desc. C. Series} \\     $(7, 5, 4, 2, 1, 0)$    } \tabularnewline
\hline
\end{tabular}
\end{center}

\

\noindent  $\mu:=\left\{\begin{array}{l}
[e_1,e_2]=e_3,    [e_1,e_3]=e_4,  [e_1,e_4]=e_6, [e_1,e_6]=e_7,  [e_2,e_3] =e_5, \\
{[e_2,e_5]}=e_6,  [e_2,e_6]=e_7,  [e_3,e_4]=-e_7, [e_3,e_5]=e_7
\end{array}
\right.$

\

\noindent Let $g \in \mathrm{GL}_7(\RR)$ and $g^{-1}$ its inverse

\noindent	\begin{tabular}{cc}
  $g=$ & $g^{-1}=$ \\
\scalebox{0.5}{$\left( \begin {array}{ccccccc}
1&-1&0&0&0&0&0\\
\sqrt {2}&\sqrt {2}&0&0&0&0&0\\
0&0&\frac{\sqrt {611}\sqrt {2}}{47}&0&0&0&0\\
0&0&0&\frac{\sqrt {13}\sqrt {5}}{47} &\frac{\sqrt {13}\sqrt {5}}{47}&0&0\\
0&0&0&\frac{\sqrt {13}\sqrt {5}\sqrt {2}}{47} & -\frac{\sqrt {13}\sqrt {5}\sqrt {2}}{47}&0&0\\
0&0&0&0&0&{\frac {13}{2209}}\sqrt {47}\sqrt {5}\sqrt {2}&0\\
0&0&0&0&0&0&{\frac {65\sqrt{3}}{2209}}
\end {array} \right)$}
 &
\scalebox{0.52}{$\left( \begin {array}{ccccccc}
\frac{1}{2} & \frac{\sqrt{2}}{4}&0&0&0&0&0\\
 -\frac{1}{2}& \frac{\sqrt{2}}{4} &0&0&0&0&0\\
 0&0&\frac{\sqrt {611} \sqrt {2}}{26} &0&0&0&0\\
 0&0&0&\frac{47\sqrt {13}\sqrt {5}}{130} & \frac {47\sqrt {13}\sqrt {5} \sqrt {2}}{260} &0&0\\
  0&0&0& \frac {47\sqrt {13} \sqrt {5}}{130}  &- \frac {47\sqrt {13} \sqrt {5} \sqrt {2}}{260} &0&0\\
0&0&0&0&0& \frac {47\sqrt {47} \sqrt {5}\sqrt {2}}{130} &0\\
0&0&0&0&0&0& \frac {2209\sqrt {3}}{195}
\end {array} \right)$}
		\end{tabular}

\

\noindent The action of $g$ over $\mu$ gives a isomorphic Lie algebra law, $g\cdot \mu=\widetilde{\mu}$

\

\noindent $\widetilde{\mu}:=\left\{\begin{array}{l}
[e_1, e_2] = \frac{\sqrt{611}}{94}e_3, [e_1,e_3] = \frac{\sqrt{235}}{47}e_5, [e_1, e_5] = \frac{\sqrt{611}}{94}e_6,  [e_2, e_3] = \frac{\sqrt{235}}{94}e_4, \\
{[e_2, e_4]} = \frac{\sqrt{611}}{94}e_6,  [e_2, e_6] = \frac{\sqrt{705}}{94}e_7, [e_3, e_5] = -\frac{\sqrt{705}}{94}e_7
\end{array}\right.$

\

\noindent By straight calculation, it is easy to see that the moment map of $\widetilde{\mu}$ with respect to the ordered basis $e_1,...,e_7$ is:

\noindent \begin{eqnarray*}
  \mm(\widetilde{\mu}) &=& \diag\left(-{\frac {23}{47}},-{\frac {23}{47}},-{\frac {27}{94}},-{\frac {4}{47}
},-{\frac {4}{47}},{\frac {11}{94}},{\frac {15}{47}}\right)\\
                       &=&-\frac{65}{94}\mathrm{Id} + \frac{19}{94}\diag(1,1,2,3,3,4,5)
\end{eqnarray*}

\

\noindent Since $\diag(1,1,2,3,3,4,5)$ is a derivation of the Lie algebra $(\RR^7, \widetilde{\mu})$ then by \cite[Theorem 2.2]{FERNANDEZ-CULMA1}, it is an Einstein Nilradical.\newline
\noindent $||\mathcal{S}_{\beta}||^2=\frac{65}{94}\thickapprox 0.692$

\

\textbf{Another way to prove that $1.17$ is an Einstein Nilradical}

\

\noindent Let $g \in \mathrm{GL}_7(\RR)$ and $g^{-1}$ its inverse

\noindent	\begin{center} \begin{tabular}{cc}
  $g=$ & $g^{-1}=$ \\
\scalebox{1}{ $\frac{1}{8} \left( \begin{array}{ccccccc} 4&4&0&0&0&0&0\\ \noalign{\medskip}-4&4
&0&0&0&0&0\\ \noalign{\medskip}0&0&4&0&0&0&0\\ \noalign{\medskip}0&0&0
&2&2&0&0\\ \noalign{\medskip}0&0&0&-2&2&0&0\\ \noalign{\medskip}0&0&0&0
&0&2&0\\ \noalign{\medskip}0&0&0&0&0&0&1\end{array} \right)$ }
 &
\scalebox{1}{$\left( \begin {array}{ccccccc} 1&-1&0&0&0&0&0\\ \noalign{\medskip}1&1
&0&0&0&0&0\\ \noalign{\medskip}0&0&2&0&0&0&0\\ \noalign{\medskip}0&0&0
&2&-2&0&0\\ \noalign{\medskip}0&0&0&2&2&0&0\\ \noalign{\medskip}0&0&0&0
&0&4&0\\ \noalign{\medskip}0&0&0&0&0&0&8\end{array} \right) $}
		\end{tabular} \end{center}

\

\noindent The action of $g$ over $\mu$ gives a isomorphic Lie algebra law, $g\cdot \mu=\widetilde{\mu}$

\

\noindent $\widetilde{\mu}:=\left\{\begin{array}{l}
[e_1, e_2] = e_3, [e_1, e_3] = e_4, [e_1, e_4] = e_6, [e_1, e_6] = e_7, [e_2, e_3] = e_5, \\
{[e_2, e_5]} = e_6, [e_3, e_5] = e_7
\end{array}\right.$

\

\noindent $\{e_1...e_7\}$ is a nice basis of $(\RR^n , \widetilde{\mu} )$

$$U=\tiny{ \left[ \begin {array}{ccccccc} 3&0&1&1&0&1&-1\\ \noalign{\medskip}0&3
&0&1&1&0&1\\ \noalign{\medskip}1&0&3&0&0&1&0\\ \noalign{\medskip}1&1&0
&3&0&-1&1\\ \noalign{\medskip}0&1&0&0&3&0&0\\ \noalign{\medskip}1&0&1&
-1&0&3&1\\ \noalign{\medskip}-1&1&0&1&0&1&3\end {array} \right]}$$

\noindent A solution to $U_{\widetilde{\mu}}x=[1]$: $x=\frac{ 1}{65 }(13,5,13,15,20,13,15 )^T$. It is an Einstein nilradical.

\

\begin{center}
\begin{tabular}{cccccccccccc}
\hline
\multicolumn{3}{||m{1.5cm}|  }  {\centering  $1.18 $  } &
\multicolumn{3}{  m{1.5cm}|  }  {\centering  $\mathfrak{n}_{7,57}^{}$ } &
\multicolumn{3}{  m{1.5cm}|  }  {\centering  } &
\multicolumn{3}{  m{4.25cm}||}  {\centering  $(2,4,5,7)J$  } \tabularnewline
\hline
\hline
\multicolumn{3}{||m{1.5cm}|   }  {\centering $\dim$ $\Der$ \\               $13 $    } &
\multicolumn{6}{  m{3.75cm}  |   }  {\centering $\dim$ {Derived Series}  \\     $(7, 4, 0)$    } &
\multicolumn{3}{  m{4.25cm}|| }  {\centering $\dim$ {Desc. C. Series} \\     $(7, 4, 3, 1, 0)$    } \tabularnewline
\hline
\end{tabular}
\end{center}

\

\noindent  $\mu:=\left\{\begin{array}{l}
[e_1,e_2]=e_4, [e_1,e_4]=e_5, [e_1,e_5]=e_7, [e_2,e_3]=e_6+e_7, {[e_2,e_4]=e_6}
\end{array}\right.$

\

\noindent Let $g \in \mathrm{GL}_7(\RR)$ and $g^{-1}$ its inverse

\noindent	\begin{tabular}{cc}
  $g=$ & $g^{-1}=$ \\
\scalebox{0.4}{$\left( \begin{array}{ccccccc}
 2\sqrt {94}&0&0&0&0&0&0\\
0 &-12\sqrt {1551}&0&0&0&0&0\\
0 & 0 & 18 \sqrt {1034} & -54\sqrt{1034}&0&0&0\\
0 & 0 &-198\sqrt {94}   &-198\sqrt {94} &0&0&0\\
0 & 0 &                 0               &0&-2376\sqrt {47}&0&0\\
0 & 0 &                 0               &0&0&2376\sqrt {20163} &0\\
0 & 0 &                 0               &0&0&7128\sqrt {141} &-9504 \sqrt {141}\end{array} \right)$}
 &
\scalebox{0.52}{$\left( \begin{array}{ccccccc}
{\frac {\sqrt {94}}{188}}&0&0&0&0&0&0\\
0&-{\frac {\sqrt {1551}}{18612}}&0&0&0&0&0\\
0&0&{\frac {\sqrt {1034}}{74448}}&-{\frac {\sqrt {94}}{24816}}&0&0&0\\
0&0&-{\frac {\sqrt {1034}}{74448}}&-{\frac {\sqrt {94}}{74448}}&0&0&0\\
0&0&0&0&-{\frac {\sqrt {47}}{111672}}&0&0\\
0&0&0&0&0&{\frac {\sqrt {80652}}{95814576}}&0\\
0&0&0&0&0&{\frac {\sqrt {20163}}{63876384}}&-{\frac {\sqrt {141}}{1340064}}
\end{array} \right) $}
		\end{tabular}

\

\noindent The action of $g$ over $\mu$ gives a isomorphic Lie algebra law, $g\cdot \mu=\widetilde{\mu}$

\

\noindent $\widetilde{\mu}:=\left\{\begin{array}{l}
[e_1, e_2] = \frac{3\sqrt{141}}{188}e_3  + \frac{\sqrt{1551}}{188}e_4, [e_1, e_3] = \frac{3\sqrt{517}}{188}e_5,   [e_1, e_4] = \frac{3\sqrt{47}}{188}e_5,\\
{[e_1, e_5]} = \frac{\sqrt{282}}{47}e_7, [e_2, e_3] = \frac{\sqrt{94}}{94}e_7, [e_2, e_4] = \frac{\sqrt{1222}}{94}e_6
\end{array}\right.$

\

\noindent By straight calculation, it is easy to see that the moment map of $\widetilde{\mu}$ with respect to the ordered basis $e_1,...,e_7$ is:

\noindent \begin{eqnarray*}
  \mm(\widetilde{\mu}) &=& \diag\left(-\frac{33}{47}, -\frac{43}{94}, -\frac{10}{47}, -\frac{10}{47}, \frac{3}{94}, \frac{13}{47}, \frac{13}{47}\right)\\
                       &=& -\frac{89}{94}\mathrm{Id}+\frac{23}{94}\diag(1, 2, 3, 3, 4, 5, 5)
\end{eqnarray*}

\

\noindent Since $\diag(1, 2, 3, 3, 4, 5, 5)$ is a derivation of the Lie algebra $(\RR^7, \widetilde{\mu})$ then by \cite[Theorem 2.2]{FERNANDEZ-CULMA1}, it is an Einstein Nilradical.\newline
\noindent $||\mathcal{S}_{\beta}||^2=\frac{89}{94}\thickapprox 0.947$

\

\begin{center}
\begin{tabular}{cccccccccccc}
\hline
\multicolumn{3}{||m{1.5cm}|  }  {\centering  $1.19 $  } &
\multicolumn{3}{  m{1.5cm}|  }  {\centering  $\mathfrak{n}_{7,81}^{}$ } &
\multicolumn{3}{  m{1.5cm}|  }  {\centering  $L_{74}^{29}(3\CC)$ } &
\multicolumn{3}{  m{4.25cm}||}  {\centering  $(2,4,7)H$  } \tabularnewline
\hline
\hline
\multicolumn{3}{||m{1.5cm}|   }  {\centering $\dim$ $\Der$ \\               $11$    } &
\multicolumn{6}{  m{3.75cm}  |   }  {\centering $\dim$ {Derived Series}  \\     $(7, 4, 0)$    } &
\multicolumn{3}{  m{4.25cm}|| }  {\centering $\dim$ {Desc. C. Series} \\     $(7, 4, 2, 0)$    } \tabularnewline
\hline
\end{tabular}
\end{center}

\

\noindent  $ [e_1,e_2]=e_4,[e_1,e_3]=e_5,[e_1,e_4]=e_7,[e_1,e_5]=e_6,[e_2,e_4]=e_6,[e_3,e_5]=e_7 $

\

\noindent  $\{e_1...e_7\}$ is a nice basis.

$$U=\tiny{\left[ \begin {array}{cccccc} 3&1&0&1&0&0\\ \noalign{\medskip}1&3&1&0
&0&0\\ \noalign{\medskip}0&1&3&1&1&1\\ \noalign{\medskip}1&0&1&3&1&1
\\ \noalign{\medskip}0&0&1&1&3&0\\ \noalign{\medskip}0&0&1&1&0&3
\end {array} \right]
}$$

\noindent General solution to $Ux=[1]$: $x=\frac{1}{29}(7,7,1,1,9,9 )^T$. It is an Einstein nilradical.

\noindent Pre-Einstein Derivation: $\frac{13}{29} \diag(1, 1, 1, 2, 2, 3, 3) $. \newline
$||\mathcal{S}_{\beta}||^2=\frac{29}{34} \thickapprox 0.8529411765$

\

\begin{center}
\begin{tabular}{cccccccccccc}
\hline
\multicolumn{3}{||m{1.5cm}|  }  {\centering  $1.20 $  } &
\multicolumn{3}{  m{1.5cm}|  }  {\centering  $\mathfrak{n}_{7,35}^{}$ } &
\multicolumn{3}{  m{1.5cm}|  }  {\centering  $L_{74}^{23}(\mathfrak{n}_3)$ } &
\multicolumn{3}{  m{4.25cm}||}  {\centering  $(1,2,4,5,7)D$  } \tabularnewline
\hline
\hline
\multicolumn{3}{||m{1.5cm}|   }  {\centering $\dim$ $\Der$ \\               $11 $    } &
\multicolumn{6}{  m{3.75cm}  |   }  {\centering $\dim$ {Derived Series}  \\     $(7, 4, 1, 0)$    } &
\multicolumn{3}{  m{4.25cm}|| }  {\centering $\dim$ {Desc. C. Series} \\     $(7, 4, 3, 2, 1, 0)$    } \tabularnewline
\hline
\end{tabular}
\end{center}

\

\noindent  $ [e_1,e_2]=e_3,[e_1,e_5]=e_6,[e_1,e_6]=e_7,[e_2,e_3]=e_4,[e_2,e_4]=e_6,[e_2,e_5]=e_7,[e_3,e_4]=e_7$

\

\noindent  $\{e_1...e_7\}$ is a nice basis.

$$U=\tiny{  \left[ \begin {array}{ccccccc} 3&1&1&0&1&1&-1\\ \noalign{\medskip}1&3
&0&0&1&1&0\\ \noalign{\medskip}1&0&3&0&-1&1&1\\ \noalign{\medskip}0&0&0
&3&0&1&0\\ \noalign{\medskip}1&1&-1&0&3&1&1\\ \noalign{\medskip}1&1&1&
1&1&3&1\\ \noalign{\medskip}-1&0&1&0&1&1&3\end {array} \right]}$$

\noindent General solution to $Ux=[1]$: $x= (\frac{24}{47}-t,\frac{12}{47},\frac{3}{47}+t,\frac{20}{47},t,-\frac{13}{47},\frac{27}{47}-t)^T$. It is not an Einstein nilradical.

\

\noindent Pre-Einstein Derivation: $ \frac{8}{47} \diag(1, 2, 3, 5, 6, 7, 8)$

\
\begin{center}
\begin{tabular}{cccccccccccc}
\hline
\multicolumn{3}{||m{1.5cm}|  }  {\centering  $1.21$  } &
\multicolumn{3}{  m{1.5cm}|  }  {\centering  $\mathfrak{n}_{7,36}^{}$ } &
\multicolumn{3}{  m{1.5cm}|  }  {\centering  $L_{74}^{15}(\mathfrak{n}_3)$ } &
\multicolumn{3}{  m{4.25cm}||}  {\centering  $(1,2,4,5,7)F$  } \tabularnewline
\hline
\hline
\multicolumn{3}{||m{1.5cm}|   }  {\centering $\dim$ $\Der$ \\               $11 $    } &
\multicolumn{6}{  m{3.75cm}  |   }  {\centering $\dim$ {Derived Series}  \\     $(7, 4, 1, 0)$    } &
\multicolumn{3}{  m{4.25cm}|| }  {\centering $\dim$ {Desc. C. Series} \\     $(7, 4, 3, 2, 1, 0)$    } \tabularnewline
\hline
\end{tabular}
\end{center}

\

\noindent  $\mu:=\left\{\begin{array}{l}
[e_1,e_2]=e_4, [e_1,e_4]=e_5, [e_1,e_5]=e_6, [e_2,e_3]=e_6, [e_2,e_4]=e_6, \\
{[e_2,e_6]}=e_7, [e_4,e_5]=-e_7
\end{array}\right.$

\

\noindent It is not an Einstein nilradical by \cite[Theorem 2.6]{FERNANDEZ-CULMA1}. In fact,

\

\noindent  Einstein Pre-Derivation $\phi=\frac{25}{113}\diag (1, 2, 3, 3, 4, 5, 7) $. \newline
Let $X=\diag (-4, 23, -28, 10, -1, -8, 8) $. $X\in \ggo_{\phi}$ therefore $g_{t}=\exp(tX) \in \Grp_{\phi}$.

\

\noindent  $g_t\cdot\mu=\left\{\begin{array}{l}
[e_1, e_2] = \euler^{-9t}e_4,    [e_1, e_4] = \euler^{-7t}e_5,  [e_1, e_5] = \euler^{-3t}e_6,  [e_2, e_3] = \euler^{-3t}e_6,\\
{[e_2, e_4]} = \euler^{-41t}e_6, [e_2, e_6] = \euler^{-7t}e_7,  [e_4, e_5] = -\euler^{-t}e_7
\end{array}\right.$

\

\noindent The $\Grp_{\phi}-$orbit of $\mu$ is not closed because  $g_t\cdot\mu \rightarrow 0$ as $t \rightarrow \infty$.

\newpage

\section{Rank two}

\

\begin{center}
\begin{tabular}{cccccccccccc}
\hline
\multicolumn{3}{||m{1.5cm}|  }  {\centering  $2.1(i_\lambda)$ \\ $\lambda \neq 0,1$  } &
\multicolumn{3}{  m{1.5cm}|  }  {\centering  $\mathfrak{n}_{7,64}^{\lambda}$ } &
\multicolumn{3}{  m{1.5cm}|  }  {\centering  $L_{74}^{17}(3\CC)$ \\ \cite{CARLES1} } &
\multicolumn{3}{  m{4.25cm}||}  {\centering  $(1,3,5,7)M$ \\ $\xi=\lambda$, $\xi\neq0,1$  } \tabularnewline
\hline
\hline
\multicolumn{3}{||m{1.5cm}|   }  {\centering $\dim$ $\Der$ \\               $14 $    } &
\multicolumn{6}{  m{3.75cm}  |   }  {\centering $\dim$ {Derived Series}  \\     $(7, 4, 0)$    } &
\multicolumn{3}{  m{4.25cm}|| }  {\centering $\dim$ {Desc. C. Series} \\     $(7, 4, 2, 1, 0)$    } \tabularnewline
\hline
\end{tabular}
\end{center}

\

\noindent  $[e_1,e_2]=e_4,[e_1,e_3]=e_5,[e_1,e_4]=e_6,[e_1,e_6]=e_7,[e_2,e_3]=e_6, [e_2,e_5]=\lambda e_7,[e_3,e_4]=(\lambda-1)e_7$

\

\noindent  $\{e_1...e_7\}$ is a nice basis.

$$U=\tiny{ \left[ \begin {array}{ccccccc} 3&1&0&1&1&1&-1\\ \noalign{\medskip}1&3
&1&1&1&-1&1\\ \noalign{\medskip}0&1&3&0&1&0&1\\ \noalign{\medskip}1&1&0
&3&-1&1&1\\ \noalign{\medskip}1&1&1&-1&3&1&1\\ \noalign{\medskip}1&-1&0
&1&1&3&1\\ \noalign{\medskip}-1&1&1&1&1&1&3\end {array} \right]}$$

\noindent  A solution to $Ux=[1]$: $x=\frac{1}{38}(5,7,8,6,4,9,3)^T$. It is an Einstein nilradical.

\

\noindent  Maximal Torus: $D_1=\diag(1, 0, 2, 1, 3, 2, 3)$, $D_2=\diag(0, 1, 0, 1, 0, 1, 1)$.
Pre-Einstein Derivation: $\frac{2}{19} \diag(3, 5, 6, 8, 9, 11, 14)$ \newline
\noindent $||\mathcal{S}_{\beta}||^2=\frac{19}{21}\thickapprox 0.905$

\

\begin{center}
\begin{tabular}{cccccccccccc}
\hline
\multicolumn{3}{||m{1.5cm}|  }  {\centering  $2.1(i_\lambda) $ \\ $\lambda=0$ } &
\multicolumn{3}{  m{1.5cm}|  }  {\centering  $\mathfrak{n}_{7,64}^{\lambda}$ \\ $\lambda=0$ } &
\multicolumn{3}{  m{1.5cm}|  }  {\centering  $L_{75}^{1}$ } &
\multicolumn{3}{  m{4.25cm}||}  {\centering  $(2,3,5,7)B$  } \tabularnewline
\hline
\hline
\multicolumn{3}{||m{1.5cm}|   }  {\centering $\dim$ $\Der$ \\               $14$    } &
\multicolumn{6}{  m{3.75cm}  |   }  {\centering $\dim$ {Derived Series}  \\     $(7, 4, 0)$    } &
\multicolumn{3}{  m{4.25cm}|| }  {\centering $\dim$ {Desc. C. Series} \\     $(7, 4, 2, 1, 0)$    } \tabularnewline
\hline
\end{tabular}
\end{center}

\

\noindent  $[e_1,e_2]=e_4,[e_1,e_3]=e_5,[e_1,e_4]=e_6,[e_1,e_6]=e_7,[e_2,e_3]=e_6, [e_3,e_4]=-e_7$

\

\noindent  $\{e_1...e_7\}$ is a nice basis.

$$U=\tiny{ \left[ \begin {array}{cccccc} 3&1&0&1&1&-1\\ \noalign{\medskip}1&3&1&
1&1&1\\ \noalign{\medskip}0&1&3&0&1&1\\ \noalign{\medskip}1&1&0&3&-1&1
\\ \noalign{\medskip}1&1&1&-1&3&1\\ \noalign{\medskip}-1&1&1&1&1&3
\end {array} \right] }$$

\noindent  General solution to $Ux=[1]$: $x= (\frac{1}{19}+t,-\frac{1}{19},\frac{4}{19},\frac{9}{19}-t,\frac{8}{19}-t,t)^T$. It is not an Einstein nilradical.

\

\noindent  Maximal Torus: $D_1=\diag(1, 0, 2, 1, 3, 2, 3)$,  $D_2=\diag(0, 1, 0, 1, 0, 1, 1)$.
Pre-Einstein Derivation: $\frac{2}{19}\diag (3, 5, 6, 8, 9, 11, 14)$

\

\begin{center}
\begin{tabular}{cccccccccccc}
\hline
\multicolumn{3}{||m{1.5cm}|  }  {\centering  $2.1(i_\lambda) $ \\ $\lambda=1$ } &
\multicolumn{3}{  m{1.5cm}|  }  {\centering  $\mathfrak{n}_{7,64}^{\lambda}$ \\ $\lambda = 1$} &
\multicolumn{3}{  m{1.5cm}|  }  {\centering  $L_{74}^{17}(3\CC)$ \\ \cite{CARLES1} } &
\multicolumn{3}{  m{4.25cm}||}  {\centering  $(1,3,5,7)M$ \\ $\xi=1$ } \tabularnewline
\hline
\hline
\multicolumn{3}{||m{1.5cm}|   }  {\centering $\dim$ $\Der$ \\               $14 $    } &
\multicolumn{6}{  m{3.75cm}  |   }  {\centering $\dim$ {Derived Series}  \\     $(7, 4, 0)$    } &
\multicolumn{3}{  m{4.25cm}|| }  {\centering $\dim$ {Desc. C. Series} \\     $(7, 4, 2, 1, 0)$    } \tabularnewline
\hline
\end{tabular}
\end{center}

\

\noindent  $[e_1,e_2]=e_4,[e_1,e_3]=e_5,[e_1,e_4]=e_6,[e_1,e_6]=e_7,[e_2,e_3]=e_6, [e_2,e_5]=e_7$

\

\noindent  $\{e_1...e_7\}$ is a nice basis.

$$U=\tiny{\left[ \begin {array}{cccccc} 3&1&0&1&1&1\\ \noalign{\medskip}1&3&1&1
&1&-1\\ \noalign{\medskip}0&1&3&0&1&0\\ \noalign{\medskip}1&1&0&3&-1&1
\\ \noalign{\medskip}1&1&1&-1&3&1\\ \noalign{\medskip}1&-1&0&1&1&3
\end {array} \right]}$$

\noindent  General solution to $Ux=[1]$: $x=(\frac{1}{19}, -\frac{1}{19}+t, \frac{4}{19}, \frac{9}{19}-t, \frac{8}{19}-t, t)$ and a particular solution is $x=\frac{1}{38} (2,7,8,9,7,9 )^T$. It is an Einstein nilradical.

\

\noindent  Maximal Torus: $D_1=\diag(1, 0, 2, 1, 3, 2, 3)$,  $D_2=\diag(0, 1, 0, 1, 0, 1, 1)$.
Pre-Einstein Derivation: $\frac{2}{19}\diag (3, 5, 6, 8, 9, 11, 14)$ \newline
\noindent $||\mathcal{S}_{\beta}||^2=\frac{19}{21}\thickapprox 0.905$

\

\begin{center}
\begin{tabular}{cccccccccccc}
\hline
\multicolumn{3}{||m{1.5cm}|  }  {\centering  $2.1(ii)$  } &
\multicolumn{3}{  m{1.5cm}|  }  {\centering  $\mathfrak{n}_{7,69}^{}$ } &
\multicolumn{3}{  m{1.5cm}|  }  {\centering  $L_{74}^{17}(3\CC)$ \\ $b=0$} &
\multicolumn{3}{  m{4.25cm}||}  {\centering  $(1,3,5,7)D$  } \tabularnewline
\hline
\hline
\multicolumn{3}{||m{1.5cm}|   }  {\centering $\dim$ $\Der$ \\               $14 $    } &
\multicolumn{6}{  m{3.75cm}  |   }  {\centering $\dim$ {Derived Series}  \\     $(7, 4, 0)$    } &
\multicolumn{3}{  m{4.25cm}|| }  {\centering $\dim$ {Desc. C. Series} \\     $(7, 4, 2, 1, 0)$    } \tabularnewline
\hline
\end{tabular}
\end{center}

\

\noindent  $[e_1,e_2]=e_4,[e_1,e_3]=e_5,[e_1,e_4]=e_6,[e_1,e_6]=e_7,[e_2,e_5]=e_7, [e_3,e_4]=e_7$

\

\noindent$\{e_1...e_7\}$ is a nice basis.

$$U=\tiny{\left[ \begin {array}{cccccc} 3&1&0&1&1&-1\\ \noalign{\medskip}1&3&1&
1&-1&1\\ \noalign{\medskip}0&1&3&0&0&1\\ \noalign{\medskip}1&1&0&3&1&1
\\ \noalign{\medskip}1&-1&0&1&3&1\\ \noalign{\medskip}-1&1&1&1&1&3
\end {array} \right]}$$

\noindent A solution to $Ux=[1]$: $x=\frac{1}{19} (4,4,4,1,5,3 )^T$. It is an Einstein nilradical.

\

\noindent Maximal Torus: $D_1=\diag(1, 0, 2, 1, 3, 2, 3)$,  $D_2=\diag(0, 1, 0, 1, 0, 1, 1)$.
Pre-Einstein Derivation: $\frac{2}{19}\diag (3, 5, 6, 8, 9, 11, 14)$ \newline
\noindent $||\mathcal{S}_{\beta}||^2=\frac{19}{21}\thickapprox 0.905$

\

\begin{center}
\begin{tabular}{cccccccccccc}
\hline
\multicolumn{3}{||m{1.5cm}|  }  {\centering  $2.1(iii) $  } &
\multicolumn{3}{  m{1.5cm}|  }  {\centering  $\mathfrak{n}_{7,101}^{}$ } &
\multicolumn{3}{  m{1.5cm}|  }  {\centering  $L_{74}^{1}(3\CC)$ } &
\multicolumn{3}{  m{4.25cm}||}  {\centering  $(1,3,5,7)A$  } \tabularnewline
\hline
\hline
\multicolumn{3}{||m{1.5cm}|   }  {\centering $\dim$ $\Der$ \\               $14 $    } &
\multicolumn{6}{  m{3.75cm}  |   }  {\centering $\dim$ {Derived Series}  \\     $(7, 3, 0)$    } &
\multicolumn{3}{  m{4.25cm}|| }  {\centering $\dim$ {Desc. C. Series} \\     $(7, 3, 2, 1, 0)$    } \tabularnewline
\hline
\end{tabular}
\end{center}

\

\noindent  $[e_1,e_2]=e_4,[e_1,e_4]=e_6,[e_1,e_6]=e_7,[e_2,e_3]=e_6,[e_2,e_5]=e_7, [e_3,e_4]=-e_7$

\

\noindent  $\{e_1...e_7\}$ is a nice basis.

$$U=\tiny{ \left[ \begin {array}{cccccc} 3&0&1&1&1&-1\\ \noalign{\medskip}0&3&0&
1&0&1\\ \noalign{\medskip}1&0&3&-1&1&1\\ \noalign{\medskip}1&1&-1&3&1&
1\\ \noalign{\medskip}1&0&1&1&3&1\\ \noalign{\medskip}-1&1&1&1&1&3
\end {array} \right]}$$

\noindent  General solution to $Ux=[1]$: $x=(\frac{1}{19}+t, \frac{4}{19}, \frac{8}{19}-t, \frac{7}{19}-t, \frac{1}{19}, t)$ and a particular solution is $x=\frac{1}{19} (4,4,5,4,1,3 )^T$. It is an Einstein nilradical.

\

\noindent  Maximal Torus: $D_1=\diag(1, 0, 2, 1, 3, 2, 3)$,  $D_2=\diag(0, 1, 0, 1, 0, 1, 1)$.
Pre-Einstein Derivation: $\frac{2}{19}\diag (3, 5, 6, 8, 9, 11, 14)$ \newline
\noindent $||\mathcal{S}_{\beta}||^2=\frac{19}{21}\thickapprox 0.905$

\

\begin{center}
\begin{tabular}{cccccccccccc}
\hline
\multicolumn{3}{||m{1.5cm}|  }  {\centering  $2.1(iv) $  } &
\multicolumn{3}{  m{1.5cm}|  }  {\centering   } &
\multicolumn{3}{  m{1.5cm}|  }  {\centering  $L_{74}^{17}(3\CC)$ \\ $b=-1$} &
\multicolumn{3}{  m{4.25cm}||}  {\centering  $(1,3,7)D$  } \tabularnewline
\hline
\hline
\multicolumn{3}{||m{1.5cm}|   }  {\centering $\dim$ $\Der$ \\               $ 14 $    } &
\multicolumn{6}{  m{3.75cm}  |   }  {\centering $\dim$ {Derived Series}  \\     $(7, 3, 0)$    } &
\multicolumn{3}{  m{4.25cm}|| }  {\centering $\dim$ {Desc. C. Series} \\     $(7, 3, 1, 0)$    } \tabularnewline
\hline
\end{tabular}
\end{center}

\

\noindent  $[e_1,e_3]=e_5,[e_1,e_4]=e_6,[e_1,e_6]=e_7,[e_2,e_3]=e_6,[e_2,e_5]=e_7, [e_3,e_4]=e_7$

\

\noindent  $\{e_1...e_7\}$ is a nice basis.

$$U=\tiny{\left[ \begin {array}{cccccc} 3&1&1&1&-1&1\\ \noalign{\medskip}1&3&0&
1&0&1\\ \noalign{\medskip}1&0&3&-1&1&1\\ \noalign{\medskip}1&1&-1&3&1&
1\\ \noalign{\medskip}-1&0&1&1&3&1\\ \noalign{\medskip}1&1&1&1&1&3
\end {array} \right] }$$

\noindent  General solution to $Ux=[1]$: $x= (\frac{8}{19}-t,\frac{4}{19},\frac{1}{19}+t,t,\frac{9}{19}-t,-\frac{1}{19})^T$. It is not an Einstein nilradical.

\

\noindent  Maximal Torus: $D_1=\diag(1, 0, 2, 1, 3, 2, 3)$,  $D_2=\diag(0, 1, 0, 1, 0, 1, 1)$.
Pre-Einstein Derivation: $\frac{2}{19}\diag (3, 5, 6, 8, 9, 11, 14)$

\

\begin{center}
\begin{tabular}{cccccccccccc}
\hline
\multicolumn{3}{||m{1.5cm}|  }  {\centering  $2.1(v)$  } &
\multicolumn{3}{  m{1.5cm}|  }  {\centering   } &
\multicolumn{3}{  m{1.5cm}|  }  {\centering  $L_{74}^{34}(3\CC)$ } &
\multicolumn{3}{  m{4.25cm}||}  {\centering  $(2,4,7)Q$  } \tabularnewline
\hline
\hline
\multicolumn{3}{||m{1.5cm}|   }  {\centering $\dim$ $\Der$ \\               $14 $    } &
\multicolumn{6}{  m{3.75cm}  |   }  {\centering $\dim$ {Derived Series}  \\     $(7, 4, 0)$    } &
\multicolumn{3}{  m{4.25cm}|| }  {\centering $\dim$ {Desc. C. Series} \\     $(7, 4, 2, 0)$    } \tabularnewline
\hline
\end{tabular}
\end{center}

\

\noindent  $[e_1,e_2]=e_4,[e_1,e_3]=e_5,[e_1,e_4]=e_6,[e_2,e_3]=e_6,[e_2,e_5]=e_7, [e_3,e_4]=e_7$

\

\noindent  $\{e_1...e_7\}$ is a nice basis.

$$U=\tiny{\left[ \begin {array}{cccccc} 3&1&0&1&1&-1\\ \noalign{\medskip}1&3&1&
1&-1&1\\ \noalign{\medskip}0&1&3&1&0&1\\ \noalign{\medskip}1&1&1&3&1&1
\\ \noalign{\medskip}1&-1&0&1&3&1\\ \noalign{\medskip}-1&1&1&1&1&3
\end {array} \right]}$$

\noindent  General solution to $Ux=[1]$: $x= (\frac{1}{19}+t,\frac{8}{19}-t,\frac{4}{19},-\frac{1}{19},\frac{9}{19}-t,t)^T$. It is not an Einstein nilradical.

\

\noindent  Maximal Torus: $D_1=\diag(1, 0, 2, 1, 3, 2, 3)$,  $D_2=\diag(0, 1, 0, 1, 0, 1, 1)$.
Pre-Einstein Derivation: $\frac{2}{19}\diag (3, 5, 6, 8, 9, 11, 14)$

\

\begin{center}
\begin{tabular}{cccccccccccc}
\hline
\multicolumn{3}{||m{1.5cm}|  }  {\centering  $2.2$  } &
\multicolumn{3}{  m{1.5cm}|  }  {\centering  $\mathfrak{n}_{7,95}^{}$ } &
\multicolumn{3}{  m{1.5cm}|  }  {\centering  } &
\multicolumn{3}{  m{4.25cm}||}  {\centering  $(1,4,7)D$  } \tabularnewline
\hline
\hline
\multicolumn{3}{||m{1.5cm}|   }  {\centering $\dim$ $\Der$ \\               $15$    } &
\multicolumn{6}{  m{3.75cm}  |   }  {\centering $\dim$ {Derived Series}  \\     $(7, 4, 0)$    } &
\multicolumn{3}{  m{4.25cm}|| }  {\centering $\dim$ {Desc. C. Series} \\     $(7, 4, 1, 0)$    } \tabularnewline
\hline
\end{tabular}
\end{center}

\

\noindent  $\mu:=\left\{\begin{array}{l}
[e_1, e_2] = e_5,   [e_1, e_3] = e_6,    [e_1, e_4] = 2e_7,  [e_2, e_3] = e_4, {[e_2, e_6]} = e_7,  \\
{[e_3, e_5]} = -e_7,  [e_3, e_6] = e_7
\end{array}
\right.$

\

\noindent It is not an Einstein nilradical by \cite[Theorem 2.6]{FERNANDEZ-CULMA1}. In fact,

\

\noindent  Einstein Pre-Derivation $\phi=\frac{1}{2}\diag (1,1,1,2,2,2,3) $. \newline
Let $X=\diag (1,-1,0,-1,0,1,0) $. $X\in \ggo_{\phi}$ therefore $g_{t}=\exp(tX) \in \Grp_{\phi}$.

\

\noindent  $g_t\cdot\mu=\left\{\begin{array}{l}
[e_1, e_2] = e_5,    [e_1, e_3] = e_6,   [e_1, e_4] = 2e_7, [e_2, e_3] = e_4,{[e_2, e_6]} = e_7, \\
{[e_3, e_5]} = -e_7,  [e_3, e_6] = \euler^{-t}e_7
\end{array}\right.$

\

\noindent The $\Grp_{\phi}-$orbit of $\mu$ is not closed because  $g_t\cdot\mu \rightarrow \widetilde{\mu}$ as $t \rightarrow \infty$ and $(\RR^7,\widetilde{\mu})$ is non-isomorphic to $(\RR^7,\mu)$; $\dim \Der (\RR^7,\mu) = 15$ and $\dim \Der (\RR^7,\widetilde{\mu}) = 17$

\

\begin{center}
\begin{tabular}{cccccccccccc}
\hline
\multicolumn{3}{||m{1.5cm}|  }  {\centering  $2.3$  } &
\multicolumn{3}{  m{1.5cm}|  }  {\centering  $\mathfrak{n}_{7,10}^{}$ } &
\multicolumn{3}{  m{1.5cm}|  }  {\centering  $L_{76}^{1}$ } &
\multicolumn{3}{  m{4.25cm}||}  {\centering  $(1,2,3,4,5,7)A$  } \tabularnewline
\hline
\hline
\multicolumn{3}{||m{1.5cm}|   }  {\centering $\dim$ $\Der$ \\               $ 13$    } &
\multicolumn{6}{  m{3.75cm}  |   }  {\centering $\dim$ {Derived Series}  \\     $(7, 5, 0)$    } &
\multicolumn{3}{  m{4.25cm}|| }  {\centering $\dim$ {Desc. C. Series} \\     $(7, 5, 4, 3, 2, 1, 0)$    } \tabularnewline
\hline
\end{tabular}
\end{center}

\

\noindent  $[e_1,e_2]=e_3,[e_1,e_3]=e_4,[e_1,e_4]=e_5,[e_1,e_5]=e_6,[e_1,e_6]=e_7$

\

\noindent  $\{e_1...e_7\}$ is a nice basis.

$$U=\tiny{\left[ \begin {array}{ccccc} 3&0&1&1&1\\ \noalign{\medskip}0&3&0&1&1
\\ \noalign{\medskip}1&0&3&0&1\\ \noalign{\medskip}1&1&0&3&0
\\ \noalign{\medskip}1&1&1&0&3\end {array} \right]}$$

\noindent  General solution to $Ux=[1]$: $x= \frac{1}{37}(5,8,9,8,5)^T$. It is an Einstein nilradical.

\

\noindent  Maximal Torus: $D_1=\diag(1, 0, 1, 2, 3, 4, 5 )$,  $D_2=\diag(0, 1, 1, 1, 1, 1, 1 )$.
Pre-Einstein Derivation: $\frac{2}{37}\diag (1, 16, 17, 18, 19, 20, 21 )$ \newline
\noindent $||\mathcal{S}_{\beta}||^2=\frac{37}{35}\thickapprox 1.06$

\

\begin{center}
\begin{tabular}{cccccccccccc}
\hline
\multicolumn{3}{||m{1.5cm}|  }  {\centering  $2.4$  } &
\multicolumn{3}{  m{1.5cm}|  }  {\centering  $\mathfrak{n}_{7,13}^{}$ } &
\multicolumn{3}{  m{1.5cm}|  }  {\centering  $L_{74}^{37}(\mathfrak{n}_3)$ } &
\multicolumn{3}{  m{4.25cm}||}  {\centering  $(2,3,4,5,7)C$  } \tabularnewline
\hline
\hline
\multicolumn{3}{||m{1.5cm}|   }  {\centering $\dim$ $\Der$ \\               $12 $    } &
\multicolumn{6}{  m{3.75cm}  |   }  {\centering $\dim$ {Derived Series}  \\     $(7, 5, 1, 0)$    } &
\multicolumn{3}{  m{4.25cm}|| }  {\centering $\dim$ {Desc. C. Series} \\     $(7, 5, 4, 3, 2, 0)$    } \tabularnewline
\hline
\end{tabular}
\end{center}

\

\noindent  $[e_1,e_2]=e_3,[e_1,e_3]=e_4,[e_1,e_4]=e_5,[e_1,e_5]=e_6,[e_2,e_5]=-e_7, [e_3,e_4]=e_7$

\

\noindent  $\{e_1...e_7\}$ is a nice basis.

$$U=\tiny{\left[ \begin {array}{cccccc} 3&0&1&1&1&-1\\ \noalign{\medskip}0&3&0&
1&0&0\\ \noalign{\medskip}1&0&3&0&-1&1\\ \noalign{\medskip}1&1&0&3&1&0
\\ \noalign{\medskip}1&0&-1&1&3&1\\ \noalign{\medskip}-1&0&1&0&1&3
\end {array} \right]}$$

\noindent  General solution to $Ux=[1]$: $x= \frac{1}{52}(12,16,13,4,12,13)^T$. It is an Einstein nilradical.

\

\noindent  Maximal Torus: $D_1=\diag(1, 0, 1, 2, 3, 4, 3 )$,  $D_2=\diag(0, 1, 1, 1, 1, 1, 2 )$.
Pre-Einstein Derivation: $\frac{ 7}{ 52}\diag (1, 4, 5, 6, 7, 8, 11 )$ \newline
\noindent $||\mathcal{S}_{\beta}||^2=\frac{26}{35}\thickapprox 0.743$

\

\begin{center}
\begin{tabular}{cccccccccccc}
\hline
\multicolumn{3}{||m{1.5cm}|  }  {\centering  $2.5$  } &
\multicolumn{3}{  m{1.5cm}|  }  {\centering  $\mathfrak{n}_{7,24}^{}$ } &
\multicolumn{3}{  m{1.5cm}|  }  {\centering  $L_{74}^{19}(\mathfrak{n}_3)$ } &
\multicolumn{3}{  m{4.25cm}||}  {\centering  $(1,2,4,5,7)H$  } \tabularnewline
\hline
\hline
\multicolumn{3}{||m{1.5cm}|   }  {\centering $\dim$ $\Der$ \\               $12 $    } &
\multicolumn{6}{  m{3.75cm}  |   }  {\centering $\dim$ {Derived Series}  \\     $(7, 5, 1, 0)$    } &
\multicolumn{3}{  m{4.25cm}|| }  {\centering $\dim$ {Desc. C. Series} \\     $(7, 5, 4, 2, 1, 0)$    } \tabularnewline
\hline
\end{tabular}
\end{center}

\

\noindent  $[e_1,e_2]=e_3,[e_1,e_3]=e_4,[e_1,e_5]=e_6,[e_1,e_6]=e_7,[e_2,e_3]=e_5, [e_2,e_4]=e_6,[e_3,e_4]=e_7$

\

\noindent  $\{e_1...e_7\}$ is a nice basis.

$$U=\tiny{\left[ \begin {array}{ccccccc} 3&0&1&1&0&1&-1\\ \noalign{\medskip}0&3
&1&1&1&-1&0\\ \noalign{\medskip}1&1&3&0&-1&1&0\\ \noalign{\medskip}1&1
&0&3&0&-1&1\\ \noalign{\medskip}0&1&-1&0&3&1&1\\ \noalign{\medskip}1&-
1&1&-1&1&3&1\\ \noalign{\medskip}-1&0&0&1&1&1&3\end {array} \right]}$$

\noindent  General solution to $Ux=[1]$: $x=(t_1, -\frac{1}{5}+t_1+t_2, \frac{3}{5}-t_1-t_2, \frac{2}{5}-t_1, \frac{3}{5}-t_1-t_2, t_2, t_1)$ a particular solution is $x= \frac{1}{5}( 1,1,1,1,1,1,1)^T$. It is an Einstein nilradical.

\

\noindent  Maximal Torus: $D_1=\diag(1, 0, 1, 2, 1, 2, 3 )$,  $D_2=\diag(0, 1, 1, 1, 2, 2, 2 )$.
Pre-Einstein Derivation: $\frac{1 }{5 }\diag ( 1, 2, 3, 4, 5, 6, 7)$ \newline
\noindent $||\mathcal{S}_{\beta}||^2=\frac{5}{7}\thickapprox 0.714$

\

\begin{center}
\begin{tabular}{cccccccccccc}
\hline
\multicolumn{3}{||m{1.5cm}|  }  {\centering  $2.6 $  } &
\multicolumn{3}{  m{1.5cm}|  }  {\centering  $\mathfrak{n}_{7,28}^{}$ } &
\multicolumn{3}{  m{1.5cm}|  }  {\centering  $L_{74}^{31}(\mathfrak{n}_3)$ } &
\multicolumn{3}{  m{4.25cm}||}  {\centering  $(2,3,4,5,7)B$  } \tabularnewline
\hline
\hline
\multicolumn{3}{||m{1.5cm}|   }  {\centering $\dim$ $\Der$ \\               $12 $    } &
\multicolumn{6}{  m{3.75cm}  |   }  {\centering $\dim$ {Derived Series}  \\     $(7, 5, 1, 0)$    } &
\multicolumn{3}{  m{4.25cm}|| }  {\centering $\dim$ {Desc. C. Series} \\     $(7, 5, 4, 2, 1, 0)$    } \tabularnewline
\hline
\end{tabular}
\end{center}

\

\noindent  $[e_1,e_2]=e_3,[e_1,e_3]=e_4,[e_1,e_4]=e_6,[e_2,e_3]=e_5,[e_2,e_6]=e_7, [e_3,e_4]=-e_7$

\

\noindent  $\{e_1...e_7\}$ is a nice basis.

$$U=\tiny{ \left[ \begin {array}{cccccc} 3&0&1&0&1&-1\\ \noalign{\medskip}0&3&0&
1&0&0\\ \noalign{\medskip}1&0&3&0&-1&1\\ \noalign{\medskip}0&1&0&3&1&1
\\ \noalign{\medskip}1&0&-1&1&3&1\\ \noalign{\medskip}-1&0&1&1&1&3
\end {array} \right]}$$

\noindent  A solution to $Ux=[1]$: $x= \frac{1}{52}(13,16,13,4,12,12)^T$. It is an Einstein nilradical.

\

\noindent  Maximal Torus: $D_1=\diag( 1, 0, 1, 2, 1, 3, 3)$,  $D_2=\diag( 0, 1, 1, 1, 2, 1, 2)$.
Pre-Einstein Derivation: $\frac{ 1}{ 52}\diag (10, 23, 33, 43, 56, 53, 76 )$ \newline
\noindent $||\mathcal{S}_{\beta}||^2=\frac{26}{35}\thickapprox 0.743$

\

\begin{center}
\begin{tabular}{cccccccccccc}
\hline
\multicolumn{3}{||m{1.5cm}|  }  {\centering  $2.7 $  } &
\multicolumn{3}{  m{1.5cm}|  }  {\centering  $\mathfrak{n}_{7,31}^{}$ } &
\multicolumn{3}{  m{1.5cm}|  }  {\centering  $L_{75}^{35}$ } &
\multicolumn{3}{  m{4.25cm}||}  {\centering  $(2,3,4,5,7)A$  } \tabularnewline
\hline
\hline
\multicolumn{3}{||m{1.5cm}|   }  {\centering $\dim$ $\Der$ \\               $13 $    } &
\multicolumn{6}{  m{3.75cm}  |   }  {\centering $\dim$ {Derived Series}  \\     $(7, 5, 0)$    } &
\multicolumn{3}{  m{4.25cm}|| }  {\centering $\dim$ {Desc. C. Series} \\     $(7, 5, 4, 2, 1, 0)$    } \tabularnewline
\hline
\end{tabular}
\end{center}

\

\noindent  $[e_1,e_2]=e_3,[e_1,e_3]=e_4,[e_1,e_4]=e_6,[e_1,e_6]=e_7,[e_2,e_3]=e_5]$

\

\noindent  $\{e_1...e_7\}$ is a nice basis.

$$U=\tiny{ \left[ \begin {array}{ccccc} 3&0&1&1&0\\ \noalign{\medskip}0&3&0&1&1
\\ \noalign{\medskip}1&0&3&0&0\\ \noalign{\medskip}1&1&0&3&0
\\ \noalign{\medskip}0&1&0&0&3\end {array} \right] }$$

\noindent  General solution to $Ux=[1]$: $x= \frac{1}{18}( 3,3,5,4,5)^T$. It is an Einstein nilradical.

\

\noindent  Maximal Torus: $D_1=\diag( 1, 0, 1, 2, 1, 3, 4)$,  $D_2=\diag( 0, 1, 1, 1, 2, 1, 1)$.
Pre-Einstein Derivation: $\frac{1 }{18 }\diag ( 3, 10, 13, 16, 23, 19, 22)$ \newline
\noindent $||\mathcal{S}_{\beta}||^2=\frac{9}{10}\thickapprox 0.9$

\

\begin{center}
\begin{tabular}{cccccccccccc}
\hline
\multicolumn{3}{||m{1.5cm}|  }  {\centering  $2.8$  } &
\multicolumn{3}{  m{1.5cm}|  }  {\centering  $\mathfrak{n}_{7,32}^{}$ } &
\multicolumn{3}{  m{1.5cm}|  }  {\centering  $L_{75}^{38}$ } &
\multicolumn{3}{  m{4.25cm}||}  {\centering  $(2,4,5,7)M$  } \tabularnewline
\hline
\hline
\multicolumn{3}{||m{1.5cm}|   }  {\centering $\dim$ $\Der$ \\               $ 13$    } &
\multicolumn{6}{  m{3.75cm}  |   }  {\centering $\dim$ {Derived Series}  \\     $(7, 5, 0)$    } &
\multicolumn{3}{  m{4.25cm}|| }  {\centering $\dim$ {Desc. C. Series} \\     $(7, 5, 4, 2, 0)$    } \tabularnewline
\hline
\end{tabular}
\end{center}

\

\noindent  $[e_1,e_2]=e_3,[e_1,e_3]=e_4,[e_1,e_4]=e_6,[e_1,e_5]=e_7,[e_2,e_3]=e_5,[e_2,e_4]=e_7$

\

\noindent $\{e_1...e_7\}$ is a nice basis.

$$U=\tiny{  \left[ \begin {array}{cccccc} 3&0&1&1&0&1\\ \noalign{\medskip}0&3&0&1
&1&-1\\ \noalign{\medskip}1&0&3&1&0&1\\ \noalign{\medskip}1&1&1&3&-1&1
\\ \noalign{\medskip}0&1&0&-1&3&1\\ \noalign{\medskip}1&-1&1&1&1&3
\end {array} \right]}$$

\noindent General solution to $Ux=[1]$: $x=(\frac{1}{6}, \frac{1}{12}+t, \frac{1}{6}, \frac{1}{3}-t, \frac{5}{12}-t, t)$ and a particular solution is $x= \frac{1}{24}(4,6,4,4,6,4)^T$. It is an Einstein nilradical.

\

\noindent Maximal Torus: $D_1=\diag(1, 0, 1, 2, 1, 3, 2 )$,  $D_2=\diag( 0, 1, 1, 1, 2, 1, 2)$.
Pre-Einstein Derivation: $\frac{1 }{12 }\diag (3, 5, 8, 11, 13, 14, 16 )$ \newline
\noindent $||\mathcal{S}_{\beta}||^2=\frac{6}{7}\thickapprox 0.857$

\

\begin{center}
\begin{tabular}{cccccccccccc}
\hline
\multicolumn{3}{||m{1.5cm}|  }  {\centering  $2.9$  } &
\multicolumn{3}{  m{1.5cm}|  }  {\centering  $\mathfrak{n}_{7,33}^{}$ } &
\multicolumn{3}{  m{1.5cm}|  }  {\centering  $L_{75}^{39}$ } &
\multicolumn{3}{  m{4.25cm}||}  {\centering  $(2,4,5,7)L$  } \tabularnewline
\hline
\hline
\multicolumn{3}{||m{1.5cm}|   }  {\centering $\dim$ $\Der$ \\               $12 $    } &
\multicolumn{6}{  m{3.75cm}  |   }  {\centering $\dim$ {Derived Series}  \\     $(7, 5, 0)$    } &
\multicolumn{3}{  m{4.25cm}|| }  {\centering $\dim$ {Desc. C. Series} \\     $(7, 5, 4, 2, 0)$    } \tabularnewline
\hline
\end{tabular}
\end{center}

\

\noindent  $[e_1,e_2]=e_3,[e_1,e_3]=e_4,[e_1,e_4]=e_6,[e_2,e_3]=e_5,[e_2,e_5]=e_7$

\

\noindent  $\{e_1...e_7\}$ is a nice basis.

$$U=\tiny{\left[ \begin {array}{ccccc} 3&0&1&0&1\\ \noalign{\medskip}0&3&0&1&0
\\ \noalign{\medskip}1&0&3&0&0\\ \noalign{\medskip}0&1&0&3&0
\\ \noalign{\medskip}1&0&0&0&3\end {array} \right]}$$

\noindent  General solution to $Ux=[1]$: $x= \frac{1 }{28 }( 4,7,8,7,8)^T$. It is an Einstein nilradical.

\noindent  Maximal Torus: $D_1=\diag(1, 0, 1, 2, 1, 3, 1 )$,  $D_2=\diag(0, 1, 1, 1, 2, 1, 3 )$.
Pre-Einstein Derivation: $\frac{9 }{28 }\diag (1, 1, 2, 3, 3, 4, 4 )$ \newline
\noindent $||\mathcal{S}_{\beta}||^2=\frac{14}{17}\thickapprox 0.824$

\

\begin{center}
\begin{tabular}{cccccccccccc}
\hline
\multicolumn{3}{||m{1.5cm}|  }  {\centering  $2.10 $  } &
\multicolumn{3}{  m{1.5cm}|  }  {\centering  $\mathfrak{n}_{7,39}^{}$ } &
\multicolumn{3}{  m{1.5cm}|  }  {\centering  $L_{74}^{8}(\mathfrak{n}_3)$ } &
\multicolumn{3}{  m{4.25cm}||}  {\centering  $(1,3,4,5,7)C$  } \tabularnewline
\hline
\hline
\multicolumn{3}{||m{1.5cm}|   }  {\centering $\dim$ $\Der$ \\               $ 12$    } &
\multicolumn{6}{  m{3.75cm}  |   }  {\centering $\dim$ {Derived Series}  \\     $(7, 4, 1, 0)$    } &
\multicolumn{3}{  m{4.25cm}|| }  {\centering $\dim$ {Desc. C. Series} \\     $(7, 4, 3, 2, 1, 0)$    } \tabularnewline
\hline
\end{tabular}
\end{center}

\

\noindent  $[e_1,e_2]=e_4,[e_1,e_3]=e_7,[e_1,e_4]=e_5,[e_1,e_5]=e_6,[e_2,e_6]=e_7, [e_4,e_5]=-e_7$

\

\noindent  $\{e_1...e_7\}$ is a nice basis.

$$U=\tiny{\left[ \begin {array}{cccccc} 3&1&0&1&1&-1\\ \noalign{\medskip}1&3&1&
1&1&1\\ \noalign{\medskip}0&1&3&0&0&0\\ \noalign{\medskip}1&1&0&3&-1&1
\\ \noalign{\medskip}1&1&0&-1&3&1\\ \noalign{\medskip}-1&1&0&1&1&3
\end {array} \right] }$$

\noindent  General solution to $Ux=[1]$: $x=(t,-\frac{1}{5},\frac{2}{5},\frac{3}{5}-t,\frac{3}{5}-t,t)^T$. It is not an Einstein nilradical.

\

\noindent  Maximal Torus: $D_1=\diag(1, 0, 2, 1, 2, 3, 3 )$,  $D_2=\diag(0, 1, 2, 1, 1, 1, 2 )$.
Pre-Einstein Derivation: $\frac{ 1}{ 5}\diag (1, 2, 6, 3, 4, 5, 7 )$

\

\begin{center}
\begin{tabular}{cccccccccccc}
\hline
\multicolumn{3}{||m{1.5cm}|  }  {\centering  $2.11$  } &
\multicolumn{3}{  m{1.5cm}|  }  {\centering  $\mathfrak{n}_{7,58}^{}$ } &
\multicolumn{3}{  m{1.5cm}|  }  {\centering  } &
\multicolumn{3}{  m{4.25cm}||}  {\centering  $(2,4,5,7)E$  } \tabularnewline
\hline
\hline
\multicolumn{3}{||m{1.5cm}|   }  {\centering $\dim$ $\Der$ \\               $14 $    } &
\multicolumn{6}{  m{3.75cm}  |   }  {\centering $\dim$ {Derived Series}  \\     $(7, 4, 0)$    } &
\multicolumn{3}{  m{4.25cm}|| }  {\centering $\dim$ {Desc. C. Series} \\     $(7, 4, 3, 1, 0)$    } \tabularnewline
\hline
\end{tabular}
\end{center}

\

\noindent  $\mu:=\left\{\begin{array}{l}[e_1,e_2]=e_4,  [e_1,e_4]=e_5, [e_1,e_5]=e_7,  [e_2,e_3]=e_6, {[e_2,e_4]}=e_6 \end{array}\right.$

\

\noindent Let $g \in \mathrm{GL}_7(\RR)$ and $g^{-1}$ its inverse

\noindent	\begin{tabular}{cc}
  $g=$ & $g^{-1}=$ \\
\scalebox{0.59}{$\left( \begin {array}{ccccccc}
1&0&0&0&0&0&0\\
0&-\frac{4\sqrt {1254}}{9}&0&0&0&0&0\\
0&0&1&-\frac{8}{3}&0&0&0\\
0&0&-\frac{4\sqrt {7}}{3}&-\frac{4\sqrt {7}}{3}&0&0&0\\
0&0&0&0&-{\frac {2\sqrt {1463}}{57}}&0&0\\
0&0&0&0&0&{\frac {88\sqrt {42}}{27}}&0\\
0&0&0&0&0&0&-{\frac {\sqrt {770}}{57}}
\end {array} \right)$}
 &
\scalebox{0.6}{$\left( \begin {array}{ccccccc}
1&0&0&0&0&0&0\\
0&-{\frac {3\sqrt {1254}}{1672}}&0&0&0&0&0\\
0&0&\frac{3}{11}  & -{\frac {6\sqrt {7}}{77}}&0&0&0\\
0&0&-\frac{3}{11} & -{\frac {9\sqrt {7}}{308}}&0&0&0\\
0&0&0&0&-{\frac {3\sqrt {1463}}{154}}&0&0\\
0&0&0&0&0&{\frac {9\sqrt{42}}{1232}}&0\\
0&0&0&0&0&0&-{\frac {57\sqrt{770}}{770}}\end {array} \right)$}
		\end{tabular}

\

\noindent The action of $g$ over $\mu$ gives a isomorphic Lie algebra law, $g\cdot \mu=\widetilde{\mu}$

\

\noindent $\widetilde{\mu}:=\left\{\begin{array}{l}
[e_1,e_2]={\frac {\sqrt{1254}}{209}} e_3 + {\frac {\sqrt{8778}}{418}}e_4, [e_1,e_3]={\frac {2\sqrt{1463}}{209}}e_5, [e_1,e_4]={\frac {3\sqrt{209}}{418}}e_5,\\
 {[e_1,e_5]}=\frac{\sqrt{190}}{38}e_7, [e_2,e_4]=\frac{\sqrt{209}}{38}e_6
\end{array}\right.$

\

\noindent By straight calculation, it is easy to see that the moment map of $\widetilde{\mu}$ with respect to the ordered basis $e_1,...,e_7$ is:

\noindent \begin{eqnarray*}
  \mm(\widetilde{\mu}) &=& \diag(-{\frac {27}{38}},-{\frac {17}{38}},-{\frac {4}{19}},-{\frac {4}{19}}
,\frac{1}{38},{\frac {11}{38}},{\frac {5}{19}})\\
                       &=&-\frac{18}{19}\mathrm{Id} + \frac{1}{38} \diag ( 9,19,28,28,37,47,46)
\end{eqnarray*}

\

\noindent Since $\diag ( 9,19,28,28,37,47,46)$ is a derivation of the Lie algebra $(\RR^7, \widetilde{\mu})$ then by \cite[Theorem 2.2]{FERNANDEZ-CULMA1}, it is an Einstein Nilradical.\newline
\noindent $||\mathcal{S}_{\beta}||^2=\frac{18}{19}\thickapprox 0.947$

\

\begin{center}
\begin{tabular}{cccccccccccc}
\hline
\multicolumn{3}{||m{1.5cm}|  }  {\centering  $2.12 $  } &
\multicolumn{3}{  m{1.5cm}|  }  {\centering  $\mathfrak{n}_{7,91}^{}$ } &
\multicolumn{3}{  m{1.5cm}|  }  {\centering  $L_{74}^{37}(3\CC)$ } &
\multicolumn{3}{  m{4.25cm}||}  {\centering  $(2,4,7)E$  } \tabularnewline
\hline
\hline
\multicolumn{3}{||m{1.5cm}|   }  {\centering $\dim$ $\Der$ \\               $14 $    } &
\multicolumn{6}{  m{3.75cm}  |   }  {\centering $\dim$ {Derived Series}  \\     $(7, 4, 0)$    } &
\multicolumn{3}{  m{4.25cm}|| }  {\centering $\dim$ {Desc. C. Series} \\     $(7, 4, 2, 0)$    } \tabularnewline
\hline
\end{tabular}
\end{center}

\

\noindent  $[e_1,e_2]=e_4,[e_1,e_3]=e_5,[e_1,e_4]=e_6,[e_2,e_4]=e_7,[e_3,e_5]=e_7$

\

\noindent  $\{e_1...e_7\}$ is a nice basis.

$$U=\tiny{\left[ \begin {array}{ccccc} 3&1&0&0&0\\ \noalign{\medskip}1&3&1&0&0
\\ \noalign{\medskip}0&1&3&1&0\\ \noalign{\medskip}0&0&1&3&1
\\ \noalign{\medskip}0&0&0&1&3\end {array} \right] }$$

\noindent General solution to $Ux=[1]$: $x= \frac{1 }{18 }( 5,3,4,3,5)^T$. It is an Einstein nilradical.

\

\noindent Maximal Torus: $D_1=\diag(1, 0, 0, 1, 1, 2, 1 )$,  $D_2=\diag(0, 1, 1, 1, 1, 1, 2 )$.
Pre-Einstein Derivation: $\frac{1 }{9 }\diag (3, 5, 5, 8, 8, 11, 13 )$ \newline
\noindent $||\mathcal{S}_{\beta}||^2=\frac{9}{10}\thickapprox 0.9$

\

\begin{center}
\begin{tabular}{cccccccccccc}
\hline
\multicolumn{3}{||m{1.5cm}|  }  {\centering  $2.13$  } &
\multicolumn{3}{  m{1.5cm}|  }  {\centering  $\mathfrak{n}_{7,37}^{}$ } &
\multicolumn{3}{  m{1.5cm}|  }  {\centering  $L_{74}^{16}(\mathfrak{n}_3)$ } &
\multicolumn{3}{  m{4.25cm}||}  {\centering  $(1,2,4,5,7)C$  } \tabularnewline
\hline
\hline
\multicolumn{3}{||m{1.5cm}|   }  {\centering $\dim$ $\Der$ \\               $12 $    } &
\multicolumn{6}{  m{3.75cm}  |   }  {\centering $\dim$ {Derived Series}  \\     $(7, 4, 1, 0)$    } &
\multicolumn{3}{  m{4.25cm}|| }  {\centering $\dim$ {Desc. C. Series} \\     $(7, 4, 3, 2, 1, 0)$    } \tabularnewline
\hline
\end{tabular}
\end{center}

\

\noindent  $[e_1,e_2]=e_4,[e_1,e_4]=e_5,[e_1,e_5]=e_6,[e_2,e_3]=e_6,[e_2,e_6]=e_7, [e_4,e_5]=-e_7 $

\

\noindent  $\{e_1...e_7\}$ is a nice basis.

$$U=\tiny{\left[ \begin{array}{cccccc} 3&0&1&1&1&-1\\ \noalign{\medskip}0&3&0&0
&0&0\\ \noalign{\medskip}1&0&3&1&-1&1\\ \noalign{\medskip}1&0&1&3&0&0
\\ \noalign{\medskip}1&0&-1&0&3&1\\ \noalign{\medskip}-1&0&1&0&1&3
\end{array} \right] }$$

\noindent  A solution to $Ux=[1]$: $x= \frac{ 1}{60 }(12,20,12,12,15,15 )^T$. It is an Einstein nilradical.

\

\noindent  Maximal Torus: $D_1=\diag(1, 0, 3, 1, 2, 3, 3 )$,  $D_2=\diag( 0, 1, 0, 1, 1, 1, 2)$.
Pre-Einstein Derivation: $\frac{1 }{60 }\diag ( 16, 21, 48, 37, 53, 69, 90)$ \newline
\noindent $||\mathcal{S}_{\beta}||^2=\frac{30}{43}\thickapprox 0.698$

\

\begin{center}
\begin{tabular}{cccccccccccc}
\hline
\multicolumn{3}{||m{1.5cm}|  }  {\centering  $2.14$  } &
\multicolumn{3}{  m{1.5cm}|  }  {\centering  $\mathfrak{n}_{7,45}^{}$ } &
\multicolumn{3}{  m{1.5cm}|  }  {\centering  $L_{75}^{9}$ } &
\multicolumn{3}{  m{4.25cm}||}  {\centering  $(1,2,3,5,7)A$  } \tabularnewline
\hline
\hline
\multicolumn{3}{||m{1.5cm}|   }  {\centering $\dim$ $\Der$ \\               $12$    } &
\multicolumn{6}{  m{3.75cm}  |   }  {\centering $\dim$ {Derived Series}  \\     $(7, 4, 0)$    } &
\multicolumn{3}{  m{4.25cm}|| }  {\centering $\dim$ {Desc. C. Series} \\     $(7, 4, 3, 2, 1, 0)$    } \tabularnewline
\hline
\end{tabular}
\end{center}

\

\noindent  $[e_1,e_2]=e_4,[e_1,e_4]=e_5,[e_1,e_5]=e_6,[e_1,e_6]=e_7,[e_2,e_3]=e_5,[e_3,e_4]=-e_6,[e_3,e_5]=-e_7$

\

\noindent  $\{e_1...e_7\}$ is a nice basis.

$$U=\tiny{\left[ \begin {array}{ccccccc} 3&0&1&1&1&-1&0\\ \noalign{\medskip}0&3
&0&1&1&1&-1\\ \noalign{\medskip}1&0&3&0&-1&1&1\\ \noalign{\medskip}1&1
&0&3&0&-1&1\\ \noalign{\medskip}1&1&-1&0&3&1&0\\ \noalign{\medskip}-1&
1&1&-1&1&3&1\\ \noalign{\medskip}0&-1&1&1&0&1&3\end {array} \right] }$$

\noindent  A solution to $Ux=[1]$: $x= \frac{1 }{5 }(1,1,1,1,1,1,1 )^T$. It is an Einstein nilradical.

\

\noindent  Maximal Torus: $D_1=\diag( 1, 0, 2, 1, 2, 3, 4)$,  $D_2=\diag(0, 1, 0, 1, 1, 1, 1 )$.
Pre-Einstein Derivation: $\frac{1 }{5 }\diag (1, 3, 2, 4, 5, 6, 7 )$ \newline
\noindent $||\mathcal{S}_{\beta}||^2=\frac{5}{7}\thickapprox 0.714$

\

\begin{center}
\begin{tabular}{cccccccccccc}
\hline
\multicolumn{3}{||m{1.5cm}|  }  {\centering  $2.15$  } &
\multicolumn{3}{  m{1.5cm}|  }  {\centering  $\mathfrak{n}_{7,48}^{}$ } &
\multicolumn{3}{  m{1.5cm}|  }  {\centering  $L_{75}^{10}$ } &
\multicolumn{3}{  m{4.25cm}||}  {\centering  $(1,2,4,5,7)A$  } \tabularnewline
\hline
\hline
\multicolumn{3}{||m{1.5cm}|   }  {\centering $\dim$ $\Der$ \\               $ 13$    } &
\multicolumn{6}{  m{3.75cm}  |   }  {\centering $\dim$ {Derived Series}  \\     $(7, 4, 0)$    } &
\multicolumn{3}{  m{4.25cm}|| }  {\centering $\dim$ {Desc. C. Series} \\     $(7, 4, 3, 2, 1, 0)$    } \tabularnewline
\hline
\end{tabular}
\end{center}

\

\noindent  $[e_1,e_2]=e_4,[e_1,e_4]=e_5,[e_1,e_5]=e_6,[e_1,e_6]=e_7,[e_2,e_3]=e_6, [e_3,e_4]=-e_7$

\

\noindent  $\{e_1...e_7\}$ is a nice basis.

$$U=\tiny{\left[ \begin {array}{cccccc} 3&0&1&1&1&-1\\ \noalign{\medskip}0&3&0&
1&0&1\\ \noalign{\medskip}1&0&3&0&1&0\\ \noalign{\medskip}1&1&0&3&-1&1
\\ \noalign{\medskip}1&0&1&-1&3&1\\ \noalign{\medskip}-1&1&0&1&1&3
\end {array} \right] }$$

\noindent A solution to $Ux=[1]$: $x= \frac{1 }{5 }(1,1,1,1,1,1 )^T$. It is an Einstein nilradical.

\

\noindent Maximal Torus: $D_1=\diag( 1, 0, 3, 1, 2, 3, 4)$,  $D_2=\diag(0, 1, 0, 1, 1, 1, 1 )$.
Pre-Einstein Derivation: $\frac{ 1}{ 5}\diag (1, 3, 3, 4, 5, 6, 7 )$ \newline
\noindent $||\mathcal{S}_{\beta}||^2=\frac{5}{6}\thickapprox 0.833$

\

\begin{center}
\begin{tabular}{cccccccccccc}
\hline
\multicolumn{3}{||m{1.5cm}|  }  {\centering  $2.16$  } &
\multicolumn{3}{  m{1.5cm}|  }  {\centering  $\mathfrak{n}_{7,51}^{}$ } &
\multicolumn{3}{  m{1.5cm}|  }  {\centering  $L_{75}^{11}$ } &
\multicolumn{3}{  m{4.25cm}||}  {\centering  $(1,3,4,5,7)A$  } \tabularnewline
\hline
\hline
\multicolumn{3}{||m{1.5cm}|   }  {\centering $\dim$ $\Der$ \\               $ 14$    } &
\multicolumn{6}{  m{3.75cm}  |   }  {\centering $\dim$ {Derived Series}  \\     $(7, 4, 0)$    } &
\multicolumn{3}{  m{4.25cm}|| }  {\centering $\dim$ {Desc. C. Series} \\     $(7, 4, 3, 2, 1, 0)$    } \tabularnewline
\hline
\end{tabular}
\end{center}

\

\noindent  $[e_1,e_2]=e_4,[e_1,e_4]=e_5,[e_1,e_5]=e_6,[e_1,e_6]=e_7,[e_2,e_3]=e_7$

\

\noindent  $\{e_1...e_7\}$ is a nice basis.

$$U=\tiny{ \left[ \begin {array}{ccccc} 3&0&1&1&1\\ \noalign{\medskip}0&3&0&1&0
\\ \noalign{\medskip}1&0&3&0&0\\ \noalign{\medskip}1&1&0&3&1
\\ \noalign{\medskip}1&0&0&1&3\end {array} \right]}$$

\noindent  General solution to $Ux=[1]$: $x= \frac{1 }{27 }(3,8,8,3,7 )^T$. It is an Einstein nilradical.

\

\noindent  Maximal Torus: $D_1=\diag(1, 0, 4, 1, 2, 3, 4 )$,  $D_2=\diag(0, 1, 0, 1, 1, 1, 1 )$.
Pre-Einstein Derivation: $\frac{1 }{ 27}\diag (5, 17, 20, 22, 27, 32, 37 )$ \newline
\noindent $||\mathcal{S}_{\beta}||^2=\frac{27}{29}\thickapprox 0.931$

\

\begin{center}
\begin{tabular}{cccccccccccc}
\hline
\multicolumn{3}{||m{1.5cm}|  }  {\centering  $2.17$  } &
\multicolumn{3}{  m{1.5cm}|  }  {\centering  $\mathfrak{n}_{7,56}^{}$ } &
\multicolumn{3}{  m{1.5cm}|  }  {\centering  $L_{74}^{40}(3\CC)$ } &
\multicolumn{3}{  m{4.25cm}||}  {\centering  $(2,3,5,7)C$  } \tabularnewline
\hline
\hline
\multicolumn{3}{||m{1.5cm}|   }  {\centering $\dim$ $\Der$ \\               $13 $    } &
\multicolumn{6}{  m{3.75cm}  |   }  {\centering $\dim$ {Derived Series}  \\     $(7, 4, 0)$    } &
\multicolumn{3}{  m{4.25cm}|| }  {\centering $\dim$ {Desc. C. Series} \\     $(7, 4, 3, 1, 0)$    } \tabularnewline
\hline
\end{tabular}
\end{center}

\

\noindent  $[e_1,e_2]=e_4,[e_1,e_4]=e_5,[e_1,e_5]=e_7,[e_2,e_3]=e_5,[e_2,e_4]=e_6, [e_3,e_4]=-e_7$

\

\noindent  $\{e_1...e_7\}$ is a nice basis.

$$U=\tiny{\left[ \begin {array}{cccccc} 3&0&1&1&0&-1\\ \noalign{\medskip}0&3&0&
1&1&1\\ \noalign{\medskip}1&0&3&-1&0&1\\ \noalign{\medskip}1&1&-1&3&1&
1\\ \noalign{\medskip}0&1&0&1&3&1\\ \noalign{\medskip}-1&1&1&1&1&3
\end {array} \right] }$$

\noindent  A solution to $Ux=[1]$: $x= \frac{ 1}{24 }(6,4,6,4,4,4 )^T$. It is  an Einstein nilradical.

\

\noindent  Maximal Torus: $D_1=\diag(1, 0, 2, 1, 2, 1, 3 )$,  $D_2=\diag(0, 1, 0, 1, 1, 2, 1 )$.
Pre-Einstein Derivation: $\frac{1 }{12 }\diag (4, 5, 8, 9, 13, 14, 17 )$ \newline
\noindent $||\mathcal{S}_{\beta}||^2=\frac{6}{7}\thickapprox 0.857$

\

\begin{center}
\begin{tabular}{cccccccccccc}
\hline
\multicolumn{3}{||m{1.5cm}|  }  {\centering  $2.18$  } &
\multicolumn{3}{  m{1.5cm}|  }  {\centering  $\mathfrak{n}_{7,60}^{}$ } &
\multicolumn{3}{  m{1.5cm}|  }  {\centering  $L_{75}^{59}$ } &
\multicolumn{3}{  m{4.25cm}||}  {\centering  $(2,4,5,7)H$  } \tabularnewline
\hline
\hline
\multicolumn{3}{||m{1.5cm}|   }  {\centering $\dim$ $\Der$ \\               $15 $    } &
\multicolumn{6}{  m{3.75cm}  |   }  {\centering $\dim$ {Derived Series}  \\     $(7, 4, 0)$    } &
\multicolumn{3}{  m{4.25cm}|| }  {\centering $\dim$ {Desc. C. Series} \\     $(7, 4, 3, 1, 0)$    } \tabularnewline
\hline
\end{tabular}
\end{center}

\

\noindent  $[e_1,e_2]=e_4,[e_1,e_4]=e_5,[e_1,e_5]=e_7,[e_2,e_3]=e_7,[e_2,e_4]=e_6$

\

\noindent  $\{e_1...e_7\}$ is a nice basis.

$$U=\tiny{\left[ \begin {array}{ccccc} 3&0&1&1&0\\ \noalign{\medskip}0&3&0&0&1
\\ \noalign{\medskip}1&0&3&1&0\\ \noalign{\medskip}1&0&1&3&1
\\ \noalign{\medskip}0&1&0&1&3\end {array} \right] }$$

\noindent  General solution to $Ux=[1]$: $x= \frac{1 }{68 }(15,18,15,8,14 )^T$. It is an Einstein nilradical.

\

\noindent  Maximal Torus: $D_1=\diag( 1, 0, 3, 1, 2, 1, 3)$,  $D_2=\diag(0, 1, 0, 1, 1, 2, 1 )$.
Pre-Einstein Derivation: $\frac{1 }{68 }\diag ( 20, 31, 60, 51, 71, 82, 91)$ \newline
\noindent $||\mathcal{S}_{\beta}||^2=\frac{34}{35}\thickapprox 0.971$

\

\begin{center}
\begin{tabular}{cccccccccccc}
\hline
\multicolumn{3}{||m{1.5cm}|  }  {\centering  $2.19$  } &
\multicolumn{3}{  m{1.5cm}|  }  {\centering  $\mathfrak{n}_{7,61}^{}$ } &
\multicolumn{3}{  m{1.5cm}|  }  {\centering  $L_{75}^{48}$ } &
\multicolumn{3}{  m{4.25cm}||}  {\centering  $(2,4,5,7)G$  } \tabularnewline
\hline
\hline
\multicolumn{3}{||m{1.5cm}|   }  {\centering $\dim$ $\Der$ \\               $ 15$    } &
\multicolumn{6}{  m{3.75cm}  |   }  {\centering $\dim$ {Derived Series}  \\     $(7, 4, 0)$    } &
\multicolumn{3}{  m{4.25cm}|| }  {\centering $\dim$ {Desc. C. Series} \\     $(7, 4, 3, 1, 0)$    } \tabularnewline
\hline
\end{tabular}
\end{center}

\

\noindent  $[e_1,e_2]=e_4,[e_1,e_3]=e_6,[e_1,e_4]=e_5,[e_1,e_5]=e_7,[e_2,e_4]=e_6$

\

\noindent  $\{e_1...e_7\}$ is a nice basis.

$$U=\tiny{\left[ \begin {array}{ccccc} 3&1&0&1&0\\ \noalign{\medskip}1&3&1&1&1
\\ \noalign{\medskip}0&1&3&0&1\\ \noalign{\medskip}1&1&0&3&0
\\ \noalign{\medskip}0&1&1&0&3\end {array} \right]}$$

\noindent  General solution to $Ux=[1]$: $x= \frac{1 }{4 }(1,0,1,1,1)^T$. It is not an Einstein nilradical.

\noindent  Maximal Torus: $D_1=\diag(1, 0, 0, 1, 2, 1, 3 )$,  $D_2=\diag(0, 1, 2, 1, 1, 2, 1 )$.
Pre-Einstein Derivation: $\frac{1 }{4 }\diag (1, 2, 4, 3, 4, 5, 5 )$

\

\begin{center}
\begin{tabular}{cccccccccccc}
\hline
\multicolumn{3}{||m{1.5cm}|  }  {\centering  $2.20 $  } &
\multicolumn{3}{  m{1.5cm}|  }  {\centering  $\mathfrak{n}_{7,76}^{}$ } &
\multicolumn{3}{  m{1.5cm}|  }  {\centering  $L_{75}^{46}$ } &
\multicolumn{3}{  m{4.25cm}||}  {\centering  $(2,4,5,7)F$  } \tabularnewline
\hline
\hline
\multicolumn{3}{||m{1.5cm}|   }  {\centering $\dim$ $\Der$ \\               $16$    } &
\multicolumn{6}{  m{3.75cm}  |   }  {\centering $\dim$ {Derived Series}  \\     $(7, 4, 0)$    } &
\multicolumn{3}{  m{4.25cm}|| }  {\centering $\dim$ {Desc. C. Series} \\     $(7, 4, 2, 1, 0)$    } \tabularnewline
\hline
\end{tabular}
\end{center}

\

\noindent  $[e_1,e_2]=e_4,[e_1,e_3]=e_5,[e_1,e_5]=e_6,[e_1,e_6]=e_7,[e_3,e_5]=e_7$

\

\noindent  $\{e_1...e_7\}$ is a nice basis.

$$U=\tiny{ \left[ \begin {array}{ccccc} 3&1&1&1&0\\ \noalign{\medskip}1&3&0&1&0
\\ \noalign{\medskip}1&0&3&0&1\\ \noalign{\medskip}1&1&0&3&1
\\ \noalign{\medskip}0&0&1&1&3\end {array} \right]}$$

\noindent  General solution to $Ux=[1]$: $x= \frac{ 1}{37 }(5,9,8,5,8)^T$. It is an Einstein nilradical.

\

\noindent Maximal Torus: $D_1=\diag(1, 0, 2, 1, 3, 4, 5 )$,  $D_2=\diag(0, 1, 0, 1, 0, 0, 0 )$.
Pre-Einstein Derivation: $\frac{ 2}{37 }\diag (5, 16, 10, 21, 15, 20, 25)$
\noindent $||\mathcal{S}_{\beta}||^2=\frac{37}{35}\thickapprox 1.06$

\

\begin{center}
\begin{tabular}{cccccccccccc}
\hline
\multicolumn{3}{||m{1.5cm}|  }  {\centering  $2.21 $  } &
\multicolumn{3}{  m{1.5cm}|  }  {\centering  $\mathfrak{n}_{7,77}^{}$ } &
\multicolumn{3}{  m{1.5cm}|  }  {\centering  $L_{75}^{2}$ } &
\multicolumn{3}{  m{4.25cm}||}  {\centering  $(2,4,5,7)C$  } \tabularnewline
\hline
\hline
\multicolumn{3}{||m{1.5cm}|   }  {\centering $\dim$ $\Der$ \\               $ 16$    } &
\multicolumn{6}{  m{3.75cm}  |   }  {\centering $\dim$ {Derived Series}  \\     $(7, 4, 0)$    } &
\multicolumn{3}{  m{4.25cm}|| }  {\centering $\dim$ {Desc. C. Series} \\     $(7, 4, 2, 1, 0)$    } \tabularnewline
\hline
\end{tabular}
\end{center}

\

\noindent  $[e_1,e_2]=e_4,[e_1,e_3]=e_5,[e_1,e_4]=e_6,[e_1,e_6]=e_7,[e_2,e_3]=e_7$

\

\noindent  $\{e_1...e_7\}$ is a nice basis.

$$U=\tiny{ \left[ \begin {array}{ccccc} 3&1&0&1&1\\ \noalign{\medskip}1&3&1&1&1
\\ \noalign{\medskip}0&1&3&0&0\\ \noalign{\medskip}1&1&0&3&1
\\ \noalign{\medskip}1&1&0&1&3\end {array} \right] }$$

\noindent  General solution to $Ux=[1]$: $x= \frac{ 1}{31 }( 6,1,10,6,6)^T$.It is an Einstein nilradical.

\

\noindent  Maximal Torus: $D_1=\diag(1, 0, 3, 1, 4, 2, 3 )$,  $D_2=\diag(0, 1, 0, 1, 0, 1, 1 )$.
Pre-Einstein Derivation: $\frac{1 }{31 }\diag (8, 19, 24, 27, 32, 35, 43 )$ \newline
\noindent $||\mathcal{S}_{\beta}||^2=\frac{31}{29}\thickapprox 1.07$

\

\begin{center}
\begin{tabular}{cccccccccccc}
\hline
\multicolumn{3}{||m{1.5cm}|  }  {\centering  $2.22 $  } &
\multicolumn{3}{  m{1.5cm}|  }  {\centering  $\mathfrak{n}_{7,79}^{}$ } &
\multicolumn{3}{  m{1.5cm}|  }  {\centering  $L_{75}^{57}$ } &
\multicolumn{3}{  m{4.25cm}||}  {\centering  $(2,4,5,7)I$  } \tabularnewline
\hline
\hline
\multicolumn{3}{||m{1.5cm}|   }  {\centering $\dim$ $\Der$ \\               $14 $    } &
\multicolumn{6}{  m{3.75cm}  |   }  {\centering $\dim$ {Derived Series}  \\     $(7,4,0)$    } &
\multicolumn{3}{  m{4.25cm}|| }  {\centering $\dim$ {Desc. C. Series} \\     $(7, 4, 2, 1, 0)$    } \tabularnewline
\hline
\end{tabular}
\end{center}

\

\noindent  $[e_1,e_3]=e_4,[e_1,e_4]=e_6,[e_1,e_6]=e_7,[e_2,e_3]=e_5,[e_3,e_4]=e_7 $

\

\noindent  $\{e_1...e_7\}$ is a nice basis.

$$U=\tiny{\left[ \begin {array}{ccccc} 3&0&1&1&0\\ \noalign{\medskip}0&3&0&0&1
\\ \noalign{\medskip}1&0&3&0&1\\ \noalign{\medskip}1&0&0&3&1
\\ \noalign{\medskip}0&1&1&1&3\end {array} \right]}$$

\noindent  General solution to $Ux=[1]$: $x= \frac{ 1}{19 }(3,6,5,5,1)^T$. It is an Einstein nilradical.

\

\noindent  Maximal Torus: $D_1=\diag(1, 0, 2, 3, 2, 4, 5)$,  $D_2=\diag(0, 1, 0, 0, 1, 0, 0 )$.
Pre-Einstein Derivation: $\frac{1 }{19 }\diag ( 5, 14, 10, 15, 24, 20, 25)$ \newline
\noindent $||\mathcal{S}_{\beta}||^2=\frac{19}{20}\thickapprox 0.950$

\

\begin{center}
\begin{tabular}{cccccccccccc}
\hline
\multicolumn{3}{||m{1.5cm}|  }  {\centering  $2.23 $  } &
\multicolumn{3}{  m{1.5cm}|  }  {\centering  $\mathfrak{n}_{7,107}^{}$ } &
\multicolumn{3}{  m{1.5cm}|  }  {\centering  $L_{74}^{6}(3\CC)$ } &
\multicolumn{3}{  m{4.25cm}||}  {\centering  $(1, 3, 7)?$ \\ \cite{CARLES1}  } \tabularnewline
\hline
\hline
\multicolumn{3}{||m{1.5cm}|   }  {\centering $\dim$ $\Der$ \\               $ 13$    } &
\multicolumn{6}{  m{3.75cm}  |   }  {\centering $\dim$ {Derived Series}  \\     $(7, 3, 0)$    } &
\multicolumn{3}{  m{4.25cm}|| }  {\centering $\dim$ {Desc. C. Series} \\     $(7, 3, 1, 0)$    } \tabularnewline
\hline
\end{tabular}
\end{center}

\

\noindent  $[e_1,e_4]=e_6,[e_1,e_6]=e_7,[e_2,e_3]=e_5,[e_2,e_5]=e_7,[e_3,e_4]=e_7$

\

\noindent  $\{e_1...e_7\}$ is a nice basis.

$$U=\tiny{\left[ \begin {array}{ccccc} 3&0&0&0&1\\ \noalign{\medskip}0&3&0&1&1
\\ \noalign{\medskip}0&0&3&0&1\\ \noalign{\medskip}0&1&0&3&1
\\ \noalign{\medskip}1&1&1&1&3\end {array} \right]}$$

\noindent  General solution to $Ux=[1]$: $x= \frac{ 1}{11 }(4,3,4,3,-1 )^T$. It is not an Einstein nilradical.

\

\noindent  Maximal Torus: $D_1=\diag( 1, 0, 2, 0, 2, 1, 2)$,  $D_2=\diag(0, 1, 0, 2, 1, 2, 2 )$.
Pre-Einstein Derivation: $\frac{ 4}{11 }\diag (1, 1, 2, 2, 3, 3, 4 )$

\

\begin{center}
\begin{tabular}{cccccccccccc}
\hline
\multicolumn{3}{||m{1.5cm}|  }  {\centering  $2.24 $  } &
\multicolumn{3}{  m{1.5cm}|  }  {\centering  $\mathfrak{n}_{7,73}^{}$ } &
\multicolumn{3}{  m{1.5cm}|  }  {\centering  $L_{75}^{17}$ } &
\multicolumn{3}{  m{4.25cm}||}  {\centering  $(2,3,5,7)A$  } \tabularnewline
\hline
\hline
\multicolumn{3}{||m{1.5cm}|   }  {\centering $\dim$ $\Der$ \\               $13 $    } &
\multicolumn{6}{  m{3.75cm}  |   }  {\centering $\dim$ {Derived Series}  \\     $(7, 4, 0)$    } &
\multicolumn{3}{  m{4.25cm}|| }  {\centering $\dim$ {Desc. C. Series} \\     $(7, 4, 2, 1, 0)$    } \tabularnewline
\hline
\end{tabular}
\end{center}

\

\noindent  $\mu:=\left\{\begin{array}{l}
[e_1, e_2] = e_4,  [e_1, e_4] = e_6,  [e_1, e_5] = -e_7,  [e_1, e_6] = e_7, {[e_2, e_3]} = e_5,  \\
{[e_3, e_4]} = e_7
\end{array}\right.$

\

\noindent Let $g \in \mathrm{GL}_7(\RR)$ and $g^{-1}$ its inverse

\noindent	\begin{tabular}{cc}
  $g=$ & $g^{-1}=$ \\
\scalebox{0.8125}{$\left( \begin {array}{ccccccc}
1&0&0&0&0&0&0 \\0&1&0&0&0&0&0 \\ 0&0&\frac{\sqrt {38}}{19}&0&0&0&0\\ 0&0&0&-\frac{\sqrt {38}}{19}&0&0&0\\
0&0&0&0&\frac{\sqrt {3}}{19}&-\frac{\sqrt {3}}{19}&0\\0&0&0&0&\frac{1}{19} & \frac{1}{19} &0\\ 0&0&0&0&0&0&-{\frac {\sqrt {114}}{361}}
\end {array} \right) $}
 &
\scalebox{0.76}{$\left(\begin {array}{ccccccc}
1&0&0&0&0&0&0\\
0&1&0&0&0&0&0\\
0&0&\frac{\sqrt {38}}{2}&0&0&0&0\\
0&0&0&-\frac{\sqrt {38}}{2}&0&0&0\\
0&0&0&0&{\frac {19\sqrt {3}}{6}}& \frac{19}{2}&0\\
0&0&0&0&-{\frac {19\sqrt {3}}{6}}& \frac{19}{2}&0\\
0&0&0&0&0&0&-{\frac {19\sqrt {114}}{6}}
\end {array} \right)$}
		\end{tabular}

\

\noindent The action of $g$ over $\mu$ gives a isomorphic Lie algebra law, $g\cdot \mu=\widetilde{\mu}$

\

\noindent $\widetilde{\mu}:=\left\{\begin{array}{l}
[e_1,e_2]=-\frac{\sqrt {38}}{19}e_4,  [e_1,e_4]=\frac{\sqrt{114}}{38}e_5 - \frac{\sqrt {38}}{38}e_6,
{[e_1,e_5]}=\frac{\sqrt {38}}{19}e_7, \\
{[e_2,e_3]}=\frac{\sqrt{114}}{38}e_5 + \frac{\sqrt {38}}{38}e_6,
{[e_3,e_4]}=\frac{\sqrt {114}}{38}e_7
\end{array}\right.$

\

\noindent By straight calculation, it is easy to see that the moment map of $\widetilde{\mu}$ with respect to the ordered basis $e_1,...,e_7$ is:

\noindent \begin{eqnarray*}
  \mm(\widetilde{\mu}) &=& \diag(-{\frac {12}{19}},-{\frac {8}{19}},-{\frac {7}{19}},-{\frac {3}{19}},
\frac{2}{19},\frac{2}{19},{\frac {7}{19}})\\
                       &=& -\frac{17}{19}\mathrm{Id} + \frac{1}{19}\diag(5,9,10,14,19,19,24)
\end{eqnarray*}

\

\noindent Since $\diag(5,9,10,14,19,19,24)$ is a derivation of the Lie algebra $(\RR^7, \widetilde{\mu})$ then by \cite[Theorem 2.2]{FERNANDEZ-CULMA1}, it is an Einstein Nilradical.\newline
\noindent $||\mathcal{S}_{\beta}||^2=\frac{17}{19}\thickapprox 0.895$

\

\begin{center}
\begin{tabular}{cccccccccccc}
\hline
\multicolumn{3}{||m{1.5cm}|  }  {\centering  $ 2.25$  } &
\multicolumn{3}{  m{1.5cm}|  }  {\centering  $\mathfrak{n}_{7,100}^{}$ } &
\multicolumn{3}{  m{1.5cm}|  }  {\centering  $L_{75}^{18}$ } &
\multicolumn{3}{  m{4.25cm}||}  {\centering  $(1,3,5,7)B$  } \tabularnewline
\hline
\hline
\multicolumn{3}{||m{1.5cm}|   }  {\centering $\dim$ $\Der$ \\               $14 $    } &
\multicolumn{6}{  m{3.75cm}  |   }  {\centering $\dim$ {Derived Series}  \\     $(7, 3, 0)$    } &
\multicolumn{3}{  m{4.25cm}|| }  {\centering $\dim$ {Desc. C. Series} \\     $(7, 3, 2, 1, 0)$    } \tabularnewline
\hline
\end{tabular}
\end{center}

\

\noindent  $\mu:=\left\{\begin{array}{l}
[e_1, e_2] = e_5,  [e_1, e_5] = e_6,  [e_1, e_6] = e_7,  [e_2, e_3] = e_6, {[e_3, e_4]} = -e_7,  \\
{[e_3, e_5]} = -e_7
\end{array}\right.$

\

\noindent Let $g \in \mathrm{GL}_7(\RR)$ and $g^{-1}$ its inverse

\noindent	\begin{tabular}{cc}
  $g=$ & $g^{-1}=$ \\
\scalebox{0.8125}{$\left( \begin {array}{ccccccc}
1&0&0&0&0&0&0\\
0&0&\frac{\sqrt {38}}{19}&0&0&0&0\\
0&1&0&0&0&0&0\\
0&0&0&\frac{\sqrt {38}}{38}&\frac{\sqrt {38}}{19}&0&0\\
0&0&0&\frac{\sqrt{114}}{38}&0&0&0\\
0&0&0&0&0&-\frac{\sqrt {3}}{19}&0\\
0&0&0&0&0&0&-{\frac {\sqrt {114}}{361}}
\end {array} \right)$}
 &
\scalebox{0.72}{$\left( \begin{array}{ccccccc}
1&0&0&0&0&0&0\\
0&0&1&0&0&0&0\\
0&\frac{\sqrt {38}}{2}&0&0&0&0&0\\
0&0&0&0&\frac{\sqrt {114}}{3}&0&0\\
0&0&0&\frac{\sqrt {38}}{2}&-\frac{\sqrt {114}}{6}&0&0\\
0&0&0&0&0&-{\frac {19\sqrt {3}}{3}}&0\\
0&0&0&0&0&0&-{\frac {19\sqrt {114}}{6}}
\end{array} \right)$}
		\end{tabular}

\

\noindent The action of $g$ over $\mu$ gives a isomorphic Lie algebra law, $g\cdot \mu=\widetilde{\mu}$

\

\noindent $\widetilde{\mu}:=\left\{\begin{array}{l}
[e_1,e_3]= \frac{\sqrt {38}}{19}e_4,  [e_1,e_4]=-\frac{\sqrt{114}}{38}e_6,  [e_1,e_5]=\frac{\sqrt {38}}{38}e_6, [e_1,e_6]=\frac{\sqrt {38}}{19}e_7,\\
{[e_2,e_3]}=\frac{\sqrt{114}}{38}e_6, [e_2,e_4]=\frac{\sqrt{114}}{38}e_7,   [e_2,e_5]=\frac{\sqrt {38}}{38}e_7
\end{array}
\right.$

\

\noindent By straight calculation, it is easy to see that the moment map of $\widetilde{\mu}$ with respect to the ordered basis $e_1,...,e_7$ is:

\noindent \begin{eqnarray*}
  \mm(\widetilde{\mu}) &=& \diag (-{\frac {12}{19}},-{\frac {7}{19}},-{\frac {7}{19}},-\frac{2}{19},-\frac{2}{19},{
\frac {3}{19}},{\frac {8}{19}})\\
                       &=&-\frac{17}{19}\mathrm{Id} + \frac{5}{19}\diag(1,2,2,3,3,4,5)
\end{eqnarray*}

\

\noindent Since $\diag(1,2,2,3,3,4,5)$ is a derivation of the Lie algebra $(\RR^7, \widetilde{\mu})$ then by \cite[Theorem 2.2]{FERNANDEZ-CULMA1}, it is an Einstein Nilradical.\newline
\noindent $||\mathcal{S}_{\beta}||^2=\frac{17}{19}\thickapprox 0.895$

\

\begin{center}
\begin{tabular}{cccccccccccc}
\hline
\multicolumn{3}{||m{1.5cm}|  }  {\centering  $2.26$  } &
\multicolumn{3}{  m{1.5cm}|  }  {\centering  $\mathfrak{n}_{7,85}^{}$ } &
\multicolumn{3}{  m{1.5cm}|  }  {\centering  $L_{74}^{38}(3\CC)$ } &
\multicolumn{3}{  m{4.25cm}||}  {\centering  $(2,4,7)J$  } \tabularnewline
\hline
\hline
\multicolumn{3}{||m{1.5cm}|   }  {\centering $\dim$ $\Der$ \\               $13 $    } &
\multicolumn{6}{  m{3.75cm}  |   }  {\centering $\dim$ {Derived Series}  \\     $(7, 4, 0)$    } &
\multicolumn{3}{  m{4.25cm}|| }  {\centering $\dim$ {Desc. C. Series} \\     $(7, 4, 2, 0)$    } \tabularnewline
\hline
\end{tabular}
\end{center}

\

\noindent $\widetilde{\mu}:=\left\{\begin{array}{l}
[e_1,e_2]=e_4,[e_1,e_3]=e_5,[e_1,e_5]=e_6,[e_2,e_5]=e_7,[e_3,e_4]=e_7, [e_3,e_5]=e_6
\end{array}\right.$

\

\noindent Let $g \in \mathrm{GL}_7(\RR)$ and $g^{-1}$ its inverse

\noindent	\begin{tabular}{cc}
  $g=$ & $g^{-1}=$ \\
\scalebox{0.8125}{$\left( \begin {array}{ccccccc}
2&0&1&0&0&0&0\\ 0&1&0&0&0&0&0\\ 0&0&\sqrt {3}&0&0&0&0\\ 0&0&0&\frac{\sqrt{114}}{19}&0&0&0\\
0&0&0&0&\frac{2\sqrt {114}}{19}  &0&0\\ 0&0&0&0&0&\frac {12}{19}&0\\
0&0&0&0&0&0& \frac {6}{19}\end {array} \right)$}
&
\scalebox{0.765}{$\left( \begin{array}{ccccccc}
\frac{1}{2}&0&-\frac{\sqrt {3}}{6}&0&0&0&0\\
0&1&0&0&0&0&0\\ 0&0&\frac{\sqrt{3}}{3}&0&0&0&0\\
0&0&0&\frac{\sqrt {114}}{6}&0&0&0\\0&0&0&0&\frac{\sqrt {114}}{12}&0&0\\ 0&0&0&0&0&{\frac {19}{12}}&0\\
0&0&0&0&0&0&{\frac {19}{6}}\end{array} \right)$}
		\end{tabular}

\

\noindent The action of $g$ over $\mu$ gives a isomorphic Lie algebra law, $g\cdot \mu=\widetilde{\mu}$

\

\noindent $\widetilde{\mu}:=\left\{\begin{array}{l}
{[e_1,e_2]}=\frac{\sqrt {114}}{38}e_4, [e_1,e_3]=\frac{1}{19}\sqrt { 38}e_5, [e_1,e_5]=\frac{\sqrt {114}}{38}e_6,\\
{[e_2,e_3]}=\frac{\sqrt { 38}}{38}e_4, [e_2,e_5]=\frac{1}{38}\sqrt {114}e_7, [e_3,e_4]=\frac{\sqrt { 38}}{19}e_7,\\
{[e_3,e_5]}=\frac{\sqrt {38}}{38}e_6
\end{array}
\right.$

\

\noindent By straight calculation, it is easy to see that the moment map of $\widetilde{\mu}$ with respect to the ordered basis $e_1,...,e_7$ is:

\noindent \begin{eqnarray*}
  \mm(\widetilde{\mu}) &=& \diag (-{\frac {10}{19}},-{\frac {7}{19}},-{\frac {10}{19}},0,-{\frac {3}{19
}},{\frac {4}{19}},{\frac {7}{19}})\\
                       &=&-\frac{17}{19}\mathrm{Id} + \frac{1}{19}\diag(7, 10, 7, 17, 14, 21, 24)
\end{eqnarray*}

\

\noindent Since $\diag(  7, 10, 7, 17, 14, 21, 24 )$ is a derivation of the Lie algebra $(\RR^7, \widetilde{\mu})$ then by \cite[Theorem 2.2]{FERNANDEZ-CULMA1}, it is an Einstein Nilradical.\newline
\noindent $||\mathcal{S}_{\beta}||^2=\frac{17}{19}\thickapprox 0.895$

\

\begin{center}
\begin{tabular}{cccccccccccc}
\hline
\multicolumn{3}{||m{1.5cm}|  }  {\centering  $2.27 $  } &
\multicolumn{3}{  m{1.5cm}|  }  {\centering  $\mathfrak{n}_{7,106}^{}$ } &
\multicolumn{3}{  m{1.5cm}|  }  {\centering   } &
\multicolumn{3}{  m{4.25cm}||}  {\centering  $(2,5,7)I$  } \tabularnewline
\hline
\hline
\multicolumn{3}{||m{1.5cm}|   }  {\centering $\dim$ $\Der$ \\               $17$    } &
\multicolumn{6}{  m{3.75cm}  |   }  {\centering $\dim$ {Derived Series}  \\     $(7, 3, 0)$    } &
\multicolumn{3}{  m{4.25cm}|| }  {\centering $\dim$ {Desc. C. Series} \\     $(7, 3, 2, 0)$    } \tabularnewline
\hline
\end{tabular}
\end{center}

\

\noindent  $\mu:=\left\{\begin{array}{l} [e_1,e_2]=e_5,  [e_1,e_3]=e_7,  [e_1,e_5]=e_6,  [e_2,e_4]=e_7,  [e_2,e_5]=e_7 \end{array}\right.$

\

\noindent Let $g \in \mathrm{GL}_7(\RR)$ and $g^{-1}$ its inverse

\noindent	\begin{tabular}{cc}
  $g=$ & $g^{-1}=$ \\
\scalebox{0.92}{$\left( \begin{array}{ccccccc}
1&0&0&0&0&0&0\\
0&\frac{\sqrt {26}}{2}&0&0&0&0&0\\
0&0&\frac{\sqrt {78}}{2}&0&0&0&0\\
0&0&0&1&\frac{1}{2}&0&0\\
0&0&0&0&\frac{\sqrt {3}}{2}&0&0\\
0&0&0&0&0&\frac{\sqrt {78}}{26}&0\\
0&0&0&0&0&0&\frac{\sqrt {3}}{2}
\end {array} \right)$}
 &
\scalebox{0.86}{$\left( \begin {array}{ccccccc}
1&0&0&0&0&0&0\\
0&\frac{\sqrt {26}}{13}&0&0&0&0&0\\
0&0&\frac{\sqrt {78}}{39}&0&0&0&0\\
0&0&0&1&-\frac{\sqrt {3}}{3}&0&0\\
0&0&0&0&\frac{2\sqrt {3}}{3}&0&0\\
0&0&0&0&0&\frac{\sqrt {78}}{3}&0\\
0&0&0&0&0&0&\frac{2\sqrt {3}}{3}
\end {array} \right)$}
		\end{tabular}

\

\noindent The action of $g$ over $\mu$ gives a isomorphic Lie algebra law, $g\cdot \mu=\widetilde{\mu}$

\

\noindent $\widetilde{\mu}:=\left\{\begin{array}{l}
[e_1,e_2]=\frac{\sqrt {26}}{26}e_4 + \frac{\sqrt {78}}{26}e_5,  [e_1,e_3]=\frac{\sqrt {26}}{26}e_7,  [e_1,e_5]=\frac{\sqrt {26}}{13}e_6,\\
{[e_2,e_4]}=\frac{\sqrt {78}}{26}e_7,  {[e_2,e_5]}=\frac{\sqrt {26}}{26}e_7
\end{array}\right.$

\

\noindent By straight calculation, it is easy to see that the moment map of $\widetilde{\mu}$ with respect to the ordered basis $e_1,...,e_7$ is:

\noindent \begin{eqnarray*}
  \mm(\widetilde{\mu}) &=& \diag (-{\frac {9}{13}},-{\frac {8}{13}},-\frac{1}{13},-\frac{2}{13},-\frac{2}{13},{\frac {4}{13}},{
\frac {5}{13}})\\
                       &=& -\frac{15}{13}\mathrm{Id} + \frac{1}{13} \diag (6,7,14,13,13,19,20)
\end{eqnarray*}

\

\noindent Since $ \diag (6,7,14,13,13,19,20)  $ is a derivation of the Lie algebra $(\RR^7, \widetilde{\mu})$ then by \cite[Theorem 2.2]{FERNANDEZ-CULMA1}, it is an Einstein Nilradical.\newline
\noindent $||\mathcal{S}_{\beta}||^2=\frac{15}{13}\thickapprox 1.15$

\

\begin{center}
\begin{tabular}{cccccccccccc}
\hline
\multicolumn{3}{||m{1.5cm}|  }  {\centering  $2.28 $  } &
\multicolumn{3}{  m{1.5cm}|  }  {\centering  $\mathfrak{n}_{7,112}^{}$ } &
\multicolumn{3}{  m{1.5cm}|  }  {\centering  $L_{74}^{11}(3\CC)$ } &
\multicolumn{3}{  m{4.25cm}||}  {\centering  $(1,4,7)B$  } \tabularnewline
\hline
\hline
\multicolumn{3}{||m{1.5cm}|   }  {\centering $\dim$ $\Der$ \\               $16$    } &
\multicolumn{6}{  m{3.75cm}  |   }  {\centering $\dim$ {Derived Series}  \\     $(7, 3, 0)$    } &
\multicolumn{3}{  m{4.25cm}|| }  {\centering $\dim$ {Desc. C. Series} \\     $(7, 3, 1, 0)$    } \tabularnewline
\hline
\end{tabular}
\end{center}

\

\noindent  $[e_1,e_2]=e_5,[e_1,e_3]=e_6,[e_1,e_6]=e_7,[e_2,e_5]=e_7,[e_3,e_4]=e_7$

\

\noindent  $\{e_1...e_7\}$ is a nice basis.

$$U=\tiny{ \left[ \begin {array}{ccccc} 3&1&1&0&0\\ \noalign{\medskip}1&3&0&0&1
\\ \noalign{\medskip}1&0&3&1&1\\ \noalign{\medskip}0&0&1&3&1
\\ \noalign{\medskip}0&1&1&1&3\end {array} \right]}$$

\noindent  General solution to $Ux=[1]$: $x= \frac{ 1}{37 }(8,8,5,9,5)^T$. It is an Einstein nilradical.

\

\noindent  Maximal Torus: $D_1=\diag(1, 0, -1, 2, 1, 0, 1 )$,  $D_2=\diag(0, 1, 2, 0, 1, 2, 2 )$.
Pre-Einstein Derivation: $\frac{ 4}{37 }\diag (4, 5, 6, 8, 9, 10, 14 )$ \newline
\noindent $||\mathcal{S}_{\beta}||^2=\frac{37}{35}\thickapprox 1.06$

\

\begin{center}
%
 \right] }$$

\noindent  General solution to $Ux=[1]$: $x= \frac{ 1}{41 }(14,12,-1,6,12 )^T$. It is not an Einstein nilradical.

\

\noindent  Maximal Torus: $D_1=\diag(1, 0, 2, -1, 1, 2, 1 )$,  $D_2=\diag( 0, 1, 0, 2, 1, 1, 2)$.
Pre-Einstein Derivation: $\frac{1 }{41 }\diag (15,22,30,29,37,52,59)$

\

\begin{center}
%
 \right]}$$

\noindent  General solution to $Ux=[1]$: $x= \frac{ 1}{27 }(8,8,3,3,7 )^T$. It is an Einstein nilradical.

\

\noindent  Maximal Torus: $D_1=\diag(1, 2, 0, 5, 3, 4, 5 )$,  $D_2=\diag(0, 0, 1, -1, 0, 0, 0 )$.
Pre-Einstein Derivation: $\frac{4 }{27 }\diag (2, 4, 5, 5, 6, 8, 10)$ \newline
\noindent $||\mathcal{S}_{\beta}||^2=\frac{27}{29}\thickapprox 0.931$

\

\begin{center}
%
 \right] }$$

\noindent  General solution to $Ux=[1]$: $x= \frac{1 }{39 }(3,13,9,12,9 )^T$. It is an Einstein nilradical.

\noindent  Maximal Torus: $D_1=\diag( 1, 0, 3, 1, 3, 2, 3)$,  $D_2=\diag(0, 1, -1, 1, 0, 1, 1 )$.
Pre-Einstein Derivation: $\frac{1 }{39 }\diag (14, 15, 27, 29, 42, 43, 57 )$ \newline
\noindent $||\mathcal{S}_{\beta}||^2=\frac{39}{46}\thickapprox 0.848$

\

\begin{center}
%
 \right] }$$

\noindent  General solution to $Ux=[1]$: $x= \frac{ 1}{27 }(7,3,8,3,8 )^T$. It is an Einstein nilradical.

\

\noindent  Maximal Torus: $D_1=\diag(1, 0, 1, 1, 2, 2, 3 )$,  $D_2=\diag( 0, 2, 1, 2, 1, 2, 2)$.
Pre-Einstein Derivation: $\frac{2 }{27 }\diag (3, 10, 8, 13, 11, 16, 19 )$ \newline
\noindent $||\mathcal{S}_{\beta}||^2=\frac{27}{29}\thickapprox 0.931$

\

\begin{center}
%
 \right] }$$

\noindent  General solution to $Ux=[1]$: $x= \frac{ 1}{33 }(6,11,6,9,9 )^T$. It is an Einstein nilradical.

\

\noindent  Maximal Torus: $D_1=\diag(2, 0, 3, 2, 3, 4, 6 )$,  $D_2=\diag( 0, 1, 0, 1, 1, 1, 1)$.
Pre-Einstein Derivation: $\frac{1 }{33 }\diag (10, 18, 15, 28, 33, 38, 48 )$ \newline
\noindent $||\mathcal{S}_{\beta}||^2=\frac{33}{41}\thickapprox 0.805$

\

\begin{center}
%
 \right]}$$

\noindent  General solution to $Ux=[1]$: $x= \frac{ 1}{47 }(12,11,2,15,15 )^T$. It is an Einstein nilradical.

\

\noindent  Maximal Torus: $D_1=\diag(2, 0, 1, 2, 3, 2, 4 )$,  $D_2=\diag( 0, 2, 1, 2, 1, 4, 2)$.
Pre-Einstein Derivation: $\frac{1 }{47 }\diag ( 22, 20, 21, 42, 43, 62, 64)$ \newline
\noindent $||\mathcal{S}_{\beta}||^2=\frac{47}{55}\thickapprox 0.854$

\

\begin{center}
%
 \right] }$$

\noindent  A solution to $Ux=[1]$: $x= \frac{ 1}{13 }(2,4,2,2,3,2 )^T$. It is an Einstein nilradical.

\

\noindent  Maximal Torus: $D_1=\diag(1, 0, -1, 1, 0, 1, 0 )$,  $D_2=\diag(0, 1, 2, 1, 2, 2, 3 )$.
Pre-Einstein Derivation: $\frac{1 }{13 }\diag (5, 6, 7, 11, 12, 17, 18 )$ \newline
\noindent $||\mathcal{S}_{\beta}||^2=\frac{13}{15}\thickapprox 0.867$

\

\begin{center}
%
 \right] }$$

\noindent General solution to $Ux=[1]$: $x= \frac{ 1}{23 }(1,4,6,4,6 )^T$. It is an Einstein nilradical.

\

\noindent Maximal Torus: $D_1=\diag(1, 0, 2, 3, 3, 2, 4)$,  $D_2=\diag(0, 1, -2, -3, -2, -1, -3 )$.
Pre-Einstein Derivation: $\frac{ 1}{ 23}\diag (18, 13, 10, 15, 28, 23, 33 )$ \newline
\noindent $||\mathcal{S}_{\beta}||^2=\frac{23}{21}\thickapprox 1.10$

\

\begin{center}
%

\

\noindent The action of $g$ over $\mu$ gives a isomorphic Lie algebra law, $g\cdot \mu=\widetilde{\mu}$

\

\noindent $\widetilde{\mu}:=\left\{\begin{array}{l}[e_1, e_2] = e_3, [e_1, e_3] = e_5, [e_1, e_4] = e_6, [e_1, e_5] = e_7, [e_2, e_3] = -e_6, \\
{[e_2, e_4]} = e_5, [e_2, e_6] = -e_7, [e_3, e_4] = 2e_7
\end{array}\right.$

\

\noindent $\{e_1...e_7\}$ is a nice basis of $(\RR^n , \widetilde{\mu} )$

$$U=\tiny{ \left[ %
 \right]}$$

\noindent A solution to $U_{\widetilde{\mu}}x=[1]$: $x=\frac{ 1}{22 }( 4,4,3,3,4,1,5,2)^T$. It is an Einstein nilradical.

\

\noindent Maximal Torus: $D_1=\diag( 1, 1, 2, 2, 3, 3, 4 )$,  \newline $D_2=\diag \left( \left(%
 \right]}$$

General solution to $Ux=[1]$: $x= \frac{ 1}{16 }(4,3,2,3,2 )^T$. It is an Einstein nilradical.

\

\noindent  Maximal Torus: $D_1=\diag(1, 0, 1, 2, 0, 2, 1 )$,  $D_2=\diag(0, 1, 1, 0, 2, 1, 2 )$.
Pre-Einstein Derivation: $\frac{ 7}{16 }\diag (1, 1, 2, 2, 2, 3, 3 )$ \newline
\noindent $||\mathcal{S}_{\beta}||^2=\frac{8}{7}\thickapprox 1.14$

\

\begin{center}
%
 \right]}$$

\noindent  General solution to $Ux=[1]$: $x= \frac{1 }{16 }( 2,3,2,4,3)^T$. It is an Einstein nilradical.

\

\noindent  Maximal Torus: $D_1=\diag(1, 0, 2, 1, 3, 2, 4 )$,  $D_2=\diag( 0, 1, 0, 1, 0, 1, 0)$.
Pre-Einstein Derivation: $\frac{ 1}{16 }\diag (5, 11, 10, 16, 15, 21, 20 )$ \newline
\noindent $||\mathcal{S}_{\beta}||^2=\frac{8}{7}\thickapprox 1.14$

\

\begin{center}
%
 \right]}$$

\noindent  General solution to $Ux=[1]$: $x= \frac{1 }{23 }(6,4,4,6,1)^T$. It is an Einstein nilradical.

\

\noindent  Maximal Torus: $D_1=\diag(1, 0, 1, 2, 2, 1, 3 )$,  $D_2=\diag( 0, 1, 1, 0, 1, 2, 0)$.
Pre-Einstein Derivation: $\frac{1 }{23 }\diag (9, 10, 19, 18, 28, 29, 27)$ \newline
\noindent $||\mathcal{S}_{\beta}||^2=\frac{23}{21}\thickapprox 1.10$

\

\begin{center}
%
 \right] }$$

\noindent  A solution to $Ux=[1]$: $x= \frac{ 1}{7 }(1,2,1,1,1,2)^T$. It is an Einstein nilradical.

\

\noindent  Maximal Torus: $D_1=\diag( 1, 0, 1, 0, 2, 1, 2)$,  $D_2=\diag(0, 1, 1, 2, 1, 2, 2 )$.
Pre-Einstein Derivation: $\frac{ 1}{7 }\diag (2, 3, 5, 6, 7, 8, 10 )$ \newline
\noindent $||\mathcal{S}_{\beta}||^2=\frac{7}{8}\thickapprox 0.875$

\

\begin{center}
%
 \right]}$$

\noindent General solution to $Ux=[1]$: $x= \frac{ 1}{41 }(6,12,12,-1,14 )^T$. It is not an Einstein nilradical.

\noindent Maximal Torus: $D_1=\diag( 1, 2, 0, 3, 1, 2, 5)$,  $D_2=\diag(0, 0, 1, 0, 1, 1, 0 )$.
Pre-Einstein Derivation: $\frac{1 }{41 }\diag (11, 22, 30, 33, 41, 52, 55 )$

\

\begin{center}
%
 \right]}$$

\noindent  General solution to $Ux=[1]$: $x= \frac{1 }{37 }(8,8,5,5,9 )^T$. It is an Einstein nilradical.

\

\noindent   Maximal Torus: $D_1=\diag(2, 0, 1, 2, 3, 5, 4 )$,  $D_2=\diag(0, 2, 1, 2, 1, 1, 2 )$.
Pre-Einstein Derivation: $\frac{1 }{ 37}\diag (11, 29, 20, 40, 31, 42, 51 )$ \newline
\noindent $||\mathcal{S}_{\beta}||^2=\frac{37}{35}\thickapprox 1.06$

\

\begin{center}
%
 \right]}$$

\noindent  General solution to $Ux=[1]$: $x= \frac{1 }{37 }(5,9,8,5,8 )^T$. It is an Einstein nilradical.

\

\noindent  Maximal Torus: $D_1=\diag( 1, 0, -1, 1, 0, -1, 1)$,  $D_2=\diag( 0, 1, 2, 1, 2, 3, 2)$.
Pre-Einstein Derivation: $\frac{1 }{37 }\diag ( 15, 19, 23, 34, 38, 42, 53)$ \newline
\noindent $||\mathcal{S}_{\beta}||^2=\frac{37}{35}\thickapprox 1.06$

\

\begin{center}
%
 \right]}$$

\noindent General solution to $Ux=[1]$: $x= \frac{ 1}{14 }(3,1,4,3,1 )^T$. It is an Einstein nilradical.

\

\noindent Maximal Torus: $D_1=\diag(1, 0, 3, 2, 1, 2, 3 )$,  $D_2=\diag(0, 1, -1, 0, 1, 1, 0 )$.
Pre-Einstein Derivation: $\frac{1 }{ 14}\diag (6, 7, 11, 12, 13, 19, 18 )$ \newline
\noindent $||\mathcal{S}_{\beta}||^2=\frac{7}{6}\thickapprox 1.17$

\newpage

\section{Rank three}

\

\begin{center}
%
 \right]}}$$

\noindent  A solution to $Ux=[1]$: $x= \frac{ 1}{ 8}(1,2,1,1,2,1 )^T$. It is an Einstein nilradical.

\

\noindent  Maximal Torus: $D_1=\diag(1, 0, 0, 1, 1, 0, 1 )$,  $D_2=\diag(0, 1, 0, 1, 0, 1, 1 )$, $D_3=\diag(0, 0, 1, 0, 1, 1, 1 )$.
Pre-Einstein Derivation: $\frac{1 }{2 }\diag ( 1, 1, 1, 2, 2, 2, 3)$
\noindent $||\mathcal{S}_{\beta}||^2=1$

\

\begin{center}
%
 \right]}}$$

\noindent  General solution to $Ux=[1]$: $x= (t,0,\frac{1}{2}-t,\frac{1}{2}-t,t )^T$. It is not an Einstein nilradical.

\

\noindent  Maximal Torus: $D_1=\diag(1, 0, 0, 1, 1, 0, 1 )$,  $D_2=\diag(0, 1, 0, 1, 0, 1, 1 )$, $D_3=\diag( 0, 0, 1, 0, 1, 1, 1)$.
Pre-Einstein Derivation: $\frac{ 1}{2 }\diag ( 1, 1, 1, 2, 2, 2, 3)$

\

\begin{center}
%
 \right]}}$$

\noindent General solution to $Ux=[1]$: $x=(0,t,\frac{1}{2}-t,\frac{1}{2}-t,t )^T$. It is not an Einstein nilradical.

\

\noindent Maximal Torus: $D_1=\diag(1, 0, 0, 1, 1, 0, 1 )$,  $D_2=\diag( 0, 1, 0, 1, 0, 1, 1)$, $D_3=\diag(0, 0, 1, 0, 1, 1, 1 )$.
Pre-Einstein Derivation: $\frac{1 }{2 }\diag (1, 1, 1, 2, 2, 2, 3 )$

\

\begin{center}
%
 \right]}}$$

\noindent  General solution to $Ux=[1]$: $x= (t,\frac{1}{2}-t,0,\frac{1}{2}-t,t)^T$. It is not an Einstein nilradical.

\

\noindent  Maximal Torus: $D_1=\diag( 1, 0, 0, 1, 1, 0, 1)$,  $D_2=\diag(0, 1, 0, 1, 0, 1, 1 )$, $D_3=\diag(0, 0, 1, 0, 1, 1, 1 )$.
Pre-Einstein Derivation: $\frac{ 1}{2 }\diag ( 1, 1, 1, 2, 2, 2, 3)$

\

\begin{center}
%
 \right]}}$$

\noindent  General solution to $Ux=[1]$: $x= \frac{1 }{13 }(3,1,4,3 )^T$. It is an Einstein nilradical.

\

\noindent  Maximal Torus: $D_1=\diag( 1, 0, 0, 1, 1, 2, 3)$,  $D_2=\diag( 0, 1, 0, 1, 0, 1, 1)$, $D_3=\diag( 0, 0, 1, 0, 1, 0, 0)$.
Pre-Einstein Derivation: $\frac{2 }{13 }\diag (1, 5, 6, 6, 7, 7, 8 )$ \newline
\noindent $||\mathcal{S}_{\beta}||^2=\frac{13}{11}\thickapprox 1.18$

\

\begin{center}
%
 \right]}}$$

\noindent  General solution to $Ux=[1]$: $x= \frac{1 }{21 }( 3,7,6,6)^T$. It is an Einstein nilradical.

\

\noindent  Maximal Torus: $D_1=\diag( 1, 0, 0, 1, 0, 2, 3)$,  $D_2=\diag( 0, 1, 0, 1, 1, 1, 1)$, $D_3=\diag( 0, 0, 1, 0, 1, 0, 0)$.
Pre-Einstein Derivation: $\frac{ 1}{21 }\diag ( 5, 12, 15, 17, 27, 22, 27)$ \newline
\noindent $||\mathcal{S}_{\beta}||^2=\frac{21}{22}\thickapprox 0.954$

\

\begin{center}
%
 \right]}}$$

\noindent  General solution to $Ux=[1]$: $x= \frac{1 }{12 }(3,3,4,4)^T$. It is an Einstein nilradical.

\

\noindent  Maximal Torus: $D_1=\diag( 1, 0, 0, 1, 1, 1, 1)$,  $D_2=\diag( 0, 1, 0, 1, 0, 2, 0)$, $D_3=\diag(0, 0, 1, 0, 1, 0, 2 )$.
Pre-Einstein Derivation: $\frac{ 1}{12 }\diag (6, 5, 5, 11, 11, 16, 16 )$ \newline
\noindent $||\mathcal{S}_{\beta}||^2=\frac{6}{7}\thickapprox 0.857$

\

\begin{center}
%
 \right]}}$$

\noindent  A solution to $Ux=[1]$: $x= \frac{1 }{20 }(5,5,4,4,4)^T$. It is an Einstein nilradical.

\

\noindent  Maximal Torus: $D_1=\diag(1, 0, 0, 1, 1, 1, 1 )$,  $D_2=\diag(0, 1, 0, 1, 0, 2, 1)$, $D_3=\diag( 0, 0, 1, 0, 1, 0, 1)$.
Pre-Einstein Derivation: $\frac{1 }{20 }\diag ( 10, 7, 11, 17, 21, 24, 28)$ \newline
\noindent $||\mathcal{S}_{\beta}||^2=\frac{10}{11}\thickapprox 0.909$

\

\begin{center}
%
 \right]}}$$

\noindent  General solution to $Ux=[1]$: $x= \frac{ 1}{13 }(1,3,4,3 )^T$. It is an Einstein nilradical.

\

\noindent  Maximal Torus: $D_1=\diag( 1, 0, 0, 1, 1, 0, 2)$,  $D_2=\diag(0, 1, 0, 1, 0, 1, 0 )$, $D_3=\diag( 0, 0, 1, 0, 1, 1, 1)$.
Pre-Einstein Derivation: $\frac{1 }{13 }\diag ( 5, 9, 7, 14, 12, 16, 17)$ \newline
\noindent $||\mathcal{S}_{\beta}||^2=\frac{13}{11}\thickapprox 1.18$

\

\begin{center}
%
 \right]}}$$

\noindent  General solution to $Ux=[1]$: $x= \frac{ 1}{3 }(1,1,0,1 )^T$. It is not an Einstein nilradical.

\

\noindent  Maximal Torus: $D_1=\diag(1, 0, 0, 2, 1, 2, 2 )$,  $D_2=\diag( 0, 1, 0, 0, 1, 1, 0)$, $D_3=\diag( 0, 0, 1, 0, 0, 0, 1)$.
Pre-Einstein Derivation: $\frac{1 }{ 3}\diag ( 1, 2, 2, 2, 3, 4, 4)$

\

\begin{center}
%
 \right]}}$$

\noindent  General solution to $Ux=[1]$: $x= \frac{ 1}{5 }(1,1,1,1)^T$. It is an Einstein nilradical.

\

\noindent  Maximal Torus: $D_1=\diag(1, 0, 0, 2, 1, 1, 2 )$,  $D_2=\diag( 0, 1, 0, 0, 1, 0, 1)$, $D_3=\diag( 0, 0, 1, 0, 0, 1, 0)$.
Pre-Einstein Derivation: $\frac{ 1}{5 }\diag ( 2, 3, 4, 4, 5, 6, 7)$ \newline
\noindent $||\mathcal{S}_{\beta}||^2=\frac{5}{4}=1.25 $

\

\begin{center}
%
 \right]}}$$

\noindent  General solution to $Ux=[1]$: $x= \frac{1 }{13 }( 3,4,3,1)^T$. It is an Einstein nilradical.

\noindent  Maximal Torus: $D_1=\diag( 1, 0, 0, 2, 1, 0, 2)$,  $D_2=\diag(0, 1, 0, 0, 1, 1, 1 )$, $D_3=\diag(0, 0, 1, 0, 0, 1, 0 )$.
Pre-Einstein Derivation: $\frac{ 2}{13 }\diag ( 3, 3, 5, 6, 6, 8, 9)$ \newline
\noindent $||\mathcal{S}_{\beta}||^2=\frac{13}{11}\thickapprox 1.18$

\

\begin{center}
%
 \right]}}$$

\noindent  A solution to $Ux=[1]$: $x= \frac{ 1}{20 }(4,4,4,5,5)^T$. It is an Einstein nilradical.

\

\noindent  Maximal Torus: $D_1=\diag(1, 0, 0, 1, 1, 1, 1 )$,  $D_2=\diag( 0, 1, 0, -1, 1, 0, 1)$, $D_3=\diag(0, 0, 1, 1, 0, 1, 1 )$.
Pre-Einstein Derivation: $\frac{1 }{20 }\diag ( 12, 7, 11, 16, 19, 23, 30)$ \newline
\noindent $||\mathcal{S}_{\beta}||^2=\frac{10}{11}\thickapprox $

\

\begin{center}
%
 \right]}}$$

\noindent  General solution to $Ux=[1]$: $x= \frac{ 1}{13 }(3,1,4,3 )^T$. It is an Einstein nilradical.

\

\noindent  Maximal Torus: $D_1=\diag(1, 0, 0, 1, 1, 1, 2 )$,  $D_2=\diag( 0, 1, 0, -1, 1, 0, 1)$, $D_3=\diag(0, 0, 1, 1, 0, 1, 0 )$.
Pre-Einstein Derivation: $\frac{ 1}{13 }\diag (5, 7, 12, 10, 12, 17, 17 )$ \newline
\noindent $||\mathcal{S}_{\beta}||^2=\frac{13}{11}\thickapprox 1.18 $

\

\begin{center}
%
 \right]}}$$

\noindent  General solution to $Ux=[1]$: $x= \frac{ 1}{ 8}(2,1,1,2)^T$. It is an Einstein nilradical.

\

\noindent  Maximal Torus: $D_1=\diag(1, 0, 0, 1, 1, 1, 1 )$,  $D_2=\diag(0, 1, 0, -1, 1, 0, -1 )$, $D_3=\diag(0, 0, 1, 1, 0, 1, 2 )$.
Pre-Einstein Derivation: $\frac{ 5}{ 8}\diag ( 1, 1, 1, 1, 2, 2, 2)$ \newline
\noindent $||\mathcal{S}_{\beta}||^2=\frac{4}{3}\thickapprox 1.33$

\

\begin{center}
%
 \right]}}$$

\noindent  General solution to $Ux=[1]$: $x= \frac{ 1}{ 21}(7,6,3,6)^T$. It is an Einstein nilradical.

\

\noindent  Maximal Torus: $D_1=\diag( 1, 0, 0, 1, 1, 2, 1)$,  $D_2=\diag( 0, 1, 0, 2, 1, 1, 2)$, $D_3=\diag( 0, 0, 1, -1, 0, 0, 0 )$.
Pre-Einstein Derivation: $\frac{1 }{21 }\diag ( 8, 11, 15, 15, 19, 27, 30)$\newline
\noindent $||\mathcal{S}_{\beta}||^2=\frac{21}{22}\thickapprox 0.954$

\

\begin{center}
%
 \right]}}$$

\noindent General solution to $Ux=[1]$: $x= \frac{ 1}{13 }(1,4,3,3 )^T$. It is an Einstein nilradical.

\

\noindent Maximal Torus: $D_1=\diag( 1, 0, 0, 2, 1, 1, 2)$,  $D_2=\diag( 0, 1, 0, -1, 1, 0, 0)$, $D_3=\diag(0, 0, 1, 1, 0, 1, 1 )$.
Pre-Einstein Derivation: $\frac{1 }{13 }\diag (5, 9, 9, 10, 14, 14, 19 )$ \newline
\noindent $||\mathcal{S}_{\beta}||^2=\frac{13}{11}\thickapprox 1.18 $

\

\begin{center}
%
 \right]}}$$

\noindent  General solution to $Ux=[1]$: $x= \frac{1 }{11 }(3,2,3,2)^T$. It is an Einstein nilradical.

\

\noindent  Maximal Torus: $D_1=\diag(1, 0, 0, 1, 1, 1, 1 )$,  $D_2=\diag( 0, 1, 0, 2, 1, 0, 2)$, $D_3=\diag( 0, 0, 1, -1, 0, 1, 0)$.
Pre-Einstein Derivation: $\frac{ 1}{11 }\diag (6, 5, 7, 9, 11, 13, 16 )$\newline
\noindent $||\mathcal{S}_{\beta}||^2=\frac{11}{10}=1.10$

\

\begin{center}
%
 \right]}}$$

\noindent  General solution to $Ux=[1]$: $x= \frac{ 1}{12 }(4,3,4,3)^T$. It is an Einstein nilradical.

\

\noindent  Maximal Torus: $D_1=\diag(1, 0, 0, 2, 1, 2, 2 )$,  $D_2=\diag( 0, 1, 0, 1, 1, 1, 1)$, $D_3=\diag(0, 0, 1, -2, 0, -1, 0 )$.
Pre-Einstein Derivation: $\frac{1 }{12 }\diag ( 5, 8, 5, 8, 13, 13, 18)$ \newline
\noindent $||\mathcal{S}_{\beta}||^2=\frac{6}{7}\thickapprox 0.857$

\

\begin{center}
%
 \right]}}$$

\noindent  General solution to $Ux=[1]$: $x= \frac{1 }{21 }(6,7,3,6)^T$. It is an Einstein nilradical.

\

\noindent  Maximal Torus: $D_1=\diag(1, 0, 0, 3, 1, 2, 3 )$,  $D_2=\diag(0, 1, 0, 1, 1, 1, 1 )$, $D_3=\diag( 0, 0, 1, -1, 0, 0, 0)$.
Pre-Einstein Derivation: $\frac{5 }{21 }\diag (1, 3, 3, 3, 4, 5, 6 )$\newline
\noindent $||\mathcal{S}_{\beta}||^2=\frac{21}{22}\thickapprox 0.954$

\

\begin{center}
%
 \right]}}$$

\noindent  General solution to $Ux=[1]$: $x= \frac{ 1}{13 }(4,3,1,3 )^T$. It is an Einstein nilradical.

\

\noindent  Maximal Torus: $D_1=\diag(1, 0, 0, 2, 2, 1, 2 )$,  $D_2=\diag( 0, 1, 0, 1, 0, 1, 1)$, $D_3=\diag( 0, 0, 1, -1, 0, 0, 0)$.
Pre-Einstein Derivation: $\frac{2 }{13 }\diag (3, 4, 5, 5, 6, 7, 10 )$ \newline
\noindent $||\mathcal{S}_{\beta}||^2=\frac{13}{11}\thickapprox 1.18$

\

\begin{center}
%
 \right]}}$$

\noindent  General solution to $Ux=[1]$: $x= \frac{ 1}{8 }(2,1,1,2)^T$. It is an Einstein nilradical.

\

\noindent  Maximal Torus: $D_1=\diag( 1, 0, 0, 1, 0, 1, 1)$,  $D_2=\diag(0, 1, 0, 1, -1, 1, 0 )$, $D_3=\diag(0, 0, 1, -1, 2, 0, 1 )$.
Pre-Einstein Derivation: $\frac{ 1}{8 }\diag ( 5, 6, 6, 5, 6, 11, 11)$ \newline
\noindent $||\mathcal{S}_{\beta}||^2=\frac{4}{3}\thickapprox 1.33$

\

\begin{center}
%
 \right]}}}$$

\noindent  General solution to $Ux=[1]$: $x= \frac{1 }{5 }( 1,1,1,1)^T$. It is an Einstein nilradical.

\

\noindent  Maximal Torus: $D_1=\diag( 1, 0, 0, 1, 1, 2, 2)$,  $D_2=\diag(0, 1, 0, 1, 0, 1, 0 )$, $D_3=\diag(0, 0, 1, 0, 1, 0, 1 )$.
Pre-Einstein Derivation: $\frac{1 }{ 5}\diag (1, 4, 4, 5, 5, 6, 6 )$ \newline
\noindent $||\mathcal{S}_{\beta}||^2=\frac{5}{4}=1.25 $

\

\begin{center}
%
 \right]}}$$

\noindent  General solution to $Ux=[1]$: $x= \frac{1 }{21 }(6,3,6,7 )^T$. It is an Einstein nilradical.

\

\noindent  Maximal Torus: $D_1=\diag( 1, 0, 0, 1, 1, 2, 1)$,  $D_2=\diag(0, 1, 0, 1, 0, 1, 0 )$, $D_3=\diag( 0, 0, 1, 0, 1, 0, 2)$.
Pre-Einstein Derivation: $\frac{1 }{21 }\diag (6, 15, 11, 21, 17, 27, 28 )$ \newline
\noindent $||\mathcal{S}_{\beta}||^2=\frac{21}{22}\thickapprox 0.954$

\

\begin{center}
%
 \right]}}$$

\noindent  A solution to $Ux=[1]$: $x= \frac{1 }{20 }(4,5,4,4,5 )^T$. It is an Einstein nilradical.

\

\noindent  Maximal Torus: $D_1=\diag( 1, 0, 0, 1, 1, 1, 2)$,  $D_2=\diag(0, 1, 0, 1, 0, 1, 0 )$, $D_3=\diag( 0, 0, 1, 0, 1, 1, 1)$.
Pre-Einstein Derivation: $\frac{1 }{20 }\diag (7, 12, 10, 19, 17, 29, 24 )$ \newline
\noindent $||\mathcal{S}_{\beta}||^2=\frac{10}{11}\thickapprox 0.909$

\

\begin{center}
%
 \right]}}$$

\noindent  General solution to $Ux=[1]$: $x= \frac{1 }{11 }( 3,2,2,3)^T$. It is an Einstein nilradical.

\

\noindent  Maximal Torus: $D_1=\diag( 1, 0, 1, 0, 2, 1, 1)$,  $D_2=\diag( 0, 1, 1, 0, 1, 2, 0)$, $D_3=\diag( 0, 0, 0, 1, 0, 0, 1)$.
Pre-Einstein Derivation: $\frac{1 }{11 }\diag ( 4, 5, 9, 9, 13, 14, 13)$ \newline
\noindent $||\mathcal{S}_{\beta}||^2=\frac{11}{10}=1.10$

\

\begin{center}
%
 \right]}}$$

\noindent  General solution to $Ux=[1]$: $x= \frac{1 }{6 }(1,1,1,1 )^T$. It is an Einstein nilradical.

\

\noindent  Maximal Torus: $D_1=\diag(1, 0, 0, 1, 1, 0, 1 )$,  $D_2=\diag( 0, 1, 0, 1, 1, 1, 2)$, $D_3=\diag( 0, 0, 1, -1, 0, 1, -1)$.
Pre-Einstein Derivation: $\frac{ 1}{ 6}\diag ( 5, 3, 4, 4, 8, 7, 7)$ \newline
\noindent $||\mathcal{S}_{\beta}||^2=\frac{3}{2}\thickapprox 1.50$

\section{Rank four}

\

\begin{center}
%
 \right]}}$$

\noindent  General solution to $Ux=[1]$: $x= \frac{1 }{ 7}(2,1,2 )^T$. It is an Einstein nilradical.

\

\noindent  Maximal Torus: $D_1=\diag(1, 0, 0, 0, 1, 1, 0 )$,  $D_2=\diag(0, 1, 0, 0, 1, 0, 0 )$, $D_3=\diag( 0, 0, 1, 0, 0, 1, 1)$, $D_4=\diag( 0, 0, 0, 1, 0, 0, 1)$. \newline
Pre-Einstein Derivation: $\frac{ 1}{ 7}\diag ( 4, 5, 4, 5, 9, 8, 9)$ \newline
\noindent $||\mathcal{S}_{\beta}||^2=\frac{7}{5}= 1.4$

\

\begin{center}
%
 \right]}}$$

\noindent  General solution to $Ux=[1]$: $x= \frac{1 }{5 }(1,1,1 )^T$. It is an Einstein nilradical.

\

\noindent  Maximal Torus: $D_1=\diag( 1, 0, 0, 0, 1, 1, 1)$,  $D_2=\diag( 0, 1, 0, 0, 1, 0, 0)$, $D_3=\diag( 0, 0, 1, 0, 0, 1, 0)$, $D_4=\diag( 0, 0, 0, 1, 0, 0, 1)$. \newline
Pre-Einstein Derivation: $\frac{ 2}{ 5}\diag ( 1, 2, 2, 2, 3, 3, 3)$ \newline
\noindent $||\mathcal{S}_{\beta}||^2=\frac{5}{3}\thickapprox 1.67$

\

\begin{center}
%
 \right]}}$$

\noindent  General solution to $Ux=[1]$: $x= \frac{ 1}{7 }(2,1,2)^T$. It is an Einstein nilradical.

\

\noindent  Maximal Torus: $D_1=\diag( 1, 0, 0, 0, 1, 1, 1)$,  $D_2=\diag(0, 1, 0, 0, 1, 1, 1 )$, $D_3=\diag(0, 0, 1, 0, -1, 0, -1 )$, $D_4=\diag( 0, 0, 0, 1, 0, 0, 1)$. \newline
Pre-Einstein Derivation: $\frac{ 1}{7 }\diag ( 5, 5, 6, 5, 4, 10, 9)$ \newline
\noindent $||\mathcal{S}_{\beta}||^2=\frac{7}{5}=1.4$

\

\begin{center}
%
 \right]}}$$

\noindent  General solution to $Ux=[1]$: $x= \frac{1 }{5 }(1,1,1)^T$. It is an Einstein nilradical.

\

\noindent  Maximal Torus: $D_1=\diag( 1, 0, 0, 0, 1, 1, 1)$,  $D_2=\diag( 0, 1, 0, 0, -1, 0, 0)$, $D_3=\diag( 0, 0, 1, 0, 0, -1, 0)$, $D_4=\diag(0, 0, 0, 1, 1, 1, 1 )$. \newline
Pre-Einstein Derivation: $\frac{4 }{5 }\diag ( 1, 1, 1, 1, 1, 1, 2)$ \newline
\noindent $||\mathcal{S}_{\beta}||^2=\frac{5}{3}\thickapprox 1.67$



\end{document}